\documentclass[twoside,11pt]{article}
\usepackage{a4wide,cite,latexsym,amsfonts,amssymb,exscale,epsfig,hyperref,url}
\usepackage[centertags,sumlimits,intlimits,namelimits,reqno]{amsmath}

\usepackage{amsthm}
\theoremstyle{definition}
\newtheorem{theorem}{Theorem}[section]
\newtheorem{lemma}[theorem]{Lemma}
\newtheorem{corollary}[theorem]{Corollary}
\newtheorem{proposition}[theorem]{Proposition}
\newtheorem{definition}[theorem]{Definition}
\newtheorem{remark}[theorem]{Remark}
\newtheorem{example}[theorem]{Example}

\usepackage{bbm}
\def\C{{\mathbbm C}}
\def\N{{\mathbbm N}}

\def\Z{{\mathbbm Z}}

\def\ssl{{\mathfrak{sl}}}

\def\u{{\mathfrak{u}}}

\def\ie{{\sl i.e.\/}}

\def\cf{{\sl c.f.\/}}

\let\phi=\varphi
\let\theta=\vartheta
\let\epsilon=\varepsilon

\def\tr{\mathop{\rm tr}\nolimits}

\def\ev{\mathop{\rm ev}\nolimits}
\def\coev{\mathop{\rm coev}\nolimits}

\def\dim{\mathop{\rm dim}\nolimits}

\def\id{\mathop{\rm id}\nolimits}

\def\End{\mathop{\rm End}\nolimits}

\def\Hom{\mathop{\rm Hom}\nolimits}

\def\im{\mathop{\rm im}\nolimits}
\def\coim{\mathop{\rm coim}\limits}

\def\lim{\mathop{\rm lim}\limits}

\let\hat=\widehat
\let\tilde=\widetilde

\def\pprime{{\prime\prime}}
\def\ppprime{{\prime\prime\prime}}
\def\hotimes{{\hat\otimes}}
\def\ract#1#2{#1\triangleleft #2}
\def\q#1{{[#1]}_q}

\numberwithin{equation}{section}

\makeatletter

\newfont{\@aidxte}{cmsy10}
\newfont{\@aidxel}{cmsy10 scaled 1095}
\newfont{\@aidxtw}{cmsy10 scaled 1200}
\newlength\@aidxtexvi
\newlength\@aidxtexvii
\newlength\@aidxelxvi
\newlength\@aidxelxvii
\newlength\@aidxtwxvi
\newlength\@aidxtwxvii
\newcommand{\alignidx}[1]{%
\@aidxtexvi=\fontdimen16\@aidxte
\@aidxtexvii=\fontdimen17\@aidxte
\@aidxelxvi=\fontdimen16\@aidxel
\@aidxelxvii=\fontdimen17\@aidxel
\@aidxtwxvi=\fontdimen16\@aidxtw
\@aidxtwxvii=\fontdimen17\@aidxtw
{\mbox{$%
\fontdimen16\@aidxte=2.9pt
\fontdimen17\@aidxte=2.9pt
\fontdimen16\@aidxel=3.1pt
\fontdimen17\@aidxel=3.1pt
\fontdimen16\@aidxtw=3.3pt
\fontdimen17\@aidxtw=3.3pt
#1$}}%
\fontdimen16\@aidxte=\@aidxtexvi
\fontdimen17\@aidxte=\@aidxtexvii
\fontdimen16\@aidxel=\@aidxelxvi
\fontdimen17\@aidxel=\@aidxelxvii
\fontdimen16\@aidxtw=\@aidxtwxvi
\fontdimen17\@aidxtw=\@aidxtwxvii}

\makeatother

\newenvironment{myenumerate}{%
\begin{enumerate}
\setlength{\partopsep}{0pt}
\setlength{\parskip}{0pt}}{\end{enumerate}}

\def\nn{\notag}

\def\emph#1{{\sl #1\/}}
\def\sym#1{{\mathcal #1}}
\def\one{\mathbbm{1}}%
\def\bar#1{\overline{#1}}%

\def\op{\mathrm{op}}

\def\Vect{\mathbf{Vect}}
\def\fdVect{\mathbf{fdVect}}

\def\coend{\mathbf{coend}}

\def\msc#1{\noindent{\small Mathematics Subject Classification (2000):
#1\par}}%
\def\keywords#1{\noindent {\small keywords: #1\par}}%

\usepackage[ps,matrix,arrow,frame,curve]{xy}
\everyxy={0;<1mm,0mm>:}

\begin{document}

\title{
Finitely semisimple spherical categories and\\
modular categories are self-dual}
\author{Hendryk Pfeiffer\thanks{E-mail: \texttt{pfeiffer@math.ubc.ca}}}
\date{\small{Department of Mathematics, The University of British Columbia,\\
1984 Mathematics Road, Vancouver, BC, V2T 1Z2, Canada}\\[1ex]
February 3, 2009}

\maketitle

\begin{abstract}

We show that every essentially small finitely semisimple $k$-linear additive
spherical category for which $k=\End(\one)$ is a field, is equivalent to its
dual over the long canonical forgetful functor. This includes the special case
of modular categories. In order to prove this result, we show that the
universal coend of the spherical category, with respect to the long forgetful
functor, is self-dual as a Weak Hopf Algebra.

\end{abstract}

\msc{%
16W30,
18D10
}

\section{Introduction}

\subsection{Spherical categories}

A spherical category is a monoidal category in which each object $X$ has a
specified left-dual $X^\ast$, in which there are canonical isomorphisms $X\to
{X^\ast}^\ast$ that can be used to equip each object with a right-dual as
well, and in which the two ways of forming the trace of a morphism
coincide. There is a coherence theorem for spherical categories generalizing
MacLane's coherence theorem for monoidal categories: it states that the
morphisms of a spherical category can be represented by string diagrams in the
sphere $S^2$ up to isotopy, hence the name ``spherical''. In particular, every
ribbon category is spherical.

The notion of a spherical category was invented by Barrett and
Westbury~\cite{BaWe99} when they developed their version~\cite{BaWe96} of the
Turaev--Viro invariant~\cite{TuVi92}. Finitely semisimple spherical categories
form the most general categories for which this invariant can be still be
formulated. Neither a braiding, a ribbon structure nor a non-degenerate
$S$-matrix are required.

\subsection{Duals of monoidal categories}

In the abstract, we claim that certain finitely semisimple spherical
categories are equivalent to their duals over some functor. What is the dual
of a monoidal category?

Assume for the moment that the monoidal category is the category $\sym{M}_H$
of finite-dim\-en\-sional modules of a Hopf algebra $H$ over some field
$k$. Then the algebra structure of $H$ determines the category underlying
$\sym{M}_H$, while the coalgebra structure of $H$ is responsible for the
monoidal structure of $\sym{M}_H$. Under certain conditions, the coalgebra
structure can be reconstructed from $\sym{M}_H$. The coalgebra structure of
$H$ also determines the category $\sym{M}^H$ of finite-dimensional comodules
of $H$. This raises a question: can we characterize the category $\sym{M}^H$
in terms of the monoidal structure of the category $\sym{M}_H$? For Hopf
algebras, Majid's notion of the dual of a monoidal category over a strong
monoidal functor~\cite{Ma91,Ma92} answers this question: under suitable
conditions, the dual of $\sym{M}_H$ over the forgetful functor
$\sym{M}_H\to\Vect_k$ is equivalent to $\sym{M}^H$. In particular, if $H$ is
finite-dimensional, then the dual of $\sym{M}_H$ over the forgetful functor
$\sym{M}_H\to\Vect_k$ is precisely $\sym{M}^H\simeq\sym{M}_{\hat H}$, the
category of finite-dimensional modules of the dual Hopf algebra $\hat H$. This
justifies the term `dual of a monoidal category'.

The following special case of this notion of the dual of a monoidal category
is much more widely known: the dual of a monoidal category $\sym{C}$ over the
identity functor on itself is the double of $\sym{C}$~\cite{Ma91,JoSt91}, also
called the categorical center.

\subsection{Tannaka--Kre\v\i n reconstruction}

The most interesting finitely semisimple spherical categories, however, do not
admit any strong monoidal functor to $\Vect_k$, and they are therefore not the
categories of modules of any Hopf algebra. Nevertheless, we can show that the
relevant spherical categories are equivalent to the categories of
comodules\footnote{This argument uses a generalization of Tannaka--Kre\v\i n
reconstruction, and one usually reconstructs the universal coacting
coalgebra.} of Weak Hopf Algebras (WHAs). Such a WHA is obtained as a
universal coend using a generalization of Tannaka--Kre\v\i n reconstruction:

\begin{theorem}
\label{thm_intro1}
Let $\sym{C}$ be an essentially small finitely semisimple $k$-linear additive
spherical category for which $k=\End(\one)$ is a field, and let
$\omega\colon\sym{C}\to\Vect_k$ be the long forgetful functor. Then
$H=\coend(\sym{C},\omega)$ is a finite-dimensional split cosemisimple
cospherical WHA for which $H_t\cap H_s\cong k$. Furthermore, $\sym{C}$ is
equivalent as a $k$-linear additive spherical category to the category
$\sym{M}^H$ of finite-dimensional right $H$-comodules.
\end{theorem}

The special case of modular categories was addressed in the
article~\cite{Pf07}. The functor $\omega$ is the long version of Hayashi's
canonical forgetful functor~\cite{Ha99}. We proceed by showing that the
reconstructed WHA is self-dual, and then generalizing Majid's work on duals of
monoidal categories to the case of WHAs. Thus equipped, we can finally prove
our claim from the abstract.

\subsection{Double triangle algebras and Depth-$2$ Frobenius extensions}

Why can we expect the reconstructed WHA to be self-dual? Ocneanu's idea of a
`double triangle algebra' gives an early idea of how to reconstruct some
algebra-like structure from a semisimple monoidal category, see, for example,
the diagrams in~\cite{BoSz97}. In hindsight, the precise definition that
captures the idea of a double triangle algebra is the notion of a Weak Hopf
Algebra~\cite{BoNi99,BoSz00}. If we replace white vertices with black vertices
in the double triangle diagrams of, say, \cite{BoSz97}, we see immediately
that we ought to get a self-dual WHA. In particular, the very first examples
of WHAs~\cite{Bo97} whose categories of modules are modular categories with
objects of non-integer Frobenius--Perron dimension, have been constructed by
this method and are indeed self-dual. It remains to generalize the argument to
work for the WHAs reconstructed from our spherical categories.

In fact, B{\"o}hm and Szlach{\'a}nyi~\cite{BoSz04} have implemented the
duality relation that leads to the self-duality of Ocneanu's double triangle
algebra in the context of abstract depth-$2$ Frobenius extensions. For every
such extension, they construct a pair of dual Hopf algebroids.

It turns out that the WHA that arises as our universal coend
$H=\coend(\sym{C},\omega)$ coincides with one of these Hopf algebroids. Since
the base of these Hopf algebroids is the endomorphism algebra $R=\End(\tilde
V)$ of the universal object
\begin{equation}
\tilde V=\bigoplus_{j\in I} V_j,
\end{equation}
the direct sum of one representative $V_j$ for each isomorphism class of
simple objects of $\sym{C}$, and since $R$ is finite-dimensional, commutative
and separable, both Hopf algebroids of~\cite{BoSz04} are in fact
WHAs. Furthermore, since $\sym{C}$ is semisimple as a pivotal category
(Definition~\ref{def_pivotalss}), one can swap $\tilde V$ and its dual
${\tilde V}^\ast$ and show that these dually paired WHAs are isomorphic.

\begin{theorem}
\label{thm_intro2}
The WHA $H$ of Theorem~\ref{thm_intro1} is self-dual as a pivotal WHA.
\end{theorem}

\begin{corollary}
\label{cor_intro3}
Every essentially small finitely semisimple $k$-linear additive spherical
category for which $k=\End(\one)$ is a field, is equivalent as a $k$-linear
additive spherical category to its dual over the long forgetful functor
$\omega\colon\sym{C}\to\Vect_k$.
\end{corollary}

We can view the isomorphism $H\to\hat H$ used to find a pair of dual bases of
$H$ and its dual $\hat H$, as a generalized Fourier transform. For modular
categories, it is known in a somewhat different context~\cite{LyMa94} that the
$S$-matrix provides a Fourier transform. In our case, however, the Fourier
transform is implemented by a generalized $6j$-symbol. It is therefore no
surprise that our result is already available for spherical categories and
even in the absence of any braiding.

\subsection{TQFTs and state sum invariants}

We are interested in spherical categories because they form the most general
categories for which the Turaev--Viro invariant is still available. Not only
do modular categories form a subclass of finitely semisimple spherical
categories, but also the doubles of a large class of finitely semisimple spherical
categories are modular~\cite{Mu03b}, and so the finitely semisimple spherical
categories can be viewed as more basic than the modular ones.

Since their doubles are modular, finitely semisimple spherical categories form
the most plausible candidate for the structure underlying some as yet unknown
3-dimensional extended Topological Quantum Field Theories (TQFTs) that would
specialize to the familiar TQFTs based on modular categories when the
manifolds have no corners. This construction would categorify the TQFTs
of~\cite{LaPf07}, based on the idea that the double of a monoidal category
categorifies the notion of the center of an algebra.

We are interested in the self-duality of our spherical categories because of
the following observation. While the Turaev--Viro invariant can be constructed
using a finitely semisimple spherical category, Kuperberg's
invariant~\cite{Ku91} is based on a certain involutory Hopf algebra. Barrett
and Westbury have shown~\cite{BaWe95} that Kuperberg's invariant for some Hopf
algebra agrees with the Turaev--Viro invariant for its category of modules
(without taking any quotient of this category modulo `negligible
morphisms'!). This observation raises a number of questions.

First, all known examples of finitely semisimple spherical categories for
which the Turaev--Viro invariant is interesting (\ie, stronger than an
invariant of homotopy type), are not the categories of modules of any Hopf
algebra. We know, however, from Theorem~\ref{thm_intro1} that they are
nevertheless the categories of comodules of some WHAs. The obvious question
is: can Kuperberg's invariant be generalized to WHAs so as to capture all
interesting examples of the Turaev--Viro invariant? Then, in order to pass
from Kuperberg to Turaev--Viro, we would take the category of comodules (no
quotient!), and in order to pass from Turaev--Viro back to Kuperberg, we
would Tannaka--Kre\v\i n reconstruct.

Second, Kuperberg's invariant for some involutory Hopf algebra $H$ and some
Heegaard splitting agrees by construction with the invariant for the dual Hopf
algebra $\hat H$ and the Poincar\'e dual Heegaard splitting. Since the
Poincar\'e dual Heegaard splitting characterizes the same $3$-manifold, the
Kuperberg invariants for $H$ and $\hat H$ always agree for any given Heegaard
splitting. How can this be explained? Once Kuperberg's invariant has been
generalized to WHAs, Theorem~\ref{thm_intro2} will provide the answer.

Third, just as the Kuperberg invariant for some Hopf algebra $H$ and for a
given Heegaard splitting agrees with the invariant for the dual Hopf algebra
$\hat H$ and the Poincar\'e dual Heegaard splitting, we can ask the analogous
question about the Turaev--Viro invariant. We would need the generalization of
the Turaev--Viro invariant to cellular complexes~\cite{GiOe02} and study the
following problem. Given the Turaev--Viro invariant for some finitely
semisimple spherical category and some cellular complex, which category do we
need in order to make the Turaev--Viro invariant for the dual cellular complex
agree with the invariant for the original complex? Corollary~\ref{cor_intro3}
will provide the answer.

Using the long forgetful functor, we can therefore improve a number of
results by M{\"u}ger on the relationship between the invariants of Kuperberg
and Turaev--Viro for the original and for the dual Hopf algebra. We refer in
particular to~\cite[Section 7]{Mu03a}.

Crane and Frenkel~\cite{CrFr94} originally proposed to categorify Kuperberg's
$3$-manifold invariant in order to get access to combinatorial $4$-manifold
invariants. Settling the open questions about spherical categories and their
reconstructed WHAs is a crucial step in the process of deciding which
invariant it is that we want to categorify.

\subsection{Overview}

The article is organized as follows. In Section~\ref{sect_prelim}, we
summarize the definitions and basic properties of WHAs and of the forgetful
functors of their categories of comodules. These are functors with a separable
Frobenius structure. In Section~\ref{sect_reconstruct}, we reconstruct WHAs
from our spherical categories and prove Theorem~\ref{thm_intro1}. We show the
self-duality of the reconstructed WHAs (Theorem~\ref{thm_intro2}) in
Section~\ref{sect_selfdual}. In Section~\ref{sect_repmon}, we introduce the
notion of a dual of a monoidal category over a functor with separable
Frobenius structure and show that the category of modules of the reconstructed
WHA is equivalent as a spherical category to the dual of the original
spherical category over the long forgetful functor. We then combine all these
ingredients and show that our spherical categories are self-dual
(Corollary~\ref{cor_intro3}). Finally, in Section~\ref{sect_example}, we
briefly sketch all these constructions for the special case of the modular
category associated with the quantum group $U_q(\mathfrak{sl}_2)$, $q$ a root
of unity, using the familiar diagrams. The reader may wish to take a quick
look at this example before reading the other sections. In the appendix, we
have collected the basic definitions and results on monoidal categories with
duals and on abelian and semisimple categories.

\section{Preliminaries}
\label{sect_prelim}

\subsection{Functors with Frobenius structure}

We use the following notation. If $\sym{C}$ is a category, we write
$X\in|\sym{C}|$ for the objects $X$ of $\sym{C}$, $\Hom(X,Y)$ for the
collection of all morphisms $f\colon X\to Y$ and $\End(X)=\Hom(X,X)$. We
ignore all set theoretic issues and tacitly assume that the $\Hom(X,Y)$ are
all sets. We denote the identity morphism of $X$ by $\id_X\colon X\to X$ and
the composition of morphisms $f\colon X\to Y$ and $g\colon Y\to Z$ by $g\circ
f\colon X\to Z$. If two objects $X,Y\in|\sym{C}|$ are isomorphic, we write
$X\cong Y$. If two categories $\sym{C}$ and $\sym{D}$ are equivalent, we write
$\sym{C}\simeq\sym{D}$. The identity functor on $\sym{C}$ is denoted by
$1_{\sym{C}}$, and $\sym{C}^\op$ is the opposite category of $\sym{C}$. The
category of vector spaces over a field $k$ is denoted by $\Vect_k$ and its
full subcategory of finite-dimensional vector spaces by $\fdVect_k$. Both are
$k$-linear abelian and symmetric monoidal, while $\fdVect_k$ is in addition
spherical. Appendix~\ref{app_duals} gives background on monoidal categories
with duals.

The forgetful functor of the category of finite-dimensional comodules of a WHA
is not necessarily strong monoidal, but it satisfies the following more
general conditions of a functor with separable Frobenius structure as defined
by Szlach{\'a}nyi~\cite{Sz05}.

\begin{definition}
\label{def_functorfrob}
Let $\sym{C}$ and $\sym{C}^\prime$ be monoidal categories. A \emph{functor
with Frobenius structure}
$(F,F_{X,Y},F_0,F^{X,Y},F^0)\colon\sym{C}\to\sym{C}^\prime$ is a functor
$F\colon\sym{C}\to\sym{C}^\prime$ that is lax monoidal as $(F,F_{X,Y},F_0)$
and oplax monoidal as $(F,F^{X,Y},F^0)$ (see Definition~\ref{def_lax}) and
that satisfies the following compatibility conditions,
\begin{equation}
\label{eq_frobenius1}
\begin{aligned}
\xymatrix{
F(X\otimes Y)\otimes^\prime FZ\ar[rr]^{F_{X\otimes Y,Z}}\ar[dd]_{F^{X,Y}\otimes^\prime\id_{FZ}}&&
F((X\otimes Y)\otimes Z)\ar[rr]^{F\alpha_{X,Y,Z}}&&
F(X\otimes(Y\otimes Z))\ar[dd]^{F^{X,Y\otimes Z}}\\
\\
(FX\otimes^\prime FY)\otimes^\prime FZ\ar[rr]_{\alpha^\prime_{FX,FY,FZ}}&&
FX\otimes^\prime(FY\otimes^\prime FZ)\ar[rr]_{\id_{FX}\otimes^\prime F_{Y,Z}}&&
FX\otimes^\prime F(Y\otimes Z),
}
\end{aligned}
\end{equation}
and
\begin{equation}
\label{eq_frobenius2}
\begin{aligned}
\xymatrix{
FX\otimes^\prime F(Y\otimes Z)\ar[rr]^{F_{X,Y\otimes Z}}\ar[dd]_{\id_{FX}\otimes^\prime F^{Y,Z}}&&
F(X\otimes (Y\otimes Z))\ar[rr]^{F\alpha^{-1}_{X,Y,Z}}&&
F((X\otimes Y)\otimes Z)\ar[dd]^{F^{X\otimes Y,Z}}\\
\\
FX\otimes^\prime(FY\otimes^\prime FZ)\ar[rr]_{{\alpha^\prime}^{-1}_{FX,FY,FZ}}&&
(FX\otimes^\prime FY)\otimes^\prime FZ\ar[rr]_{F_{X,Y}\otimes^\prime\id_{FZ}}&&
F(X\otimes Y)\otimes^\prime FZ,
}
\end{aligned}
\end{equation}
for all $X,Y,Z\in|\sym{C}|$. It is called a \emph{functor with separable
Frobenius structure} if in addition
\begin{equation}
F_{X,Y}\circ F^{X,Y} = \id_{F(X\otimes Y)},
\end{equation}
for all $X,Y\in|\sym{C}|$.
\end{definition}

Note that every strong monoidal functor between monoidal categories is a
functor with separable Frobenius structure, and that in this case, the
conditions~\eqref{eq_frobenius1} and~\eqref{eq_frobenius2} both follow from
the hexagon axiom.

\subsection{Weak Hopf Algebras and their comodules}

In this section, we summarize the relevant definitions and properties of Weak
Bialgebras (WBAs) and Weak Hopf Algebras (WHAs)
following~\cite{Ni98,BoNi99,BoSz00} and of their categories of comodules
following~\cite{Pf07}.

\begin{definition}
\label{def_wba}
A \emph{Weak Bialgebra} $(H,\mu,\eta,\Delta,\epsilon)$ over a field
$k$ is a $k$-vector space $H$ with linear maps $\mu\colon H\otimes
H\to H$ (\emph{multiplication}), $\eta\colon k\to H$ (\emph{unit}),
$\Delta\colon H\to H\otimes H$ (\emph{comultiplication}), and
$\epsilon\colon H\to k$ (\emph{counit}) such that the following
conditions hold:
\begin{myenumerate}
\item
$(H,\mu,\eta)$ is an associative unital algebra, \ie\
$\mu\circ(\mu\otimes\id_H)=\mu\circ(\id_H\otimes\mu)$ and
$\mu\circ(\eta\otimes\id_H)=\id_H=\mu\circ(\id_H\otimes\eta)$.
\item
$(H,\Delta,\epsilon)$ is a coassociative counital coalgebra, \ie\
$(\Delta\otimes\id_H)\circ\Delta=(\id_H\otimes\Delta)\circ\Delta$
and
$(\epsilon\otimes\id_H)\circ\Delta=\id_H=(\id_H\otimes\epsilon)\circ\Delta$.
\item
The following compatibility conditions hold:
\begin{eqnarray}
\label{eq_wba1}
\Delta\circ\mu
&=& (\mu\otimes\mu)\circ(\id_H\otimes\sigma_{H,H}\otimes\id_H)\circ(\Delta\otimes\Delta),\\
\label{eq_wba2}
\epsilon\circ\mu\circ(\mu\otimes\id_H)
&=& (\epsilon\otimes\epsilon)\circ(\mu\otimes\mu)\circ(\id_H\otimes\Delta\otimes\id_H)\nn\\
&=& (\epsilon\otimes\epsilon)\circ(\mu\otimes\mu)\circ(\id_H\otimes\Delta^\op\otimes\id_H),\\
\label{eq_wba3}
(\Delta\otimes\id_H)\circ\Delta\circ\eta
&=& (\id_H\otimes\mu\otimes\id_H)\circ(\Delta\otimes\Delta)\circ(\eta\otimes\eta)\nn\\
&=& (\id_H\otimes\mu^\op\otimes\id_H)\circ(\Delta\otimes\Delta)\circ(\eta\otimes\eta).
\end{eqnarray}
\end{myenumerate}
Here $\sigma_{V,W}\colon V\otimes W\to W\otimes V$, $v\otimes w\mapsto
w\otimes v$ is the transposition of the tensor factors, \ie\ the symmetric
braiding of $\Vect_k$, and by $\Delta^\op=\sigma_{H,H}\circ\Delta$ and
$\mu^\op=\mu\circ\sigma_{H,H}$, we denote the \emph{opposite comultiplication}
and \emph{opposite multiplication}, respectively. We tacitly identify the
vector spaces $(V\otimes W)\otimes U\cong V\otimes(W\otimes U)$ and $V\otimes
k\cong V\cong k\otimes V$, exploiting the coherence theorem for the monoidal
category $\Vect_k$.
\end{definition}

We use the term \emph{comultiplication} for the operation $\Delta$ in a
coalgebra, whereas \emph{coproduct} always refers to a colimit in a category.

\begin{definition}
A \emph{homomorphism} $\phi\colon H\to H^\prime$ of WBAs
$(H,\mu,\eta,\Delta,\epsilon)$ and
$(H^\prime,\mu^\prime,\Delta^\prime,\epsilon^\prime)$ over the same
field $k$ is a $k$-linear map that is a homomorphism of unital
algebras, \ie\ $\phi\circ\eta=\eta^\prime$ and
$\phi\circ\mu=\mu^\prime\circ(\phi\otimes\phi)$, as well as a
homomorphism of counital coalgebras, \ie\
$\epsilon^\prime\circ\phi=\epsilon$ and
$\Delta^\prime\circ\phi=(\phi\otimes\phi)\circ\Delta$.
\end{definition}

\begin{definition}
\label{def_sourcetarget}
Let $(H,\mu,\eta,\Delta,\epsilon)$ be a WBA. The linear maps
$\epsilon_t\colon H\to H$ (\emph{target counital map}) and
$\epsilon_s\colon H\to H$ (\emph{source counital map}) are defined by
\begin{eqnarray}
\label{eq_epsilont}
\epsilon_t&:=&(\epsilon\otimes\id_H)\circ(\mu\otimes\id_H)\circ(\id_H\otimes\sigma_{H,H})
\circ(\Delta\otimes\id_H)\circ(\eta\otimes\id_H),\\
\label{eq_epsilons}
\epsilon_s&:=&(\id_H\otimes\epsilon)\circ(\id_H\otimes\mu)\circ(\sigma_{H,H}\otimes\id_H)
\circ(\id_H\otimes\Delta)\circ(\id_H\otimes\eta).
\end{eqnarray}
\end{definition}

Both $\epsilon_t$ and $\epsilon_s$ are idempotents. A WBA
$(H,\mu,\eta,\Delta,\epsilon)$ is a bialgebra if and only if
$\Delta\circ\eta=\eta\otimes\eta$, if and only if
$\epsilon\circ\mu=\epsilon\otimes\epsilon$, if and only if
$\epsilon_s=\eta\circ\epsilon$ and if and only if
$\epsilon_t=\eta\circ\epsilon$.

\begin{proposition}
\label{prop_hsht}
Let $(H,\mu,\eta,\Delta,\epsilon)$ be a WBA.
\begin{myenumerate}
\item
The subspace $H_t:=\epsilon_t(H)$ (\emph{target base algebra}) forms
a unital subalgebra and a left coideal, \ie\
\begin{equation}
\Delta(H_t)\subseteq H\otimes H_t.
\end{equation}
\item
The subspace $H_s:=\epsilon_s(H)$ (\emph{source base algebra}) forms
a unital subalgebra and a right coideal, \ie\
\begin{equation}
\label{eq_rightcoideal}
\Delta(H_s)\subseteq H_s\otimes H.
\end{equation}
\item
The subalgebras $H_s$ and $H_s$ commute, \ie\ $xy=yx$ for all $x\in
H_t$ and $y\in H_s$.
\end{myenumerate}
\end{proposition}

\begin{definition}
A \emph{Weak Hopf Algebra} $(H,\mu,\eta,\Delta,\epsilon,S)$ is a Weak
Bialgebra $(H,\mu,\eta,\Delta,\epsilon)$ with a linear map $S\colon
H\to H$ (\emph{antipode}) that satisfies the following conditions:
\begin{eqnarray}
\label{eq_wha1}
\mu\circ(\id_H\otimes S)\circ\Delta &=& \epsilon_t,\\
\label{eq_wha2}
\mu\circ(S\otimes\id_H)\circ\Delta &=& \epsilon_s,\\
\label{eq_wha3}
\mu\circ(\mu\otimes\id_H)\circ(S\otimes\id_H\otimes S)
\circ(\Delta\otimes\id_H)\circ\Delta&=&S.
\end{eqnarray}
\end{definition}

For convenience, we write $1=\eta(1)$ and omit parentheses in products,
exploiting associativity. We also use Sweedler's notation and write
$\Delta(x)=x^\prime\otimes x^\pprime$ for the comultiplication of $x\in H$ as
an abbreviation of the expression $\Delta(x)=\sum_k a_k\otimes b_k$ with some
$a_k,b_k\in H$. Similarly, we write
$((\Delta\otimes\id_H)\circ\Delta)(x)=x^\prime\otimes x^\pprime\otimes
x^\ppprime$, exploiting coassociativity.

\begin{definition}
A \emph{homomorphism} $\phi\colon H\to H^\prime$ of WHAs is a
homomorphism of WBAs for which $\phi\circ S=S^\prime\circ\phi$.
\end{definition}

The antipode of a WHA is an algebra antihomomorphism, \ie\
$S\circ\mu=\mu^\op\circ(S\otimes S)$ and $S\circ\eta=\eta$, as well as a
coalgebra antihomomorphism, \ie\ $(S\otimes S)\circ\Delta=\Delta^\op\circ S$
and $\epsilon\circ S=\epsilon$.

We extend Sweedler's notation to the right $H$-comodules and write
$\beta(v)=v_V\otimes v_H$ for the coaction $\beta\colon V\to V\otimes H$ of
$H$ on some vector space $V$.

\begin{proposition}[see~\cite{Ni98,Pf07}]
\label{prop_monoidal}
Let $H$ be a WBA. Then the category $\sym{M}^H$ of finite-dimensional right
$H$-comodules is a monoidal category
$(\sym{M}^H,\hotimes,H_s,\alpha,\lambda,\rho)$. Here the monoidal unit
object is the source base algebra $H_s$ with the coaction
\begin{equation}
\beta_{H_s}\colon H_s\to H_s\otimes H,\qquad
x\mapsto x^\prime\otimes x^\pprime.
\end{equation}
The tensor product $V\hotimes W$ of two right $H$-comodules is the
\emph{truncated tensor product}, which is the vector space
\begin{equation}
V\hotimes W := \{\,v\otimes w\in V\otimes W\mid\quad v\otimes w=
P_{V,W}(v\otimes w)\,\}
\end{equation}
where $P_{V,W}$ is the $k$-linear idempotent
\begin{equation}
\label{eq_idempotent}
P_{V,W}\colon V\otimes W\to V\otimes W,\quad
v\otimes w\mapsto (v_V\otimes w_W)\epsilon(v_Hw_H).
\end{equation}
The coaction on $V\hotimes W$ is given by
\begin{equation}
\beta_{V\hotimes W}\colon V\hotimes W\to (V\hotimes W)\otimes H,\quad
v\otimes w \mapsto (v_V\otimes w_W)\otimes (v_Hw_H).
\end{equation}
The unit constraints of the monoidal category are
\begin{alignat}{3}
\lambda_V &\colon H_s\hotimes V\to V,&&\quad x\otimes v\mapsto v_V\epsilon(xv_H),\\
\rho_V    &\colon V\hotimes H_s\to V,&&\quad v\otimes x\mapsto v_V\epsilon(v_H\epsilon_s(x)),
\end{alignat}
and the associator is induced from that of $\Vect_k$.
\end{proposition}

\begin{proposition}[see {\cite[Proposition~5.9]{Pf07}}]
\label{prop_forgetful}
Let $(H,\mu,\eta,\Delta,\epsilon)$ be a WBA and $U\colon\sym{M}^H\to\Vect_k$
be the obvious forgetful functor. Then $(U,U_{X,Y},U_0,U^{X,Y},U^0)$ is a
$k$-linear faithful functor with separable Frobenius structure and it takes values
in $\fdVect_k$. The Frobenius structure is given by
\begin{eqnarray}
U_{X,Y}=\coim P_{X,Y}  \colon UX\otimes UY         &\to& P_{X,Y}(UX\otimes UY),\\
U_0    =\eta           \colon k                    &\to& H_s,\\
U^{X,Y}=\im P_{X,Y}    \colon P_{X,Y}(UX\otimes UY)&\to& UX\otimes UY,\\
U^0    =\epsilon|_{H_s}\colon H_s                  &\to& k
\end{eqnarray}
Here $P_{X,Y}$ denotes the idempotent of~\eqref{eq_idempotent} with its image
factorization $P_{X,Y}=\im P_{X,Y}\circ\coim P_{X,Y}$. Its image
$P_{X,Y}(UX\otimes UY)=U(X\hotimes Y)$ is the vector space underlying the
truncated tensor product. Finally, $H_s=U\one$ is the vector space underlying
the monoidal unit.
\end{proposition}

\begin{proposition}
\label{prop_autonomous}
Let $H$ be a WHA. Then $\sym{M}^H$ is left-autonomous if the left-dual of
every object $V\in|\sym{M}^H|$ is chosen to be $(V^\ast,\ev_V,\coev_V)$, where
the dual vector space $V^\ast$ is equipped with the coaction
\begin{equation}
\beta_{V^\ast}\colon V^\ast\to V^\ast\otimes H,\qquad
\theta\mapsto (v\mapsto \theta(v_V)\otimes S(v_H)),
\end{equation}
and the evaluation and coevaluation maps are given by
\begin{alignat}{3}
\ev_V   &\colon V^\ast\hotimes V\to H_s,&&\quad \theta\otimes v\to\theta(v_V)\epsilon_s(v_H),\\
\coev_V &\colon H_s\to V\hotimes V^\ast,&&\quad x\mapsto\sum_j ({(v_j)}_V\otimes \theta^j)\epsilon(x{(v_j)}_H).
\end{alignat}
Here we have used the evaluation and coevaluation maps that turn
$V^\ast$ into a left-dual of $V$ in the category $\Vect_k$:
\begin{alignat}{3}
\ev_V^{(\Vect_k)}   &\colon V^\ast\otimes V\to k,&&\quad \theta\otimes v\mapsto \theta(v),\\
\coev_V^{(\Vect_k)} &\colon k\to V\otimes V^\ast,&&\quad 1\mapsto \sum_j v_j\otimes \theta^j.
\end{alignat}
\end{proposition}

\section{Tannaka--Kre\v\i n reconstruction}
\label{sect_reconstruct}

In this section, we define the notion of a cospherical WHA. We show that the
universal coend, $H=\coend(\sym{C},\omega)$, of our spherical category
$\sym{C}$ with respect to the long forgetful functor
$\omega\colon\sym{C}\to\Vect_k$ is cospherical and that the category $\sym{C}$
is equivalent as a spherical category to the category of finite-dimensional
right $H$-comodules. This generalizes~\cite{Pf07} to our class of spherical
categories. In order to construct $\coend(\sym{C},\omega)$, we need $\sym{C}$
to be small. We assume this from now on.

\subsection{The long forgetful functor}

First, we study the properties of the long forgetful functor
$\omega\colon\sym{C}\to\Vect_k$ associated with our spherical category
$\sym{C}$. It turns out that Section~3 of~\cite{Pf07} applies to the spherical
case without any significant change, and so we keep this section brief.

\begin{definition}
\label{def_longfunctor}
Let $\sym{C}$ be a finitely semisimple $k$-linear additive spherical category
such that $k=\End(\one)$ is a field. The \emph{long forgetful functor} is the
functor
\begin{eqnarray}
\omega\colon\sym{C}\to\Vect_k,\quad X &\mapsto&\Hom(\tilde V,\tilde V\otimes X),\\
f &\mapsto& (\id_{\tilde V}\otimes f)\circ-.\nn
\end{eqnarray}
Here $\tilde V$ denotes the \emph{universal object} of $\sym{C}$,
\begin{equation}
\tilde V=\bigoplus_{j\in I} V_j,
\end{equation}
where the sum ranges over a set of representatives of the equivalence classes of
simple objects $V_j$, $j\in I$, of $\sym{C}$.
\end{definition}

\begin{proposition}
\label{prop_leftdual}
Let $\sym{C}$ be a finitely semisimple $k$-linear additive spherical category,
$k=\End(\one)$ be a field and $\omega\colon\sym{C}\to\Vect_k$ be the long
forgetful functor. Then for each $X\in|\sym{C}|$, $\omega(X)$ has a left-dual
$({\omega(X)}^\ast,\ev_{\omega(X)},\coev_{\omega(X)})$ where
${\omega(X)}^\ast=\Hom(\tilde V\otimes X,\tilde V)$,
\begin{alignat}{2}
\label{eq_hayashiev}
\ev_{\omega(X)}&\colon{\omega(X)}^\ast\otimes\omega(X)\to k,&&\quad
\theta\otimes v\mapsto\tr_{\tilde V}(D_{\tilde V}\circ\theta\circ v),\\
\label{eq_hayashicoev}
\coev_{\omega(X)}&\colon k\to\omega(X)\otimes{\omega(X)}^\ast,&&\quad
1\mapsto \sum_je_j^{(X)}\otimes e^j_{(X)}.
\end{alignat}
By $D\colon 1_\sym{C}\Rightarrow 1_\sym{C}$ we denote the natural equivalence
\begin{equation}
D_X\colon X\to X,\qquad
D_X:=\sum_{\ell=1}^{n^X}\imath_\ell^X\circ\pi_\ell^X{(\dim V_{j_\ell^X})}^{-1}.
\end{equation}
Here, $\imath_\ell^X$, $\pi_\ell^X$, $j_\ell^X$ and $n^X$ are as in
Definition~\ref{def_semisimple}(3), decomposing $X$ in a biproduct (direct
sum) of simple objects
\begin{equation}
X\cong V_{j_1^X}\oplus\cdots V_{j_{n^X}^X}
\end{equation}
where $j_1^X,\ldots,j_{n^X}^X\in I$. In~\eqref{eq_hayashicoev},
${(e^{(X)}_j)}_j$ and ${(e^j_{(X)})}_j$ denote a pair of dual bases of
$\omega(X)$ and ${\omega(X)}^\ast$ with respect to the bilinear
map~\eqref{eq_hayashiev}.
\end{proposition}

Recall that in a finitely semisimple $k$-linear additive spherical category,
the bilinear map~\eqref{eq_hayashiev} is non-degenerate
(Proposition~\ref{prop_nondegenerate}). For a morphism $f\colon X\to Y$ of
$\sym{C}$, the morphism dual to $\omega(f)=(\id_{\tilde V}\otimes f)\circ-$ is
given by
\begin{equation}
{\omega(f)}^\ast = -\circ(\id_{\tilde V}\otimes f).
\end{equation}

\begin{theorem}
Let $\sym{C}$ be a finitely semisimple $k$-linear additive spherical category
such that $k=\End(\one)$ is a field. The long forgetful functor is $k$-linear,
faithful, takes values in $\fdVect_k$ and has a separable Frobenius structure
$(\omega,\omega_{X,Y},\omega_0,\omega^{X,Y},\omega^0)$ with
\begin{alignat}{2}
\label{eq_longfrob1}
\omega_{X,Y}&\colon\omega(X)\otimes\omega(Y)\to\omega(X\otimes Y),&&\quad
f\otimes g\mapsto\alpha_{\tilde V,X,Y}\circ(f\otimes\id_Y)\circ g,\\
\omega_0&\colon k\to\omega(\one),&&\quad
1\mapsto\rho_{\tilde V}^{-1},
\end{alignat}
and
\begin{alignat}{2}
\omega^{X,Y}&\colon\omega(X\otimes Y)\to\omega(X)\otimes\omega(Y),&&\quad
h\mapsto\sum_{j,\ell}\ev_{\omega(X\otimes Y)}(e^{j\ell}_{(X\otimes Y)}\otimes h)\,
e_j^{(X)}\otimes e_{\ell}^{(Y)},\\
\omega^0&\colon\omega(\one)\to k,&&\quad
v\mapsto \ev_{\omega(\one)}(\rho_{\tilde V}\otimes v).
\end{alignat}
where we have written
\begin{equation}
e^{j\ell}_{(X\otimes Y)} := e^\ell_{(Y)}\circ(e^j_{(X)}\otimes\id_Y)\circ\alpha^{-1}_{\tilde V,X,Y}.
\end{equation}
\end{theorem}

\begin{lemma}
\label{lem_isoanti}
Let $\sym{C}$ be a finitely semisimple $k$-linear additive spherical category,
$k=\End(\one)$ be a field and $\omega\colon\sym{C}\to\Vect_k$ be the long
forgetful functor. Then there are natural isomorphisms for all
$X\in|\sym{C}|$,
\begin{eqnarray}
\Phi_X \colon \omega(X)&\to& {\omega(X^\ast)}^\ast,\nn\\
v&\mapsto& D_{\tilde V}^{-1}\circ\rho_{\tilde V}\circ(\id_{\tilde V}\otimes\bar\ev_X)
\circ\alpha_{\tilde V,X,X^\ast}\circ (v\otimes\id_{X^\ast})\circ (D_{\tilde V}\otimes\id_{X^\ast}),\\
\label{eq_psi}
\Psi_X \colon {\omega(X)}^\ast&\to& \omega(X^\ast),\nn\\
\theta&\mapsto& (\theta\otimes\id_{X^\ast})\circ\alpha^{-1}_{\tilde V,X,X^\ast}\circ(\id_{\tilde V}\otimes\coev_X)
\circ\rho^{-1}_{\tilde V}.
\end{eqnarray}
Their composites are given by
\begin{eqnarray}
\Psi_{X^\ast}\circ\Phi_X \colon\omega X&\to&\omega({X^\ast}^\ast),\nn\\
v&\mapsto& (D_{\tilde V}^{-1}\otimes\tau_X)\circ v\circ D_{\tilde V},\\
\Phi_{X^\ast}\circ\Psi_X \colon{\omega(X)}^\ast&\to&{\omega({X^\ast}^\ast)}^\ast,\nn\\
\theta&\mapsto& D_{\tilde V}^{-1}\circ\theta\circ(D_{\tilde V}\otimes\tau_X^{-1}).
\end{eqnarray}
Here $\tau_X\colon X\to{X^\ast}^\ast$ denotes the pivotal structure of
$\sym{C}$ (see~\eqref{eq_pivotal}).
\end{lemma}

The following diagrams illustrate the maps $\Phi_X$ and $\Psi_X$:
\begin{equation}
\Phi_X\biggl(
\begin{xy}
(0,0)*+{\hbox to 3em{\hfill $v$\hfill}}*\frm{-}!U="t" !D="b";
"t"+< 0em,0em>;"t"+< 0em, 2em> **\dir{-}; ?(.5)*\dir{<}; ?(.7)+< 0.7em,0em>*{\tilde V};
"b"+<-1em,0em>;"b"+<-1em,-2em> **\dir{-}; ?(.5)*\dir{>}; ?(.7)+<-0.7em,0em>*{\tilde V};
"b"+< 1em,0em>;"b"+< 1em,-2em> **\dir{-}; ?(.5)*\dir{>}; ?(.7)+< 0.7em,0em>*{X};
\end{xy}
\biggr) :=
\begin{xy}
(0,0)*+{\hbox to 3em{\hfill $v$\hfill}}*\frm{-}!U="t" !D="b";
(0,0)+<0em,3em>*+{D_{\tilde V}}*\frm{o}!U="ct" !D="cb";
(0,0)+<-1em,-3em>*+{D_{\tilde V}^{-1}}*\frm{o}!U="dt" !D="db";
"t" +< 0em,  0em>;"cb"+<0em,-0.2em> **\dir{-}; ?(.5)*\dir{<};
"ct"+< 0em,0.2em>;"ct"+<0em,   2em> **\dir{-}; ?(.5)*\dir{<}; ?(.7)+< 0.7em,0em>*{\tilde V};
"b"+<-1em,0em>;"dt"+<0em,0.2em> **\dir{-}; ?(.5)*\dir{>};
"db"+<0em,-0.2em>;"db"+<0em,-2em> **\dir{-}; ?(.5)*\dir{>}; ?(.7)+<-0.7em,0em>*{\tilde V};
"b"+< 1em,0em>;"t"+< 3em, 2em> **\crv{"b"+<1em,-1em>&"b"+<3em,-1em>&"b"+<3em,0em>};
?(.9)*\dir{>}; ?(.95)+<0.7em,0em>*{X};
\end{xy},\qquad
\Psi_X\biggl(
\begin{xy}
(0,0)*+{\hbox to 3em{\hfill $\theta$\hfill}}*\frm{-}!U="t" !D="b";
"t"+<-1em,0em>;"t"+<-1em, 2em> **\dir{-}; ?(.5)*\dir{<}; ?(.7)+<-0.7em,0em>*{\tilde V};
"t"+< 1em,0em>;"t"+< 1em, 2em> **\dir{-}; ?(.5)*\dir{<}; ?(.7)+< 0.7em,0em>*{X};
"b"+< 0em,0em>;"b"+< 0em,-2em> **\dir{-}; ?(.5)*\dir{>}; ?(.7)+< 0.7em,0em>*{\tilde V};
\end{xy}
\biggr) :=
\begin{xy}
(0,0)*+{\hbox to 3em{\hfill $\theta$\hfill}}*\frm{-}!U="t" !D="b";
"t"+<-1em,0em>;"t"+<-1em, 2em> **\dir{-}; ?(.5)*\dir{<}; ?(.7)+<-0.7em,0em>*{\tilde V};
"t"+< 1em,0em>;"b"+< 3em,-2em> **\crv{"t"+<1em,1em>&"t"+<3em,1em>&"t"+<3em,0em>};
?(.9)*\dir{<}; ?(.95)+<0.7em,0em>*{X};
"b"+< 0em,0em>;"b"+< 0em,-2em> **\dir{-}; ?(.5)*\dir{>}; ?(.7)+< 0.7em,0em>*{\tilde V};
\end{xy}
\end{equation}
These diagrams specify well-defined morphisms of $\sym{C}$ due to the
coherence theorem of~\cite{BaWe99}.

\begin{proposition}
Let $\sym{C}$ be a finitely semisimple $k$-linear additive spherical category,
$k=\End(\one)$ be a field and $\omega\colon\sym{C}\to\Vect_k$ be the long
forgetful functor. Let ${\{e^{(X)}_j\}}_j$ and ${\{e^j_{(X)}\}}_j$ form a pair
of dual bases of $\omega(X)$ and ${\omega(X)}^\ast$ with respect
to~\eqref{eq_hayashiev}. Then
\begin{equation}
\label{eq_complete1}
\sum_j e^{(X)}_j\circ e^j_{(X)} = \id_{\tilde V\otimes X}
\end{equation}
and
\begin{equation}
\label{eq_complete2}
\sum_j \Psi(e^j_{(X)})\circ \Phi(e^{(X)}_j) = \id_{\tilde V\otimes X^\ast}.
\end{equation}
\end{proposition}

\subsection{Cospherical Weak Hopf Algebras}
\label{sect_cospherical}

In this section, we define the notion of a cospherical WHA and show that its
category of finite-dimensional comodules is spherical. The special case of
cospherical Hopf algebras reduces to the definitions given in~\cite{Oe03}. For
background on autonomous, pivotal and spherical categories, we refer to
Appendix~\ref{app_duals}.

\begin{definition}
Let $(H,\mu,\eta,\Delta,\epsilon,S)$ be a WHA. A linear form $f\colon
H\to k$ is called:
\begin{enumerate}
\item
\emph{convolution invertible} if there exists some linear $\bar
f\colon H\to k$ such that $f(x^\prime)\bar
f(x^\pprime)=\epsilon(x)=\bar f(x^\prime)f(x^\pprime)$ for all $x\in H$,
\item
\emph{dual central} if $f(x^\prime)x^\pprime = x^\prime f(x^\pprime)$
for all $x\in H$,
\item
\emph{dual group-like} if it is convolution invertible and if for all
$x,y\in H$,
\begin{equation}
\label{eq_dualgrouplike}
w(x^\prime)w(y^\prime)\epsilon(x^\pprime
y^\pprime)=w(xy)=\epsilon(x^\prime
y^\prime)w(x^\pprime)w(y^\pprime).
\end{equation}
\end{enumerate}
\end{definition}

Note that $\bar f$ in (1) is uniquely determined by $f$. Every dual
group-like linear form $f\colon H\to k$ also satisfies $f(S(x))=\bar
f(x)$ and $w(\epsilon_s(x))=\epsilon(x)=w(\epsilon_t(x))$ for all
$x\in H$.

\begin{definition}
A \emph{copivotal} WHA $(H,\mu,\eta,\Delta,\epsilon,S,w)$ is a WHA
$(H,\mu,\eta,\Delta,\epsilon,S)$ with a dual group-like linear form
$w\colon H\to k$ that satisfies
\begin{equation}
\label{eq_copivotal}
S^2(x) = w(x^\prime)x^\pprime\bar w(x^\ppprime)
\end{equation}
for all $x\in H$.
\end{definition}

\begin{proposition}
\label{prop_pivotal}
Let $(H,\mu,\eta,\Delta,\epsilon,S,w)$ be a copivotal WHA. Then the
category $\sym{M}^H$ is pivotal with $\tau_V\colon V\to {V^\ast}^\ast$
given by
\begin{equation}
\label{eq_pivotalcomodule}
\tau_V(v) = \tau_V^{(\Vect_k)}(v_V)w(v_H)
\end{equation}
for all finite-dimensional right $H$-comodules $V\in|\sym{M}^H|$ and $v\in
V$. Here we denote by $\tau_V^{(\Vect_k)}\colon V\to {V^\ast}^\ast$ the
pivotal structure of $\Vect_k$ which is just the usual canonical
identification $V\cong {V^\ast}^\ast$.
\end{proposition}

\begin{proof}
The category $\sym{M}^H$ is left-autonomous by
Proposition~\ref{prop_autonomous}. The map $\tau_V$ is obviously
$k$-linear. Its inverse is given by
\begin{equation}
\tau_V^{-1}(\tau_V^{(\Vect_k)}(v)) = v_V\bar w(v_H)
\end{equation}
for all $\tau_V^{(\Vect_k)}(v)\in {V^\ast}^\ast$. In order to show that
$\tau_V$ is a morphism of right $H$-comodules, we need~\eqref{eq_copivotal}
and the fact that $H$ coacts on ${V^\ast}^\ast$ by
\begin{equation}
\label{eq_vectpivotal}
{(\tau_V^{(\Vect_k)}(v))}_V\otimes {(\tau_V^{(\Vect_k)}(v))}_H
= \tau_V^{(\Vect_k)}(v_V)\otimes S^2(v_H)
\end{equation}
for all $v\in V$. Naturality follows from the properties of a comodule and
from the naturality of $\tau_V^{(\Vect_k)}$. In order to
verify~\eqref{eq_pivotal}, we compute the dual $\tau_V^\ast$ in $\sym{M}^H$,
using~\eqref{eq_leftdual}, Propositions~\ref{prop_monoidal}
and~\ref{prop_autonomous}, \eqref{eq_copivotal}, \cite[equation (5.13)]{Pf07},
\eqref{eq_vectpivotal} and the fact that $\Vect_k$ is pivotal.
\end{proof}

Recall that in a pivotal category, the left- and right-duals are related, and
we can define traces. In general, however, left- and right-traces
$\tr_V^{(L)}(f)$ and $\tr_V^{(R)}(f)$ of some morphism $f\colon V\to V$ need
not agree (Appendix~\ref{app_duals}).

\begin{definition}
A \emph{cospherical} WHA $H$ is a copivotal WHA for which
\begin{equation}
\label{eq_cospherical}
\tr^{(L)}_V(f) = \tr^{(R)}_V(f)
\end{equation}
for all finite-dimensional right $H$-comodules $V\in|\sym{M}^H|$ and
all morphisms $f\colon V\to V$. Recall that~\eqref{eq_cospherical} is an
identity between morphisms $H_s\to H_s$.
\end{definition}

\noindent
If $H$ is a cospherical WHA, then $\sym{M}^H$ is therefore not only pivotal,
but also spherical.

\begin{example}
Every coribbon WHA~\cite[Definition~4.17]{Pf07} is cospherical.
\end{example}

\begin{proof}
Every coribbon WHA is copivotal because of~\cite[Lemma~5.4]{Pf07}. Its
category of finite-dimensional comodules is
ribbon~\cite[Proposition~5.13]{Pf07}, and so left- and right-traces agree by
Example~\ref{ex_ribboncat}.
\end{proof}

The traces in~\eqref{eq_cospherical} are the left- and right-traces in the
pivotal category $\sym{M}^H$. They can be computed as follows.

\begin{proposition}
Let $(H,\mu,\eta,\Delta,\epsilon,S,w)$ be a copivotal WHA, $V\in|\sym{M}^H|$
and $f\colon V\to V$. Then
\begin{eqnarray}
\label{eq_tracel}
\tr^{(L)}_V(f)[h] &=& \epsilon(hS(c_f^\prime))\epsilon_s(c_f^\pprime)\bar w(c_f^\ppprime)\in H_s,\\
\label{eq_tracer}
\tr^{(R)}_V(f)[h] &=& w(c_f^\prime)\epsilon(hc_f^\pprime)\epsilon_s(S(c_f^\ppprime))\in H_s,
\end{eqnarray}
for all $h\in H_s$. Here $c_f=\sum_{j,\ell=1}^nf_{j\ell}c_{\ell j}^{(V)}$
where $f_{\ell j}\in k$ are the coefficients of $f$, \ie\
$f(e_j)=\sum_{\ell=1}^ne_\ell f_{\ell j}$, and $c^{(V)}_{\ell j}$ are the
matrix elements of the right $H$-comodule $V$, \ie\
$\beta_V(e_j)=\sum_{\ell=1}^ne_\ell\otimes c_{\ell j}^{(V)}$, for some basis
${(e_j)}_{1\leq j\leq n}$ of $V$.
\end{proposition}

\begin{proof}
Compute the traces of Definition~\ref{def_traces} in the pivotal category
$\sym{M}^H$, using the left-autonomous structure of
Proposition~\ref{prop_autonomous} and the pivotal structure of
Proposition~\ref{prop_pivotal}.
\end{proof}

In the following, we study the left- and right-traces for WHAs that satisfy a
number of additional conditions. We will see below that the universal coend of
our spherical categories over the long forgetful functor satisfies these
conditions.

\begin{proposition}
\label{prop_traceink}
Let $(H,\mu,\eta,\Delta,\epsilon,S,w)$ be a copivotal WHA over $k$ with
$H_t\cap H_s\cong k$. Then there are elements $t^{(V)}_L,t^{(V)}_R\in k$ for
each $V\in|\sym{M}^H|$ such that
\begin{eqnarray}
\label{eq_tracelsimple}
\sum_{j,\ell=1}^n \epsilon_s(c_{j\ell}^{(V)})\bar w(c_{\ell j}^{(V)}) &=& t^{(L)}_V\,1,\\
\label{eq_tracersimple}
\sum_{j,\ell=1}^n w(c_{j\ell}^{(V)})\epsilon_s(S(c_{\ell j}^{(V)})) &=& t^{(R)}_V\,1.
\end{eqnarray}
\end{proposition}

\begin{proof}
Since $H_t\cap H_s\cong k$, the monoidal unit object $H_s$ is simple, \ie\
$\End(H_s)\cong k$ in
$\sym{M}^H$~\cite[Lemma~5.16]{Pf07}. Evaluate~\eqref{eq_tracel}
and~\eqref{eq_tracer} for $h=1$ and $f=\id_V$, \ie\
$f_{\ell j}=\delta_{\ell j}$ and $c_f=\sum_{j=1}^nc_{jj}^{(V)}$.
\end{proof}

\begin{theorem}
\label{thm_cospherical}
Let $H$ be a finite-dimensional split cosemisimple copivotal WHA over $k$ with
$H_t\cap H_s\cong k$. $H$ is cospherical if, and only if, $t^{(L)}_V=t^{(R)}_V$
in Proposition~\ref{prop_traceink} for all simple $V\in|\sym{M}^H|$.
\end{theorem}

\begin{proof}
Since $H$ is finite-dimensional split cosemisimple and $H_t\cap H_s\cong k$,
the category $\sym{M}^H$ is finitely semisimple as a pivotal
category~\cite[Corollary~5.17 and Proposition~5.18]{Pf07},
\cf~Definition~\ref{def_pivotalss}, in particular $H_s$ is simple in $\sym{M}^H$.

According to the definition, $H$ is cospherical if, and only if,
$\tr_V^{(L)}(f)=\tr_V^{(R)}(f)$ for all $V\in|\sym{M}^H|$ and all morphisms
$f\colon V\to V$. Since $\sym{M}^H$ is finitely semisimple, this holds if and
only if the condition is satisfied for all simple $V$. If $V$ is simple, however,
then $f=\lambda\,\id_V$ for some $\lambda\in k$. Now we take the
traces~\eqref{eq_tracel} and~\eqref{eq_tracer} for $f=\id_V$. As $H_s$ is
simple, we can evaluate at $h=1$ and obtain~\eqref{eq_tracelsimple}
and~\eqref{eq_tracersimple}. $H$ is cospherical if, and only if, the two traces
are equal for each simple $V$.
\end{proof}

\subsection{Tannaka--Kre\v\i n reconstruction}

In this section, we show that the universal coend $H=\coend(\sym{C},\omega)$
of a finitely semisimple $k$-linear additive spherical category $\sym{C}$ for
which $k=\End(\one)$ is a field, with respect to the long forgetful functor
$\omega\colon\sym{C}\to\Vect_k$, forms a cospherical WHA.

The following theorem was shown in~\cite{Pf07} for modular categories. Its
proof is exactly the same for our spherical categories.

\begin{theorem}[see {\cite[Section 4]{Pf07}}]
\label{thm_reconstruct}
Let $\sym{C}$ be a finitely semisimple $k$-linear additive spherical category,
$k=\End(\one)$ a field, and $\omega\colon\sym{C}\to\Vect_k$ be the long
forgetful functor. Then the universal coend $H=\coend(\sym{C},\omega)$ forms a
finite-dimensional split cosemisimple WHA $(H,\mu,\eta,\Delta,\epsilon,S)$
with $H_t\cap H_s\cong k$.
\end{theorem}

The operations of $H$ are defined in~\cite{Pf07} using the universal property
of the coend. Here we just review how to compute them. As a vector space,
\begin{equation}
\label{eq_coendvect}
H \cong \bigoplus_{j\in I} {\omega(V_j)}^\ast\otimes\omega(V_j)
\end{equation}
where the sum is over a set of representatives $V_j\in|\sym{C}|$, $j\in I$, of
the isomorphism classes of simple objects. We write for the homogeneous
elements of $H$,
\begin{equation}
{[\theta|v]}_X=\theta\otimes v\in{\omega(X)}^\ast\otimes\omega(X).
\end{equation}
We use this notation for all objects $X\in|\sym{C}|$. For all $X$ that do not
equal the chosen representatives of the simple objects, this becomes an
element of~\eqref{eq_coendvect} upon identifying ${[\eta\circ(\id_{\tilde
V}\otimes f)|v]}_X = {[\eta|(\id_{\tilde V}\otimes f)\circ v]}_Y$ for all
$v\in\omega (X)=\Hom(\tilde V,\tilde V\otimes X)$,
$\eta\in{\omega(Y)}^\ast=\Hom(\tilde V\otimes Y,\tilde V)$ and for all
morphisms $f\colon X\to Y$ of $\sym{C}$.

The operations of $H$ are given as follows.
\begin{eqnarray}
\Delta({[\theta|v]}_X)
&=& \sum_j {[\theta|e^{(X)}_j]}_X\otimes{[e^j_{(X)}|v]}_X,\\
\epsilon({[\theta|v]}_X)
&=& \ev_{\omega(X)}(\theta\otimes v),\\
\mu ({[\theta|v]}_X\otimes{[\zeta|w]}_Y)
&=& {[\zeta\circ(\theta\otimes\id_Y)\circ\alpha^{-1}_{\tilde V,X,Y}|
\alpha_{\tilde V,X,Y}\circ(v\otimes\id_Y)\circ w]}_{X\otimes Y},\\
\eta(1) &=& {[\rho_{\tilde V}|\rho^{-1}_{\tilde V}]}_{\one},\\
S({[\theta|v]}_X)
&=& {[\Phi_X(v)|\Psi_X(\theta)]}_{X^\ast},
\end{eqnarray}
where $v\in\omega(X)$, $\theta\in{\omega(X)}^\ast$, $w\in\omega(Y)$ and
$\zeta\in{\omega(Y)}^\ast$, $X,Y\in|\sym{C}|$, and ${(e^{(X)}_j)}_j$ and
${(e^j_{(X)})}_j$ denotes a pair of dual bases of $\omega(X)$ and
${\omega(X)}^\ast$. The maps $\Phi_X$ and $\Psi_X$ are as in
Lemma~\ref{lem_isoanti}. For convenience, we also give the target and source
counital maps and the square of the antipode,
\begin{eqnarray}
\epsilon_t({[\theta|v]}_X)
&=& {[\Phi_X(v)\circ\Psi_X(\theta)\circ\rho_{\tilde V}|\rho_{\tilde V}^{-1}]}_\one,\\
\label{eq_epsilonscomp}
\epsilon_s({[\theta|v]}_X)
&=& {[\rho_{\tilde V}|\rho_{\tilde V}^{-1}\circ\theta\circ v]}_\one,\\
S^2({[\theta|v]}_X)
&=& {[D_{\tilde V}^{-1}\circ\theta\circ(D_{\tilde V}\otimes X)|
(D_{\tilde V}^{-1}\otimes\id_X)\circ v\circ D_{\tilde V}]}_X.
\end{eqnarray}
\cite{Pf07} showed that if $\sym{C}$ is modular, then the WHA $H$ is
coribbon. We now generalize this result to our spherical categories.

\begin{proposition}
Under the assumptions of Theorem~\ref{thm_reconstruct}, the reconstructed WHA
$H=\coend(\sym{C},\omega)$ is cospherical, and its copivotal structure is
given by $w\colon H\to k$ and its convolution inverse $\bar w$ as follows.
\begin{eqnarray}
\label{eq_pivotalform}
w({[\theta|v]}_X) &=&
\ev_{\omega(X)}((D_{\tilde V}^{-1}\circ\theta\circ(D_{\tilde V}\otimes\id_X))
\otimes v),\\
\label{eq_pivotalinv}
\bar w({[\theta|v]}_X) &=&
\ev_{\omega(X)}(\theta\otimes((D_{\tilde V}^{-1}\otimes\id_X)\circ v
\circ D_{\tilde V})).
\end{eqnarray}
\end{proposition}

\begin{proof}
The triangle identities for the left-duals of Proposition~\ref{prop_leftdual}
suffice to show that $w$ and $\bar w$ are mutually convolution inverse and
that $w$ satisfies~\eqref{eq_dualgrouplike} and~\eqref{eq_copivotal}. Thus $H$
is copivotal.

In order to show that $H$ is cospherical, we note that in the usual bases of
$\omega(X)$ and ${\omega(X)}^\ast$ for arbitrary simple $X$, the matrix
elements of the coaction read
\begin{equation}
c_{\ell j}^{(X)} = {[e^\ell_{(X)}|e^{(X)}_j]}_X.
\end{equation}
We compute the left-hand sides of~\eqref{eq_tracelsimple}
and~\eqref{eq_tracersimple}, using $w$ and $\bar w$ of~\eqref{eq_pivotalform}
and~\eqref{eq_pivotalinv}, $\epsilon_s$ of~\eqref{eq_epsilonscomp}, the triangle
identities for $\ev_{\omega(X)}$ and $\coev_{\omega(X)}$, and the
conditions~\eqref{eq_complete1} and~\eqref{eq_complete2} and find
$t^{(L)}_X=\dim X=t^{(R)}_X$, and so $H$ is cospherical by
Theorem~\ref{thm_cospherical}.
\end{proof}

\subsection{Pivotal functors and traces}

In this section, we define the notion of a pivotal functor, \ie, a functor
that preserves the pivotal structure, and see how it behaves with respect to
traces.

The following proposition is well known for strong monoidal functors, but it
cannot even be formulated for functors that are merely lax or merely oplax
monoidal. It nevertheless holds for functors with a Frobenius structure.

\begin{proposition}
\label{prop_functordual}
Let $\sym{C}$ be a left-autonomous and $\sym{C}^\prime$ be a monoidal category
and
\begin{equation}
(F,F_{X,Y},F_0,F^{X,Y},F^0)\colon\sym{C}\to\sym{C}^\prime
\end{equation}
be a functor with Frobenius structure. Then for every $X\in|\sym{C}|$, the
object $FX\in|\sym{C}^\prime|$ has a left-dual
$(F(X^\ast),\ev_{F(X^\ast)},\coev_{F(X^\ast)})$ where
\begin{alignat}{3}
\ev_{F(X^\ast)}   &=&\, F^0\circ F\ev_X\circ F_{X^\ast,X}   &\colon
F(X^\ast)\otimes^\prime FX\to\one^\prime,\\
\coev_{F(X^\ast)} &=&\, F^{X,X^\ast}\circ F\coev_X\circ F_0 &\colon
\one^\prime\to FX\otimes F(X^\ast).
\end{alignat}
\end{proposition}

\begin{proof}
In order to verify the triangle identity~\eqref{eq_triangle1}, we
use~\eqref{eq_frobenius1}, naturality of $F^{X,-}$ and $F_{-,X}$ and the image
under $F$ of the triangle identity in $\sym{C}$. For the other triangle
identity~\eqref{eq_triangle2}, we need~\eqref{eq_frobenius2}, naturality of
$F_{X^\ast,-}$ and $F^{-,X^\ast}$ and the image under $F$ of the triangle
identity in $\sym{C}$.
\end{proof}

If both $\sym{C}$ and $\sym{C}^\prime$ are left-autonomous, then we can use a
standard result (Proposition~\ref{prop_dualisic}) to show that the left-duals
in $\sym{C}^\prime$ and those obtained from the duals in $\sym{C}$ by using
Proposition~\ref{prop_functordual}, are canonically isomorphic.

\begin{corollary}
\label{cor_functordual}
Let $\sym{C}$ and $\sym{C}^\prime$ be left-autonomous categories and
$(F,F_{X,Y},F_0,F^{X,Y},F^0)\colon\sym{C}\to\sym{C}^\prime$ be a functor with
Frobenius structure. Then there are natural isomorphisms $u_X\colon
F(X^\ast)\to {(FX)}^\ast$ given by
\begin{eqnarray}
\label{eq_isoux}
u_X&=&\lambda^\prime_{{(FX)}^\ast}\circ(F^0\otimes^\prime\id_{{(FX)}^\ast})
\circ(F\ev_X\otimes^\prime\id_{{(FX)}^\ast})
\circ(F_{X^\ast,X}\otimes^\prime\id_{{(FX)}^\ast})\nn\\
&&\circ{\alpha^\prime}^{-1}_{F(X^\ast),FX,{(FX)}^\ast}
\circ(\id_{F(X^\ast)}\otimes^\prime\coev_{FX})
\circ{\rho^\prime}^{-1}_{F(X^\ast)}
\end{eqnarray}
with inverse
\begin{eqnarray}
u_X^{-1}&=&\lambda^\prime_{F(X^\ast)}\circ(\ev_{FX}\otimes^\prime\id_{F(X^\ast)})
\circ{\alpha^\prime}^{-1}_{{(FX)}^\ast,FX,F(X^\ast)}
\circ(\id_{{(FX)}^\ast}\otimes^\prime F^{X,X^\ast})\nn\\
&&\circ(\id_{{(FX)}^\ast}\otimes^\prime F\coev_X)
\circ(\id_{{(FX)}^\ast}\otimes^\prime F_0)\circ{\rho^\prime}^{-1}_{{(FX)}^\ast}.
\end{eqnarray}
They satisfy
\begin{equation}
\label{eq_functordualev}
\begin{aligned}
\xymatrix{
F(X^\ast)\otimes^\prime FX\ar[rrrr]^{u_X\otimes\id_{FX}}\ar[dd]_{F_{X^\ast,X}}&&&&
{(FX)}^\ast\otimes^\prime FX\ar[dd]^{\ev_{FX}}\\
\\
F(X^\ast\otimes X)\ar[rr]_{F\ev_X}&&
F\one\ar[rr]_{F^0}&&
\one^\prime,
}
\end{aligned}
\end{equation}
and
\begin{equation}
\label{eq_functordualcoev}
\begin{aligned}
\xymatrix{
\one^\prime\ar[rrrr]^{\coev_{FX}}\ar[dd]_{F_0}&&&&
FX\otimes^\prime{(FX)}^\ast\ar[dd]^{\id_{FX}\otimes u_X^{-1}}\\
\\
F\one\ar[rr]_{F\coev_X}&&
F(X\otimes X^\ast)\ar[rr]_{F^{X,X^\ast}}&&
FX\otimes^\prime F(X^\ast).
}
\end{aligned}
\end{equation}
\end{corollary}

\noindent
Note that $\Psi_X$ of~\eqref{eq_psi} is precisely the $u^{-1}_X$ associated
with the long forgetful functor $\omega$.

\begin{definition}
\label{def_pivotalfunctor}
Let $\sym{C}$ and $\sym{C}^\prime$ be pivotal categories. A strong
monoidal functor
\begin{equation}
(F,F_{X,Y},F_0)\colon\sym{C}\to\sym{C}^\prime
\end{equation}
is called \emph{pivotal} if
\begin{equation}
\begin{aligned}
\xymatrix{
FX\ar[rr]^{F\tau_X}\ar[dd]_{\tau^\prime_{FX}}&&
F({X^\ast}^\ast)\ar[dd]^{u_{X^\ast}}\\
\\
{{(FX)}^\ast}^\ast\ar[rr]_{{(u_X)}^\ast}&&{(F(X^\ast))}^\ast
}
\end{aligned}
\end{equation}
for all $X\in|\sym{C}|$. Here, $u_X\colon F(X^\ast)\to{(FX)}^\ast$ are the
natural isomorphisms of Corollary~\ref{cor_functordual} associated with $F$.
\end{definition}

\begin{example}
Every strong monoidal ribbon functor (see, for example~\cite[Definition
A.14]{Pf07}) between ribbon categories is pivotal.
\end{example}

\begin{proof}
By~\cite[Definition~A.9, Definition~A.14, equation~(A.24)]{Pf07}, the axioms
of a strong monoidal functor and equations~\eqref{eq_functordualev}
and~\eqref{eq_functordualcoev}.
\end{proof}

\noindent
The following is a refinement of Proposition~\ref{prop_monadjunction} in the
case of a pivotal functor.

\begin{proposition}
\label{prop_ribbonadjunction}
Let $\sym{C}$ and $\sym{C}^\prime$ be left-autonomous categories and $F\dashv
G\colon\sym{C}^\prime\to\sym{C}$ be an adjoint equivalence. If $F$ is pivotal
strong monoidal, then $G$ is pivotal.
\end{proposition}

\begin{proof}
First, we observe that if the isomorphisms $u_X$ are defined as
in~\eqref{eq_isoux} for the functor $F$, the analogues of $u_X$ associated
with $G$ turn out to be equal to
\begin{equation}
\label{eq_uvsv}
v_Y=\eta^{-1}_{{(GY)}^\ast}\circ Gu_{GY}^{-1}\circ G\epsilon^\ast_Y\colon
G(Y^\ast)\to {(GY)}^\ast,
\end{equation}
$Y\in|\sym{C}^\prime|$, where $\eta\colon 1_\sym{C}\Rightarrow G\circ F$ and
$\epsilon\colon F\circ G\Rightarrow 1_{\sym{C}^\prime}$ are the unit and
counit of the adjunction, respectively.

For the proof, we need Definition~\ref{def_pivotalfunctor}, the triangle
identities of the adjunction and $v_{Y^\ast}$ and ${(v_Y)}^\ast$
from~\eqref{eq_uvsv}.
\end{proof}

By an equivalence of pivotal categories we mean an equivalence of categories
in which one functor is pivotal strong monoidal. Pivotal strong monoidal
functors $F\colon\sym{C}\to\sym{C}^\prime$ relate the traces of $\sym{C}$ and
$\sym{C}^\prime$ as follows.

\begin{proposition}
Let $\sym{C}$ and $\sym{C}^\prime$ be pivotal categories and
$(F,F_{X,Y},F_0)\colon\sym{C}\to\sym{C}^\prime$ be a pivotal strong
monoidal functor. Then for every morphism $f\colon X\to X$ of
$\sym{C}$,
\begin{equation}
\begin{aligned}
\xymatrix{
\one^\prime\ar[rr]^{F_0}\ar[dd]_{\tr^{(L)}_{FX}(Ff)}&&
F\one\ar[dd]^{F\tr^{(L)}_X(f)}\\
\\
\one^\prime\ar[rr]_{F_0}&&F\one
}
\end{aligned}
\end{equation}
and similarly for the right-trace.
\end{proposition}

\begin{proof}
By Definition~\ref{def_traces}, Definition~\ref{def_pivotalfunctor},
equations \eqref{eq_functordualev} and \eqref{eq_functordualcoev}.
\end{proof}

\begin{proposition}[analogous to {\cite[Proposition A.28]{Pf07}}]
Let $\sym{C}$ and $\sym{C}^\prime$ be semisimple $k$-linear pivotal
categories, $k=\End(\one)$ be a field, and
$(F,F_{X,Y},F_0)\colon\sym{C}\to\sym{C}^\prime$ be a pivotal strong monoidal
$k$-linear functor. Then for each morphism $f\colon X\to X$ of $\sym{C}$,
\begin{equation}
\tr^{(L)}_X(f) = \tr^{(L)}_{FX}(Ff)\qquad\mbox{and}\qquad
\tr^{(R)}_X(f) = \tr^{(R)}_{FX}(Ff).
\end{equation}
In particular, if $F$ is essentially surjective and full and $\sym{C}$ is
spherical, then $\sym{C}^\prime$ is spherical, too.
\end{proposition}

By an equivalence of spherical categories, we therefore mean an equivalence of
categories, one functor of which is pivotal strong monoidal.

\subsection{Equivalence of categories}

Here we show that the original spherical category $\sym{C}$ is equivalent as a
spherical category to $\sym{M}^H$, the category of finite-dimensional
comodules over the universal coend $H=\coend(\sym{C},\omega)$ with respect to
the long forgetful functor.

The proof of the following theorem (shown in~\cite{Pf07} for modular
categories) is exactly identical in the spherical case.

\begin{theorem}[see {\cite[Theorem~6.1]{Pf07}}]
\label{thm_monequiv}
Let $\sym{C}$ be a finitely semisimple $k$-linear additive spherical category,
$k=\End(\one)$ be a field, $\omega\colon\sym{C}\to\Vect_k$ be the long forgetful
functor and $H=\coend(\sym{C},\omega)$ be the reconstructed WHA.
\begin{myenumerate}
\item
The long forgetful functor factors through $\sym{M}^H$, \ie, the diagram
\begin{equation}
\begin{aligned}
\xymatrix{
\sym{C}\ar[rr]^F\ar[ddrr]_{\omega}&&\sym{M}^H\ar[dd]^U\\
\\
&&\Vect_k
}
\end{aligned}
\end{equation}
commutes. Here $U\colon\sym{M}^H\to\Vect_k$ is the forgetful
functor of Proposition~\ref{prop_forgetful}.
\item
The functor $F$ is $k$-linear, essentially surjective and fully
faithful.
\item
$(F,F_{X,Y},F_0)$ forms a strong monoidal functor with
\begin{alignat}{2}
F_{X,Y}&\colon FX\hotimes FY\to F(X\otimes Y),&&\quad
f\otimes g\mapsto \alpha_{\tilde V,X,Y}\circ(f\otimes\id_Y)\circ g,\\
F_0&\colon H_s\mapsto F\one,&&\quad
{[\rho_{\tilde V}|v]}_{\one}\mapsto v.
\end{alignat}
\end{myenumerate}
\end{theorem}

If $\sym{C}$ is modular, then the functor $F$ forms part of an equivalence of
ribbon categories. We now generalize this result to our spherical categories.

\begin{lemma}
Let $\sym{C}$ be a finitely semisimple $k$-linear additive spherical category,
$k=\End(\one)$ a field, and $\omega\colon\sym{C}\to\Vect_k$ be the long
forgetful functor. Then the natural isomorphisms of~\eqref{eq_isoux}
associated with $F$ are given by
\begin{alignat}{3}
\label{eq_isouxmh}
u_X     &\colon F(X^\ast)\to{(F X)}^\ast,&&\quad
\tilde v\mapsto \rho_{\tilde V}\circ(\id_{\tilde V}\otimes\ev_X)
\circ\alpha_{\tilde V,X^\ast,X}\circ(\tilde v\otimes\id_X),\\
u_X^{-1}&\colon{(F X)}^\ast\to F(X^\ast),&&\quad
\theta\mapsto(\theta\otimes\id_{X^\ast})\circ\alpha^{-1}_{\tilde V,X,X^\ast}
\circ(\id_{\tilde V}\otimes\coev_X)\circ\rho_{\tilde V}^{-1}.
\end{alignat}
\end{lemma}

\noindent
The following diagrams illustrate the maps $u_X$ and $u_X^{-1}$:

\begin{equation}
u_X\biggl(
\begin{xy}
(0,0)*+{\hbox to 3em{\hfill $\tilde v$\hfill}}*\frm{-}!U="t" !D="b";
"t"+< 0em,0em>;"t"+< 0em, 2em> **\dir{-}; ?(.5)*\dir{<}; ?(.7)+< 0.7em,0em>*{\tilde V};
"b"+<-1em,0em>;"b"+<-1em,-2em> **\dir{-}; ?(.5)*\dir{>}; ?(.7)+<-0.7em,0em>*{\tilde V};
"b"+< 1em,0em>;"b"+< 1em,-2em> **\dir{-}; ?(.5)*\dir{<}; ?(.7)+< 0.7em,0em>*{X};
\end{xy}
\biggr) :=
\begin{xy}
(0,0)*+{\hbox to 3em{\hfill $\tilde v$\hfill}}*\frm{-}!U="t" !D="b";
"t" +< 0em,  0em>;"t"+<0em,2em> **\dir{-}; ?(.5)*\dir{<}; ?(.7)+< 0.7em,0em>*{\tilde V};
"b"+<-1em,0em>;"b"+<-1em,-2em> **\dir{-}; ?(.5)*\dir{>}; ?(.7)+<-0.7em,0em>*{\tilde V};
"b"+< 1em,0em>;"t"+< 3em, 2em> **\crv{"b"+<1em,-1em>&"b"+<3em,-1em>&"b"+<3em,0em>};
?(.9)*\dir{<}; ?(.95)+<0.7em,0em>*{X};
\end{xy},\qquad
\u_X^{-1}\biggl(
\begin{xy}
(0,0)*+{\hbox to 3em{\hfill $\theta$\hfill}}*\frm{-}!U="t" !D="b";
"t"+<-1em,0em>;"t"+<-1em, 2em> **\dir{-}; ?(.5)*\dir{<}; ?(.7)+<-0.7em,0em>*{\tilde V};
"t"+< 1em,0em>;"t"+< 1em, 2em> **\dir{-}; ?(.5)*\dir{<}; ?(.7)+< 0.7em,0em>*{X};
"b"+< 0em,0em>;"b"+< 0em,-2em> **\dir{-}; ?(.5)*\dir{>}; ?(.7)+< 0.7em,0em>*{\tilde V};
\end{xy}
\biggr) :=
\begin{xy}
(0,0)*+{\hbox to 3em{\hfill $\theta$\hfill}}*\frm{-}!U="t" !D="b";
"t"+<-1em,0em>;"t"+<-1em, 2em> **\dir{-}; ?(.5)*\dir{<}; ?(.7)+<-0.7em,0em>*{\tilde V};
"t"+< 1em,0em>;"b"+< 3em,-2em> **\crv{"t"+<1em,1em>&"t"+<3em,1em>&"t"+<3em,0em>};
?(.9)*\dir{<}; ?(.95)+<0.7em,0em>*{X};
"b"+< 0em,0em>;"b"+< 0em,-2em> **\dir{-}; ?(.5)*\dir{>}; ?(.7)+< 0.7em,0em>*{\tilde V};
\end{xy}
\end{equation}

\begin{proposition}
\label{prop_pivotaleq}
Let $\sym{C}$ be a finitely semisimple $k$-linear additive spherical category,
$k=\End(\one)$ a field. Then the functor $F$ of Theorem~\ref{thm_monequiv} is pivotal.
\end{proposition}

\begin{proof}
Let $H=\coend(\sym{C},\omega)$ be the reconstructed WHA. We first compute
$\tau_{FX}\colon FX\to{{(FX)}^\ast}^\ast$ (\cf~\eqref{eq_pivotalcomodule}) in
$\sym{M}^H$ for the copivotal form of~\eqref{eq_pivotalform}:
\begin{equation}
\tau_{FX}(v)
= \tau_{FX}^{(\Vect_k)}((D_{\tilde V}\otimes\id_X)\circ v\circ D_{\tilde V}^{-1})
\end{equation}
for all $v\in FX$. We also need the dual morphism ${(u_X)}^\ast$
of~\eqref{eq_isouxmh} in $\sym{M}^H$. Therefore, we
compute~\eqref{eq_leftdual} in $\sym{M}^H$ and find
\begin{eqnarray}
{(u_X)}^\ast(\tau_{FX}^{(\Vect_k)}(v))
&=& D_{\tilde V}^{-1}\circ\rho_{\tilde V}\circ(\id_{\tilde V}\otimes\ev_{X^\ast})
\circ(\id_{\tilde V}\otimes(\tau_X\otimes\id_{X^\ast}))\nn\\
&&\circ\alpha_{\tilde V,X,X^\ast}
\circ(v\otimes\id_{X^\ast})
\circ(D_{\tilde V}\otimes\id_{X^\ast})\in {(F(X^\ast))}^\ast
\end{eqnarray}
for all $v\in FX$. Then we can verify that
\begin{eqnarray}
u_{X^\ast}(F\tau_X(v))
&=& \rho_{\tilde V}\circ(\id_{\tilde V}\otimes\ev_{X^\ast})
\circ(\id_{\tilde V}\otimes(\tau_X\otimes\id_{X^\ast}))
\circ\alpha_{\tilde V,X,X^\ast}
\circ(v\otimes\id_{X^\ast})\nn\\
&=& {(u_X)}^\ast(\tau_{FX}(v))
\end{eqnarray}
for all $v\in FX$, and so $F$ is pivotal.
\end{proof}

\noindent
Since the functor $F$ is strong monoidal, this implies
Theorem~\ref{thm_intro1}. Note that this theorem also shows that every
finitely semisimple $k$-linear additive spherical category $\sym{C}$ for which
$k=\End(\one)$ is a field, is abelian and that the long forgetful functor
$\omega\colon\sym{C}\to\Vect_k$ is exact.

\section{Self-duality}
\label{sect_selfdual}

\subsection{The coend as $\End(\tilde V^\ast\otimes\tilde V)$}

Thanks to the finite semisimplicity of our spherical categories $\sym{C}$, the
vector space underlying the coend
\begin{equation}
H=\coend(\sym{C},\omega)=\bigoplus_{j\in I}{\omega(V_j)}^\ast\otimes\omega(V_j)
\end{equation}
is just $\End(\tilde V^\ast\otimes\tilde V)$. Here, $\omega(X)=\Hom(\tilde
V,\tilde V\otimes X)$ and ${\omega(X)}^\ast=\Hom(\tilde V\otimes X,\tilde
V)$. It is instructive to see how our pair of dual bases of $\omega(X)$ and
${\omega(X)}^\ast$ for simple $X$ provides a decomposition of $\tilde
V^\ast\otimes\tilde V$ into a finite biproduct of simple objects.

\begin{proposition}
\label{prop_endvv}
Let $\sym{C}$ be a finitely semisimple $k$-linear additive spherical category,
$k=\End(\one)$ a field and $\omega\colon\sym{C}\to\Vect_k$ be the long
forgetful functor. Define maps $\pi_\ell^{(X)}\colon \tilde
V^\ast\otimes\tilde V\to X$ and $\imath_\ell^{(X)}\colon X\to\tilde
V^\ast\otimes\tilde V$ by
\begin{eqnarray}
\label{eq_defpi}
\pi_\ell^{(X)}
&:=& \lambda_X\circ(\ev_{\tilde V}\otimes \id_X)
\circ\alpha^{-1}_{\tilde V^\ast,\tilde V,X}
\circ(\id_{\tilde V^\ast}\otimes e_\ell^{(X)}),\\
\label{eq_defimath}
\imath_\ell^{(X)}
&:=& \dim X\cdot(\id_{\tilde V^\ast}\otimes D_{\tilde V})
\circ(\id_{\tilde V^\ast}\otimes e^\ell_{(X)})
\circ\alpha_{\tilde V^\ast,\tilde V,X}\nn\\
&&   \circ(\bar\coev_{\tilde V}\otimes\id_X)\circ\lambda^{-1}_X
\end{eqnarray}
where ${(e^{(X)}_\ell)}_\ell$ and ${(e^\ell_{(X)})}_\ell$ form a pair of dual
bases of $\omega(X)$ and ${\omega(X)}^\ast$, $X\in|\sym{C}|$, with respect
to~\eqref{eq_hayashiev}. By $\bar\coev_{\tilde V}$, we denote the coevaluation
map associated with the right-dual of $\tilde V$ as
in~\eqref{eq_barcoev}. Then these maps satisfy
\begin{myenumerate}
\item
The domination property,
\begin{equation}
\sum\limits_{j\in I}\sum\limits_{\ell=1}^{\dim\omega(V_j)}
\imath^{(V_j)}_\ell\circ\pi^{(V_j)}_\ell
=\id_{\tilde V^\ast\otimes\tilde V}.
\end{equation}
\item
If $X,Y\in|\sym{C}|$ are simple and $X\not\cong Y$, then for all $m,\ell$,
\begin{equation}
\pi^{(Y)}_m\circ\imath^{(X)}_\ell=0.
\end{equation}
\item
If $X\in|\sym{C}|$ is simple, then
\begin{equation}
\pi^{(X)}_m\circ\imath^{(X)}_\ell=\delta_{m\ell}\id_X.
\end{equation}
\item
Every morphism $f\colon\tilde V^\ast\otimes\tilde V\to\tilde V^\ast\otimes\tilde V$
is of the form
\begin{equation}
f=\sum\limits_{j\in I}\sum\limits_{\ell,m=1}^{\dim\omega(V_j)}
f^{(V_j)}_{\ell m}\cdot(\imath^{(V_j)}_\ell\circ\pi^{(V_j)}_m),
\end{equation}
where
\begin{equation}
f^{(X)}_{\ell m}=\tr_{X}(\pi^{(X)}_\ell\circ f\circ\imath^{(X)}_m)/\dim X.
\end{equation}
\end{myenumerate}
\end{proposition}

\begin{proof}
We need~\eqref{eq_dominance}
and~\eqref{eq_semisimpledualbases} for $\tilde V^\ast\otimes\tilde V$, the
Schur axiom of Definition~\ref{def_semisimple}(3b), the fact that any $g\colon
X\to X$ for simple $X$ is of the form $g=(\tr_X(g)/\dim_X)\,\id_X$, and the
triangle identities for the duals given in Proposition~\ref{prop_leftdual}.
\end{proof}

\noindent
The following diagrams illustrate the above definitions:
\begin{equation}
\begin{xy}
(0,0)*+{\hbox to 3em{\hfill $\pi^{(X)}_\ell$\hfill}}*\frm{-}!U="t" !D="b";
"t"+<-1em,0em>;"t"+<-1em, 2em> **\dir{-}; ?(.5)*\dir{>}; ?(.7)+<-0.7em,0em>*{\tilde V};
"t"+< 1em,0em>;"t"+< 1em, 2em> **\dir{-}; ?(.5)*\dir{<}; ?(.7)+< 0.7em,0em>*{\tilde V};
"b"+< 0em,0em>;"b"+< 0em,-2em> **\dir{-}; ?(.5)*\dir{>}; ?(.7)+< 0.7em,0em>*{X};
\end{xy} :=
\begin{xy}
(0,0)*+{\hbox to 3em{\hfill $e^{(X)}_\ell$\hfill}}*\frm{-}!U="t" !D="b";
"t"+< 0em,0em>;"t"+< 0em, 2em> **\dir{-}; ?(.5)*\dir{<}; ?(.7)+< 0.7em,0em>*{\tilde V};
"b"+<-1em,0em>;"t"+<-3em, 2em> **\crv{"b"+<-1em,-1em>&"b"+<-3em,-1em>&"b"+<-3em,0em>};
?(.9)*\dir{>}; ?(.95)+<-0.7em,0em>*{\tilde V};
"b"+< 1em,0em>;"b"+< 1em,-2em> **\dir{-}; ?(.5)*\dir{>}; ?(.7)+< 0.7em,0em>*{X};
\end{xy}
\qquad\mbox{and}\qquad
\begin{xy}
(0,0)*+{\hbox to 3em{\hfill $\imath^{(X)}_\ell$\hfill}}*\frm{-}!U="t" !D="b";
"t"+< 0em,0em>;"t"+< 0em, 2em> **\dir{-}; ?(.5)*\dir{<}; ?(.7)+< 0.7em,0em>*{X};
"b"+<-1em,0em>;"b"+<-1em,-2em> **\dir{-}; ?(.5)*\dir{<}; ?(.7)+<-0.7em,0em>*{\tilde V};
"b"+< 1em,0em>;"b"+< 1em,-2em> **\dir{-}; ?(.5)*\dir{>}; ?(.7)+< 0.7em,0em>*{\tilde V};
\end{xy} :=
\begin{xy}
(0,0)*+{\hbox to 3em{\hfill $e_{(X)}^\ell$\hfill}}*\frm{-}!U="t" !D="b";
(0,0)+<0em,-3em>*+{D_{\tilde V}}*\frm{o}!U="ct" !D="cb";
"t"+< 1em,0em>;"t"+< 1em, 2em> **\dir{-}; ?(.5)*\dir{<}; ?(.7)+< 0.7em,0em>*{X};
"t"+<-1em,0em>;"b"+<-3em,-2em> **\crv{"t"+<-1em,1em>&"t"+<-3em,1em>&"t"+<-3em,0em>};
?(.9)*\dir{<}; ?(.95)+<-0.7em,0em>*{\tilde V};
"b";"ct"+<0em,0.2em> **\dir{-}; ?(.6)*\dir{>};
"cb"+<0em,-0.2em>;"cb"+<0em,-1.5em> **\dir{-}; ?(.6)*\dir{>}; ?(.8)+<0.7em,0em>*{\tilde V};
\end{xy}
\cdot\dim X
\end{equation}

\begin{definition}
Under the assumptions of Proposition~\ref{prop_endvv}, we write $H=\End(\tilde
V^\ast\otimes\tilde V)$ and $\hat H=\End(\tilde V\otimes\tilde V^\ast)$ and
define linear isomorphisms $\phi_L,\phi_R\colon H\to\hat H$ (\emph{left-} and
\emph{right-Fourier transform}) by
\begin{eqnarray}
\phi_L(f)
&=& (\lambda_{\tilde V}\otimes\id_{\tilde V^\ast})
\circ((\bar\ev_{\tilde V}\otimes\id_{\tilde V})\otimes\id_{\tilde V^\ast})
\circ(\alpha^{-1}_{\tilde V,\tilde V^\ast,\tilde V}\otimes\id_{\tilde V^\ast})\nn\\
&&\circ((D^{-1}_{\tilde V}\otimes f)\otimes\id_{\tilde V^\ast})
\circ\alpha^{-1}_{\tilde V,\tilde V^\ast\otimes\tilde V,\tilde V^\ast}
\circ(\id_{\tilde V}\otimes\alpha^{-1}_{\tilde V^\ast,\tilde V,\tilde V^\ast})\nn\\
&&\circ(\id_{\tilde V}\otimes(\id_{\tilde V^\ast}\otimes\coev_{\tilde V}))
\circ(\id_{\tilde V}\otimes\rho^{-1}_{\tilde V^\ast}),\\
\phi_R(f)
&=& (\id_{\tilde V}\otimes\rho_{\tilde V^\ast})
\circ(\id_{\tilde V}\otimes(\id_{\tilde V^\ast}\otimes\bar\ev_{\tilde V}))
\circ(\id_{\tilde V}\otimes\alpha_{\tilde V^\ast,\tilde V,\tilde V^\ast})\nn\\
&&\circ\alpha_{\tilde V,\tilde V^\ast\otimes\tilde V,\tilde V^\ast}
\circ((\id_{\tilde V}\otimes f)\otimes D_{\tilde V^\ast}^{-1})
\circ(\alpha_{\tilde V,\tilde V^\ast,\tilde V}\otimes\id_{\tilde V^\ast})\nn\\
&&\circ((\coev_{\tilde V}\otimes\id_{\tilde V})\otimes\id_{\tilde V^\ast})
\circ(\lambda^{-1}_{\tilde V}\otimes\id_{\tilde V^\ast}).
\end{eqnarray}
\end{definition}

The maps $\phi_L$ and $\phi_R$ are identical to the maps of~\cite{BoSz04}
except for the factors $D_{\tilde V}$. Although you might suspect that these
factors are merely a consequence of the choice of the bilinear
form~\eqref{eq_hayashiev}, they are actually present because of the
pivotal form~\eqref{eq_pivotalform} and cannot be avoided by simple
redefinitions. Diagrammatically, the Fourier transforms read
\begin{equation}
\phi_L(\;
\begin{xy}
(0,0)*+{\hbox to 2em{\hfill $f$\strut\hfill}}*\frm{-}!U="t" !D="b";
"t"+<-0.5em,0em>;"t"+<-0.5em, 1em> **\dir{-}; ?(.7)*\dir{>};
"t"+< 0.5em,0em>;"t"+< 0.5em, 1em> **\dir{-}; ?(.3)*\dir{<};
"b"+<-0.5em,0em>;"b"+<-0.5em,-1em> **\dir{-}; ?(.3)*\dir{<};
"b"+< 0.5em,0em>;"b"+< 0.5em,-1em> **\dir{-}; ?(.7)*\dir{>};
\end{xy}\;) =
\begin{xy}
(12,0)*+{\hbox to 2em{\hfill $f$\strut\hfill}}*\frm{-}!U="t" !D="b";
(0,0)*+{D^{-1}_{\tilde V}}*\frm{o}!U="ct" !D="cb";
"t"+<-0.5em,0em>;"t"+<-0.5em, 1em> **\dir{-}; ?(.7)*\dir{>};
"t"+< 0.5em,0em>;"b"+< 2em,-1em>
**\crv{"t"+<0.5em,1em>&"t"+<2em,1em>&"t"+<2em,0em>}; ?(.95)*\dir{<};
"b"+<-0.5em,0em>;"cb"+<0em,-0.2em>
**\crv{"b"+<-0.5em,-1em>&"cb"+<0em,-1em>}; ?(.5)*\dir{<};
"b"+< 0.5em,0em>;"b"+< 0.5em,-1em> **\dir{-}; ?(.7)*\dir{>};
"ct"+<0em,0.2em>;"ct"+<0em,1em> **\dir{-}; ?(.5)*\dir{<};
\end{xy}
\end{equation}
and
\begin{equation}
\phi_R(\;
\begin{xy}
(0,0)*+{\hbox to 2em{\hfill $f$\strut\hfill}}*\frm{-}!U="t" !D="b";
"t"+<-0.5em,0em>;"t"+<-0.5em, 1em> **\dir{-}; ?(.7)*\dir{>};
"t"+< 0.5em,0em>;"t"+< 0.5em, 1em> **\dir{-}; ?(.3)*\dir{<};
"b"+<-0.5em,0em>;"b"+<-0.5em,-1em> **\dir{-}; ?(.3)*\dir{<};
"b"+< 0.5em,0em>;"b"+< 0.5em,-1em> **\dir{-}; ?(.7)*\dir{>};
\end{xy}\;) =
\begin{xy}
(0,0)*+{\hbox to 2em{\hfill $f$\strut\hfill}}*\frm{-}!U="t" !D="b";
(12,0)*+{D^{-1}_{\tilde V}}*\frm{o}!U="ct" !D="cb";
"t"+< 0.5em,0em>;"t"+< 0.5em, 1em> **\dir{-}; ?(.3)*\dir{<};
"t"+<-0.5em,0em>;"b"+<-2em,-1em>
**\crv{"t"+<-0.5em,1em>&"t"+<-2em,1em>&"t"+<-2em,0em>}; ?(.95)*\dir{>};
"b"+<-0.5em,0em>;"b"+<-0.5em,-1em> **\dir{-}; ?(.3)*\dir{<};
"b"+< 0.5em,0em>;"cb"+<0em,-0.2em>
**\crv{"b"+<0.5em,-1em>&"cb"+<0em,-1em>}; ?(.5)*\dir{>};
"ct"+<0em,0.2em>;"ct"+<0em,1em> **\dir{-}; ?(.5)*\dir{>};
\end{xy}
\end{equation}
where all arrows are labeled by $\tilde V$.

\begin{proposition}
\label{prop_structendvv}
Under the assumptions of Proposition~\ref{prop_endvv},
the unital algebra underlying $\coend(\sym{C},\omega)$ is isomorphic to
$H=\End(\tilde V^\ast\otimes\tilde V)$ with the convolution product
\begin{eqnarray}
f\ast g&= &
(\rho^{-1}_{\tilde V^\ast}\otimes\id_{\tilde V})
\circ((\id_{\tilde V^\ast}\otimes\bar\ev_{\tilde V})\otimes\id_{\tilde V})
\circ(\alpha_{\tilde V^\ast,\tilde V,\tilde V^\ast}\otimes\id_{\tilde V})
\circ\alpha_{\tilde V^\ast\otimes\tilde V,\tilde V^\ast,\tilde V}\nn\\
&&\circ((\id_{\tilde V^\ast}\otimes D^{-1}_{\tilde V})\otimes
(\id_{\tilde V^\ast}\otimes\id_{\tilde V}))
\circ(f\otimes g)
\circ\alpha_{\tilde V^\ast\otimes\tilde V,\tilde V^\ast,\tilde V}\nn\\
&&\circ(\alpha^{-1}_{\tilde V^\ast,\tilde V,\tilde V^\ast}\otimes\id_{\tilde V})
\circ((\id_{\tilde V^\ast}\otimes\coev_{\tilde V})\otimes\id_{\tilde V})
\circ(\rho_{\tilde V^\ast}\otimes\id_{\tilde V})
\end{eqnarray}
and unit
\begin{equation}
1_\ast = (D_{\tilde V^\ast}\otimes\id_{\tilde V})\circ\bar\coev_{\tilde
V}\circ\ev_{\tilde V}.
\end{equation}
The coalgebra structure of the coend reads in terms of $\End(\tilde
V^\ast\otimes\tilde V)$ as follows,
\begin{eqnarray}
\Delta(e^{(X)}_{m\ell})
&=& \sum_{p=1}^{\dim\omega(X)} e^{(X)}_{p\ell}\otimes e^{(X)}_{mp},\\
\epsilon(e^{(X)}_{m\ell})
&=& \delta_{\ell m},\\
\end{eqnarray}
for all simple $X\in|\sym{C}|$. Here we have written,
\begin{equation}
\label{eq_defeml}
e^{(X)}_{m\ell} := (\imath^{(X)}_\ell\circ\pi^{(X)}_m)/\dim X,
\end{equation}
with $\imath^{(X)}_\ell$ and $\pi^{(X)}_m$ as in~\eqref{eq_defpi}
and~\eqref{eq_defimath}. The antipode is given by $S=\phi_L^{-1}\circ\phi_R$.
\end{proposition}

\begin{proof}
By a direct computation using the linear isomorphism
\begin{equation}
\coend(\sym{C},\omega)\to\End(\tilde V^\ast\otimes\tilde V),\qquad
{[e^\ell_{(X)}|e^{(X)}_m]}_X\mapsto e^{(X)}_{m\ell}.
\end{equation}
\end{proof}

\noindent
Diagrammatically, the convolution algebra structure of $H=\End(\tilde
V^\ast\otimes\tilde V)$ reads,
\begin{equation}
\begin{xy}
(0,0)*+{\hbox to 2em{\hfill $f$\strut\hfill}}*\frm{-}!U="t" !D="b";
"t"+<-0.5em,0em>;"t"+<-0.5em, 1em> **\dir{-}; ?(.7)*\dir{>};
"t"+< 0.5em,0em>;"t"+< 0.5em, 1em> **\dir{-}; ?(.3)*\dir{<};
"b"+<-0.5em,0em>;"b"+<-0.5em,-1em> **\dir{-}; ?(.3)*\dir{<};
"b"+< 0.5em,0em>;"b"+< 0.5em,-1em> **\dir{-}; ?(.7)*\dir{>};
\end{xy}
\ast
\begin{xy}
(0,0)*+{\hbox to 2em{\hfill $g$\strut\hfill}}*\frm{-}!U="t" !D="b";
"t"+<-0.5em,0em>;"t"+<-0.5em, 1em> **\dir{-}; ?(.7)*\dir{>};
"t"+< 0.5em,0em>;"t"+< 0.5em, 1em> **\dir{-}; ?(.3)*\dir{<};
"b"+<-0.5em,0em>;"b"+<-0.5em,-1em> **\dir{-}; ?(.3)*\dir{<};
"b"+< 0.5em,0em>;"b"+< 0.5em,-1em> **\dir{-}; ?(.7)*\dir{>};
\end{xy}
=
\begin{xy}
(0,0)*+{\hbox to 2em{\hfill $f$\strut\hfill}}*\frm{-}!U="ft" !D="fb";
"ft"+<-0.5em,0em>;"ft"+<-0.5em, 1em> **\dir{-}; ?(.7)*\dir{>};
"fb"+<-0.5em,0em>;"fb"+<-0.5em,-1em> **\dir{-}; ?(.3)*\dir{<};
(16,0)*+{\hbox to 2em{\hfill $g$\strut\hfill}}*\frm{-}!U="gt" !D="gb";
"gt"+< 0.5em,0em>;"gt"+< 0.5em, 1em> **\dir{-}; ?(.3)*\dir{<};
"gb"+< 0.5em,0em>;"gb"+< 0.5em,-1em> **\dir{-}; ?(.7)*\dir{>};
(8,-9)*+{D^{-1}_{\tilde V}}*\frm{o}!L="dl" !R="dr";
"ft"+<0.5em,0em>;"gt"+<-0.5em,0em> **\crv{"ft"+<0.5em,1em>&"gt"+<-0.5em,1em>}
?(.5)*\dir{<};
"fb"+< 0.5em,0em>;"dl"+<0em,0em> **\crv{"fb"+< 0.5em,-1em>} ?(.6)*\dir{>};
"gb"+<-0.5em,0em>;"dr"+<0em,0em> **\crv{"gb"+<-0.5em,-1em>} ?(.4)*\dir{<};
\end{xy}
\qquad\mbox{and}\qquad
1_\ast =
\begin{xy}
(0,-4)*+{D_{\tilde V}}*\frm{o}!L="l" !R="r";
(8,8);(-8,8) **\crv{(8,0)&(-8,0)}; ?(.55)*\dir{>};
(-8,-8);"l" **\crv{(-8,-4)}; ?(.7)*\dir{>};
(8,-8);"r" **\crv{(8,-4)}; ?(.3)*\dir{<};
\end{xy}
\end{equation}
where all arrows are labeled with $\tilde V$. For our purposes, we do not need
the structure of $\End(\tilde V^\ast\otimes\tilde V)$ as a unital associative
algebra with respect to composition. Note that the left-Fourier transform
$\phi_L$ maps composition of $H$ to convolution of $\hat H$ and convolution of
$H$ to opposite composition of $\hat H$, whereas the right-Fourier transform
$\phi_R$ maps them to opposite convolution and composition, respectively.

\subsection{A non-degenerate pairing of WHAs}

In this section, we use Fourier transform in order to establish a dual pairing
of WHAs between $H$ and $\hat H$.

First, observe that all the constructions of Section~\ref{sect_reconstruct},
from the definition of the long forgetful functor
(Definition~\ref{def_longfunctor}) to the WHA structure of the universal coend
(Theorem~\ref{thm_reconstruct}) can be done with $\tilde V^\ast$ rather
than $\tilde V$ as the universal object. This is possible because taking the
dual is an involution on the set of isomorphism classes of simple objects
(Definition~\ref{def_pivotalss}). Similarly,
Proposition~\ref{prop_structendvv} is available for $\tilde V^\ast$ rather
than $\tilde V$, and so we get a WHA structure on
\begin{equation}
\hat H = \End(\tilde V\otimes\tilde V^\ast)
\cong \End(\tilde V^{\ast\ast}\otimes\tilde V^\ast).
\end{equation}
If we denote by ${(\hat e^{(X)}_m)}_m$ and ${(\hat e^\ell_{(X)})}_\ell$ a pair
of dual bases of $\hat\omega(X):=\End(\tilde V^\ast,\tilde V^\ast\otimes X)$
and ${\hat\omega(X)}^\ast:=\End(\tilde V^\ast\otimes X,\tilde V^\ast)$,
respectively, we have the (non-canonical) isomorphism of WHAs
\begin{equation}
\label{eq_isohhat}
H\to\hat H,\quad
{[e^\ell_{(X)}|e^{(X)}_m]}_X\mapsto
{[\hat e^\ell_{(X)}|\hat e^{(X)}_m]}_X,
\end{equation}
just because the construction of both WHAs in Section~\ref{sect_reconstruct}
proceeds identically except for using $\tilde V^\ast$ instead of $\tilde
V$. For all steps in the construction of $\hat H$, we put a hat on the
corresponding symbol that we use for $H$.

\begin{proposition}
\label{prop_pairing}
Under the assumptions of Proposition~\ref{prop_endvv}, there are non-degenerate
bilinear maps
\begin{alignat}{3}
\label{eq_dualpair}
\left<-;\hat-\right> &\colon H\otimes\hat H\to k,&&\quad
f\otimes\hat g\mapsto \tr_{\tilde V\otimes\tilde V^\ast}(\hat g\circ\phi_L(f)),\\
\left<\hat-;-\right> &\colon\hat H\otimes H\to k,&&\quad
\hat f\otimes g\mapsto \tr_{\tilde V^\ast\otimes\tilde V}(g\circ\phi_L^{-1}(\hat f)).
\end{alignat}
\end{proposition}

\begin{proof}
Choose the basis $(e_{j\ell}^{(X)})$ of $H$ as in~\eqref{eq_defeml}
and, analogously, a basis $\hat e^{(X)}_{j\ell}$ of $\hat H$. Compute
$\phi_L$ and $\phi_L^{-1}$ in these bases, using
Proposition~\ref{prop_endvv}(4). Then the conditions
$\phi_L^{-1}\circ\phi_L=\id_H$ and $\phi_L\circ\phi_L^{-1}=\id_{\hat H}$ imply
the following orthogonality relations:
\begin{eqnarray}
\dim Y\,\sum_{j\in I}\dim V_j\,\sum_{p,q=1}^{\dim\omega(V_j)}
\left<e^{(X)}_{\ell m};\hat e^{(V_j)}_{qp}\right>
\left<\hat e^{(V_j)}_{pq};e^{(Y)}_{sr}\right>
&=& \delta_{XY}\delta_{r\ell}\delta_{sm},\\
\dim Y\,\sum_{j\in I}\dim V_j\,\sum_{p,q=1}^{\dim\omega(V_j)}
\left<\hat e^{(X)}_{\ell m};e^{(V_j)}_{qp}\right>
\left<e^{(V_j)}_{pq};\hat e^{(Y)}_{sr}\right>
&=& \delta_{XY}\delta_{r\ell}\delta_{sm},
\end{eqnarray}
for all $\ell,m,r,s$, upon comparing coefficients, for any simple
$X,Y\in|\sym{C}|$. By $\delta_{XY}$, we mean that $\delta_{XY}=1$ if $X\cong
Y$ and $\delta_{XY}=0$ otherwise. We finally get the canonical element
\begin{eqnarray}
\label{eq_canelement}
G\colon k&\to&\hat H\otimes H,\nn\\
1&\mapsto&
\sum_{i,j\in I}\sum_{p,q=1}^{\dim\omega(V_i)}\sum_{r,s=1}^{\dim\omega(V_j)}
\hat e^{(V_i)}_{pq}\otimes e^{(V_j)}_{rs}\,
\left<\hat e^{(V_i)}_{qp};e^{(V_j)}_{sr}\right>
\dim V_i\dim V_j,
\end{eqnarray}
that satisfies the triangle identities together with $\left<-;\hat -\right>$,
and a similar one for $\left<\hat -;-\right>$, establishing non-degeneracy.
\end{proof}

The following diagrams illustrate the bilinear maps:
\begin{equation}
\left<e^{(X)}_{j\ell};\hat e^{(Y)}_{pq}\right> =
\begin{xy}
(8,24)*+{\hbox to 16mm{\hfill\strut $\pi^{(X)}_j$\hfill}}*\frm{-}!U="pt" !D="pb";
(8,8)*+{\hbox to 16mm{\hfill\strut $\imath^{(X)}_\ell$\hfill}}*\frm{-}!U="it" !D="ib";
"pb";"it" **\dir{-}; ?(.6)*\dir{>}; ?(.5)+<.6em,0em>*{X};
(20,-8)*+{\hbox to 16mm{\hfill\strut $\hat\pi^{(Y)}_p$\hfill}}*\frm{-}!U="hpt" !D="hpb";
(20,-24)*+{\hbox to 16mm{\hfill\strut $\hat\imath^{(Y)}_q$\hfill}}*\frm{-}!U="hit" !D="hib";
"hpb";"hit" **\dir{-}; ?(.6)*\dir{>}; ?(.5)+<.6em,0em>*{Y};
"ib"+<6mm,0mm>;"hpt"+<-6mm,0mm> **\dir{-}; ?(.6)*\dir{>};
"pt"+<18mm,0mm>;"pt"+<6mm,0mm> **\crv{"pt"+<18mm,6mm>&"pt"+<6mm,6mm>};
"hpt"+<6mm,0mm>;"pt"+<18mm,0mm> **\dir{-}; ?(.6)*\dir{>};
"pt"+<-6mm,0mm>;"pt"+<30mm,0mm> **\crv{"pt"+<-6mm,12mm>&"pt"+<30mm,12mm>};
"pt"+<30mm,0mm>;"hib"+<18mm,0mm> **\dir{-}; ?(.55)*\dir{>};
"hib"+<18mm,0mm>;"hib"+<6mm,0mm> **\crv{"hib"+<18mm,-6mm>&"hib"+<6mm,-6mm>};
(2,-8)*+{D^{-1}_{\tilde V}}*\frm{o}!U="dt" !D="db";
"dt";"ib"+<-6mm,0mm> **\dir{-}; ?(.6)*\dir{>};
"hib"+<-6mm,0mm>;"hib"+<-18mm,0mm> **\crv{"hib"+<-6mm,-6mm>&"hib"+<-18mm,-6mm>};
"hib"+<-18mm,0mm>;"db" **\dir{-}; ?(.6)*\dir{>};
\end{xy}\,\,/(\dim X\dim Y),
\end{equation}
\begin{equation}
\left<\hat e^{(X)}_{j\ell};e^{(Y)}_{pq}\right> =
\begin{xy}
(8,-24)*+{\hbox to 16mm{\hfill\strut $\imath^{(Y)}_q$\hfill}}*\frm{-}!U="pt" !D="pb";
(8,-8)*+{\hbox to 16mm{\hfill\strut $\pi^{(Y)}_p$\hfill}}*\frm{-}!U="it" !D="ib";
"pb";"it" **\dir{-}; ?(.4)*\dir{<}; ?(.5)+<.6em,0em>*{Y};
(20,8)*+{\hbox to 16mm{\hfill\strut $\hat\imath^{(X)}_\ell$\hfill}}*\frm{-}!U="hpt" !D="hpb";
(20,24)*+{\hbox to 16mm{\hfill\strut $\hat\pi^{(X)}_j$\hfill}}*\frm{-}!U="hit" !D="hib";
"hpb";"hit" **\dir{-}; ?(.4)*\dir{<}; ?(.5)+<.6em,0em>*{X};
"it"+<6mm,0mm>;"hpb"+<-6mm,0mm> **\dir{-}; ?(.4)*\dir{<};
"pb"+<18mm,0mm>;"pb"+<6mm,0mm> **\crv{"pb"+<18mm,-6mm>&"pb"+<6mm,-6mm>};
"hpb"+<6mm,0mm>;"pb"+<18mm,0mm> **\dir{-}; ?(.4)*\dir{<};
"pb"+<-6mm,0mm>;"pb"+<30mm,0mm> **\crv{"pb"+<-6mm,-12mm>&"pb"+<30mm,-12mm>};
"pb"+<30mm,0mm>;"hit"+<18mm,0mm> **\dir{-}; ?(.45)*\dir{<};
"hit"+<18mm,0mm>;"hit"+<6mm,0mm> **\crv{"hit"+<18mm,6mm>&"hit"+<6mm,6mm>};
(2,8)*+{D_{\tilde V}}*\frm{o}!U="dt" !D="db";
"db";"it"+<-6mm,0mm> **\dir{-}; ?(.4)*\dir{<};
"hit"+<-6mm,0mm>;"hit"+<-18mm,0mm> **\crv{"hit"+<-6mm,6mm>&"hit"+<-18mm,6mm>};
"hit"+<-18mm,0mm>;"dt" **\dir{-}; ?(.4)*\dir{<};
\end{xy}\,\,/(\dim X\dim Y).
\end{equation}
Here, all unlabeled arrows refer to the object $\tilde V$. Expanding all
occurrences of $\tilde V$ and ${\tilde V}^\ast$ as direct sums of simple
objects shows that both bilinear forms are just sums of generalized
$6j$-symbols associated with the category $\sym{C}$.

\begin{definition}
Let $H$ and $L$ be WHAs over some field $k$. A \emph{dual pairing} of $H$ and
$L$ is a non-degenerate $k$-bilinear map $g\colon H\otimes L\to k$,
$h\otimes\ell\mapsto g(h;\ell)$ such that
\begin{eqnarray}
g(\Delta h;\ell_1\otimes\ell_2) &=& g(h;\ell_1\ell_2),\\
\epsilon(h)                     &=& g(h;1),\\
g(h_1h_2;\ell)                  &=& g(h_1\otimes h_2;\Delta\ell),\\
g(1;\ell)                       &=& \epsilon(\ell),\\
g(Sh;\ell)                      &=& g(h;S\ell),
\end{eqnarray}
for all $h,h_1,h_2\in H$ and $\ell,\ell_1,\ell_2\in L$.
\end{definition}

\noindent
Note that we extend $g$ to tensor products by $g(h_1\otimes
h_2;\ell_1\otimes\ell_2):=g(h_1;\ell_1)\,g(h_2;\ell_2)$.

\begin{theorem}
\label{thm_duality}
Under the assumptions of Proposition~\ref{prop_endvv}, the map $\left<-;\hat
-\right>$ of~\eqref{eq_dualpair} forms a dual pairing of WHAs.
\end{theorem}

\begin{proof}
Verify the conditions for the bases $e^{(X)}_{j\ell}$ of $H$ and $\hat
e^{(X)}_{j\ell}$ of $\hat H$, \cf~\eqref{eq_defeml}, and use the triangle
identities for evaluation and coevaluation of Proposition~\ref{prop_leftdual}.
\end{proof}

Combining Theorem~\ref{thm_duality} with the trivial fact that there exists
the isomorphism~\eqref{eq_isohhat}, we see that the WHA $H$ is self-dual. It
is instructive to give the basis of $\hat H$ dual to our basis
${(e^{(X)}_{j\ell})}_{Xj\ell}$ of $H$ with respect to $\left<-;\hat -\right>$:
\begin{equation}
\hat E^{(X)}_{j\ell} := \dim X\,\sum_{j\in I}\dim V_j\,
\sum_{p,q=1}^{\dim\omega(V_j)} \hat e^{(V_j)}_{pq}
\left<\hat e^{(V_j)}_{qp};e^{(X)}_{\ell j}\right>\in\hat H,
\end{equation}
where $X\in|\sym{C}|$ is simple. It can be read off the canonical
element~\eqref{eq_canelement} and satisfies $\left<e^{(Z)}_{rs};\hat
E^{(X)}_{j\ell}\right>=\delta_{XZ}\delta_{rj}\delta_{s\ell}$ for all simple
$X,Z\in|\sym{C}|$ and all $j,\ell,r,s$.

\section{Self-duality of spherical categories}
\label{sect_repmon}

\subsection{Duals of monoidal categories}

We first generalize the notion of the dual of a monoidal category
of~\cite{Ma91,Ma92} from strong monoidal functors to functors with a separable
Frobenius structure.

\begin{definition}
Let $\sym{C}$ and $\sym{V}$ be monoidal categories and
$(F,F_{X,Y},F_0,F^{X,Y},F^0)\colon\sym{C}\to\sym{V}$ be a functor with
separable Frobenius structure.
\begin{myenumerate}
\item
A \emph{right $(\sym{C},F)$-module} $(V,c_V)$ consists of an object
$V\in|\sym{V}|$ and a natural transformation $c_V\colon V\otimes
F\Rightarrow F\otimes V$ such that
\begin{equation}
\label{eq_moduleax1}
\begin{aligned}
\xymatrix{
V\otimes\one\ar[rr]^{\id_V\otimes F_0}\ar[d]_{\rho_V}&&
V\otimes F\one\ar[dd]^{{(c_V)}_{\one}}\\
V\ar[d]_{\lambda_V^{-1}}\\
\one\otimes V&&
F\one\otimes V\ar[ll]^{F^0\otimes\id_V}
}
\end{aligned}
\end{equation}
and
\begin{equation}
\label{eq_moduleax2}
\begin{aligned}
\xymatrix{
V\otimes(FX\otimes FY)\ar[rr]^{\id_V\otimes F_{X,Y}}\ar[dd]_{\alpha_{V,FX,FY}^{-1}}&&
V\otimes F(X\otimes Y)\ar[rr]^{{(c_V)}_{X\otimes Y}}&&
F(X\otimes Y)\otimes V\ar[dd]^{F^{X,Y}\otimes\id_V}\\
\\
(V\otimes FX)\otimes FY\ar[dd]_{{(c_V)}_X\otimes\id_{FY}}&&
&&
(FX\otimes FY)\otimes V\ar[dd]^{\alpha_{FX,FY,V}}\\
\\
(FX\otimes V)\otimes FY\ar[rr]_{\alpha_{FX,V,FY}}&&
FX\otimes (V\otimes FY)\ar[rr]_{\id_{FX}\otimes{(c_V)}_Y}&&
FX\otimes (FY\otimes V)
}
\end{aligned}
\end{equation}
commute for all $X,Y\in|\sym{C}|$.
\item
A \emph{morphism} $\phi\colon(V,c_V)\to(W,c_W)$ of right
$(\sym{C},F)$-modules is a morphism $\phi\colon V\to W$ of $\sym{V}$ such
that
\begin{equation}
\label{eq_morphax}
\begin{aligned}
\xymatrix{
V\otimes FX\ar[rr]^{\phi\otimes\id_{FX}}\ar[dd]_{{(c_V)}_X}&&
W\otimes FX\ar[dd]^{{(c_W)}_X}\\
\\
FX\otimes V\ar[rr]_{\id_{FX}\otimes\phi}&& FX\otimes W
}
\end{aligned}
\end{equation}
commutes for all $X\in|\sym{C}|$.
\item
The category ${(\sym{C},F)}^\ast$ whose objects and morphisms are the right
$(\sym{C},F)$-modules and their morphisms, is called the \emph{full
right-dual of $\sym{C}$ over $F$}.
\item
The full subcategory ${(\sym{C},F)}^\circ$ of ${(\sym{C},F)}^\ast$ whose
objects are those right $(\sym{C},F)$-modules $(V,c_V)$ for which $c_V$ is a
natural equivalence, is called the \emph{right-dual of $\sym{C}$ over $F$}.
\item
If $\sym{V}=\Vect_k$ and $k$ is a field, we denote by
$\sym{M}_{(\sym{C},F)}$ the full subcategory of ${(\sym{C},F)}^\ast$ whose
objects are those right $(\sym{C},F)$-modules $(V,c_V)$ for which $V$ is
finite-dimensional.
\end{myenumerate}
\end{definition}

In order to keep this section brief, we do not develop the abstract theory of
the dual category nor do we say how ${(\sym{C},F)}^\circ$ and
$\sym{M}_{(\sym{C},F)}$ are related. We merely show for the case of our
spherical categories and the reconstructed WHAs $H$ that the categories
$\sym{M}_{(\sym{M}^H,U)}\cong\sym{M}_H$ are isomorphic as categories (without
extra structure). Here, $U\colon\sym{M}^H\to\Vect_k$ is the forgetful functor
of Proposition~\ref{prop_forgetful}. Since we know the WHA $H$ in detail, it
is not difficult to see in a second step how the extra structure of
$\sym{M}_H$ as a pivotal category equips $\sym{M}_{(\sym{M}^H,U)}$ with the
structure of a pivotal category.

\begin{proposition}
\label{prop_dualcategory}
Let $H$ be a WBA over the field $k$ and $U\colon\sym{M}^H\to\Vect_k$ be the
forgetful functor of the category $\sym{M}^H$ of finite-dimensional
right $H$-comodules (Proposition~\ref{prop_forgetful}).
\begin{myenumerate}
\item
There is a functor
\begin{equation}
\Phi\colon\sym{M}_{(\sym{M}^H,U)}\to\sym{M}_H
\end{equation}
given as follows. $\Phi$ assigns to every right $(\sym{M}^H,U)$-module
$(V,c_V)$ the right $H$-module whose underlying vector space is $V$ with the
action $\gamma_V\colon V\otimes UH\to V$ defined by commutativity of
\begin{equation}
\label{eq_defgamma}
\begin{aligned}
\xymatrix{
V\otimes UH\ar[rr]^{{(c_V)}_H}\ar[dd]_{\gamma_V}&&
UH\otimes V\ar[dd]^{\epsilon\otimes\id_V}\\
\\
V&&
k\otimes V\ar[ll]^{\lambda_V}
}
\end{aligned}
\end{equation}
Here we write $UH$ for the vector space underlying the WHA $H$ where $H$ is
viewed as the regular right $H$-comodule. The functor $\Phi$ assigns to each
morphism $\phi\colon(V,c_V)\to(W,c_W)$ of right-$(\sym{M}^H,U)$-modules the
underlying $k$-linear map $\phi\colon V\to W$ that forms a morphism of
right $H$-modules.
\item
There is a functor
\begin{equation}
\Psi\colon\sym{M}_H\to\sym{M}_{(\sym{M}^H,U)}
\end{equation}
given as follows. $\Psi$ assigns to each right $H$-module $\gamma_V\colon
V\otimes UH\to V$ the right $(\sym{M}^H,U)$-module with the same underlying
vector space $V$ and $c_V$ defined by commutativity of
\begin{equation}
\label{eq_defcv}
\begin{aligned}
\xymatrix{
V\otimes UX\ar[rr]^{\id_V\otimes\beta_X}\ar[dd]_{{(c_V)}_X}&&
V\otimes(UX\otimes UH)\ar[rr]^{\alpha_{V,UX,UY}^{-1}}&&
(V\otimes UX)\otimes UH\ar[dd]^{\tau_{V,UX}\otimes\id_{UH}}\\
\\
UX\otimes V&&
UX\otimes(V\otimes UH)\ar[ll]^{\id_{UX}\otimes\gamma_V}&&
(UX\otimes V)\otimes UH\ar[ll]^{\alpha_{UX,V,UH}}
}
\end{aligned}
\end{equation}
for all $X\in|\sym{M}^H|$. Here, $\beta_X\colon UX\to UX\otimes UH$ denotes
the comodule structure of $X$. The functor $\Psi$ assigns to each
morphism $\phi\colon V\to W$ of right $H$-modules the morphism
$\phi\colon(V,c_V)\to(W,c_W)$ of right $(\sym{M}^H,U)$-modules with the same
underlying $k$-linear map.
\item
The composition $\Phi\circ\Psi=1_{\sym{M}_H}$ is the identity functor.
\end{myenumerate}
\end{proposition}

\begin{proof}
For the separable Frobenius structure of $U$, see
Proposition~\ref{prop_forgetful}.
\begin{myenumerate}
\item
We claim that if $(V,c_V)$ is a right $(\sym{M}^H,U)$-module, then $V$ forms
a right $H$-module with the action $\gamma_V$ of~\eqref{eq_defgamma}.

In order to show that $\gamma_V\circ(\id_V\otimes\eta)=\rho_V$, we use the
fact that $\eta=U\imath\circ U_0$ where $\imath\colon\one=H_s\to H$ is the
inclusion; the definition~\eqref{eq_defgamma} of $\gamma_V$; the fact that
$c_V$ is natural and the inclusion $\imath\colon H_s\to H$ is a morphism of
right $H$-comodules; the identity $\epsilon\circ U\imath=U^0$; and the
axiom~\eqref{eq_moduleax1}.

In order to show that
$\gamma_V\circ(\id_V\otimes\mu)\circ\alpha_{V,UH,UH}=\gamma_V\circ(\gamma_V\otimes\id_{UH})$,
we start with the left-hand side and use the definition~\eqref{eq_defgamma};
the identity $\mu=\tilde\mu\circ U_{H,H}$ where $\tilde\mu=\mu\circ U^{H,H}$
and $H\otimes H$ denotes the tensor product of two copies of the regular
right $H$-comodule in $\sym{M}^H$; the fact that $c_V$ is natural and
$\tilde\mu\colon H\otimes H\to H$ is a morphism of right $H$-comodules; the
identity $\epsilon\circ\tilde\mu=(\epsilon\otimes\epsilon)\circ U^{H,H}$;
the axiom~\eqref{eq_moduleax2}; and twice the definition~\eqref{eq_defgamma}
of $\gamma_V$ again.

Furthermore, if $\phi\colon(V,c_V)\to(W,c_W)$ is a morphism of right
$(\sym{M}^H,U)$-modules, then $\phi\colon V\to W$ is a morphism of right
$H$-modules. In order to see this, we use the axiom~\eqref{eq_morphax} and
the definition~\eqref{eq_defgamma} for $\gamma_V$ and $\gamma_W$.
\item
We claim that if $\gamma_V\colon V\otimes UH\to V$ is a right $H$-module
structure on $V$, then $(V,c_V)$ with $c_V$ as in~\eqref{eq_defcv} forms a
right $(\sym{M}^H,U)$-module.

In order to see that $c_V$ is natural for some morphism $f\colon X\to Y$ of
$\sym{M}^H$, we need the condition that $f$ is a morphism and the
definition~\eqref{eq_defcv} for $c_V$ and $c_W$.

In order to verify~\eqref{eq_moduleax1}, we need the
definition~\eqref{eq_defcv} for $X=\one$; the identity
$\lambda_{UH}\circ(U^0\otimes\id_{UH})\circ\beta_\one\circ U_0=\eta$; and
the condition $\gamma_V\circ(\id_V\otimes\eta)=\rho_V$.

In order to verify~\eqref{eq_moduleax2}, we need the
definition~\eqref{eq_defcv} for both ${(c_V)}_X$ and ${(c_V)}_Y$ and the
condition
$\gamma_V\circ(\id_V\otimes\mu)\circ\alpha_{V,UH,UH}=\gamma_V\circ(\gamma_V\otimes\id_{UH})$.

Furthermore, if $\phi\colon V\to W$ is a morphism of right $H$-modules, then
$\phi\colon(V,c_V)\to (W,c_W)$ is a morphism of right
$(\sym{M}^H,U)$-modules. In order to see this, we need the condition that
$\phi$ is a morphism of right $H$-modules and the
definition~\eqref{eq_defcv} for $c_V$ and~$c_W$.
\item
In order to show that $\Phi\circ\Psi=1_{\sym{M}_H}$, let $\gamma_V\colon
V\otimes UH\to V$ define a right $H$-module. Then $c_V$ of~\eqref{eq_defcv}
is a right $(\sym{M}^H,U)$-module, and~\eqref{eq_defgamma} defines another
right $H$-module structure which we now call $\tilde\gamma_V\colon V\otimes
UH\to V$. We verify that $\gamma_V=\tilde\gamma_V$ by using the definition of
the regular right $H$-comodule structure, \ie, $\beta_H=\Delta$, and the
fact that $H$ forms a counital coalgebra.
\end{myenumerate}
\end{proof}

\begin{theorem}
Let $\sym{C}$ be a finitely semisimple $k$-linear additive spherical category,
$k=\End(\one)$ be a field, $\omega\colon\sym{C}\to\Vect_k$ be the long forgetful
functor, $H=\coend(\sym{C},\omega)$ and $U\colon\sym{M}^H\to\Vect_k$ be the
forgetful functor. Then the functors $\Phi$ and $\Psi$ of
Proposition~\ref{prop_dualcategory} satisfy in addition that
$\Psi\circ\Phi=1_{\sym{M}_{(\sym{M}^H,U)}}$, \ie\ the categories
$\sym{M}_H\cong\sym{M}_{(\sym{M}^H,U)}$ are isomorphic.
\end{theorem}

\begin{proof}
Whereas Proposition~\ref{prop_dualcategory} uses only the abstract properties
of WBAs, the present theorem requires some knowledge of how to reconstruct $H$
from $\sym{M}^H$. Since $H$ is split cosemisimple,
\begin{equation}
H=\bigotimes_{j\in I}{(UV_j)}^\ast\otimes UV_j
\end{equation}
is a direct sum of matrix coalgebras ${(UV_j)}^\ast\otimes UV_j$. For each
simple $X\in|\sym{C}|$, there are therefore homomorphisms of coalgebras
$\imath^X\colon {(UX)}^\ast\otimes UX\to H$ and $\pi^X\colon
H\to{(UX)}^\ast\otimes UX$ such that
\begin{equation}
\id_H=\sum_{j\in I}\imath^{V_j}\circ\pi^{V_j}\qquad\mbox{and}\qquad
\pi^{V_j}\circ\imath^{V_\ell}=\delta_{j\ell}\,\id_{{(UV_j)}^\ast\otimes UV_j}.
\end{equation}
If $H$ is viewed as the regular right $H$-comodule, we have as vector spaces:
\begin{equation}
H=\bigotimes_{j\in I} k^{\dim UV_j}\otimes UV_j,
\end{equation}
and both $\imath^X$ and $\pi^X$ are morphisms of right $H$-comodules. Since
$H$ coacts only on the right tensor factor,
\begin{equation}
\begin{aligned}
\xymatrix{
V\otimes ({(UX)}^\ast\otimes UX)\ar[rr]^{\alpha^{-1}}
\ar[dd]_{{(c_V)}_{{(UX)}^\ast\otimes UX}}&&
(V\otimes {(UX)}^\ast)\otimes UX\ar[rr]^{\tau\otimes\id}&&
({(UX)}^\ast\otimes V)\otimes UX\ar[dd]^{\alpha}\\
\\
({(UX)}^\ast\otimes UX)\otimes V&&
{(UX)}^\ast\otimes(UX\otimes V)\ar[ll]^{\alpha}&&
{(UX)}^\ast\otimes (V\otimes UX)\ar[ll]^{\id\otimes {(c_V)}_X}
}
\end{aligned}
\end{equation}
commutes for each simple $X\in|\sym{C}|$. We can therefore compute
\begin{equation}
\label{eq_cvh}
{(c_V)}_H=\sum_{j\in I}(U\imath^{V_j}\otimes\id_V)\circ{(c_V)}_{{(UV_j)}^\ast\otimes UV_j}
\circ(\id_V\otimes U\pi^{V_j})
\end{equation}
in terms of the ${(c_V)}_X$ for the simple $X\in|\sym{C}|$.

Given some right $(\sym{M}^H,U)$-module $(V,c_V)$, there is a right $H$-module
$\gamma_V\colon V\otimes UH\to V$ given by~\eqref{eq_defgamma} and another
right $(\sym{M}^H,U)$-module from~\eqref{eq_defcv} which we now call $(V,\tilde
c_V)$. Expressing ${(\tilde c_V)}_X$ first in terms of $\gamma_V$, then in
terms of ${(c_V)}_H$, and finally in terms of the ${(c_V)}_Y$ for the simple
$Y\in|\sym{C}|$ using~\eqref{eq_cvh} shows that $c_V=\tilde c_V$.
\end{proof}

\subsection{Pivotal structure}

So far, we have an isomorphism $\sym{M}_H\cong\sym{M}_{(\sym{M}^H,U)}$ of
categories (without extra structure). Putting the structure of a pivotal
category on $\sym{M}_{(\sym{M}^H,U)}$ is straightforward because we can show
that $\sym{M}_H$ is a pivotal category. The subsequent constructions are dual
to those for $\sym{M}^H$.

\begin{proposition}
Let $H$ be a WBA. Then the category $\sym{M}_H$ of finite-dimensional right
$H$-modules is a monoidal category
$(\sym{M}_H,\hotimes,H_s,\alpha,\lambda,\rho)$. The monoidal unit object is
the source base algebra $H_s$ with the action
\begin{equation}
\gamma_{H_s}\colon H_s\otimes H\to H_s,\quad x \otimes h\mapsto
\ract{x}{h}:=\epsilon_s(xh).
\end{equation}
The tensor product $V\hotimes W$ of two right $H$-modules is the vector space
\begin{equation}
V\hotimes W:=\{\,v\otimes w\in V\otimes W\mid\quad v\otimes
w=(\ract{v}{1^\prime})\otimes(\ract{w}{1^\pprime})\,\}
\end{equation}
with the action
\begin{equation}
\ract{(v\otimes w)}{h}
:= (\ract{v}{h^\prime})\otimes (\ract{w}{h^\pprime}).
\end{equation}
The unit constraints are given by
\begin{alignat}{3}
\lambda_V &\colon H_s\hotimes V\to V,&&\quad
h\otimes v\mapsto\epsilon(h1^\prime)\,(\ract{v}{1^\pprime}),\\
\rho_V    &\colon V\hotimes H_s\to V,&&\quad
v\otimes h\mapsto\ract{v}{h},
\end{alignat}
and the associator is induced from that of $\Vect_k$.
\end{proposition}

\begin{proposition}
Let $H$ be a WHA. Then the category $\sym{M}_H$ is left-autonomous if the
left-dual of every object $V\in|\sym{M}_H|$ is chosen to be
$(V^\ast,\ev_V,\coev_V)$, where the dual vector space $V^\ast$ is equipped with
the action
\begin{equation}
\gamma_{V^\ast}\colon V^\ast\otimes H\to V^\ast,\qquad
\theta\otimes h\mapsto(v\mapsto\theta(\ract{v}{(Sh)}),
\end{equation}
and evaluation and coevaluation are given by
\begin{alignat}{3}
\ev_V  &\colon V^\ast\hotimes V\mapsto H_s,&&\quad
\theta\otimes v\mapsto 1^\prime\,\theta(\ract{v}{(S1^\pprime)}),\\
\coev_V&\colon H_s\mapsto V\hotimes V^\ast,&&\quad
h\mapsto\sum_j (\ract{v_j}{(Sh)})\otimes\theta^j.
\end{alignat}
\end{proposition}

\begin{definition}
Let $H$ be a WHA. An element $m\in H$ is called \emph{group-like} if it has a
multiplicative inverse and
\begin{equation}
(m1^\prime)\otimes(m1^\pprime) = m^\prime\otimes m^\pprime
= (1^\prime m)\otimes(1^\pprime m).
\end{equation}
\end{definition}

\noindent
Note that every group-like element $m\in H$ also satisfies $m=Sm$ and
$\epsilon_s(m)=1=\epsilon_t(m)$.

\begin{definition}
A \emph{pivotal} WHA $(H,\mu,\eta,\Delta,\epsilon,S,m)$ is a WHA
$(H,\mu,\eta,\Delta,\epsilon,S)$ with a group-like element $m\in H$, called
the \emph{pivotal element}, that satisfies
\begin{equation}
S^2(x) = mxm^{-1}
\end{equation}
for all $x\in H$.
\end{definition}

\begin{proposition}
Let $(H,\mu,\eta,\Delta,\epsilon,S,m)$ be a pivotal WHA. Then the category
$\sym{M}_H$ is pivotal with $\tau_V\colon V\to {V^\ast}^\ast$ given by
\begin{equation}
\tau_V(v) = \tau_V^{\Vect_k}(\ract{v}{m})
\end{equation}
for all finite-dimensional right $H$-modules $V\in|\sym{M}_H|$ and all $v\in
V$.
\end{proposition}

\begin{definition}
A \emph{spherical} WHA $H$ is a pivotal WHA for which
$\tr_V^{(L)}(f)=\tr_V^{(R)}(f)$ for all finite-dimensional right $H$-modules
$V\in|\sym{M}_H|$ and all morphisms $f\colon V\to V$.
\end{definition}

\begin{proposition}
Let $\sym{C}$ be a finitely semisimple $k$-linear additive spherical
category, $k=\End(\one)$ be a field and $\omega\colon\sym{C}\to\Vect_k$ be the
long forgetful functor. Then $H=\coend(\sym{C},\omega)$ is a pivotal WHA with
pivotal element
\begin{equation}
\label{eq_pivotalelement}
m = {[D_{\tilde V}\circ\rho_V|\rho_V^{-1}\circ D_{\tilde V}^{-1}]}_\one.
\end{equation}
\end{proposition}

\begin{proof}
Direct computation.
\end{proof}

Note that both of the functors, $\Phi$ and $\Psi$, that form the isomorphism
$\sym{M}_H\cong\sym{M}_{(\sym{M}^H,U)}$, leave the vector spaces underlying
the objects and the linear maps underlying the morphisms unchanged. We can
therefore use the monoidal, left-autonomous and pivotal structure of
$\sym{M}_H$ to equip $\sym{M}_{(\sym{M}^H,U)}$ with the structure of
a pivotal category in such a way that the functor $\Phi$ becomes pivotal and
strict monoidal. Finally, $\Phi$ is $k$-linear. We therefore get the following theorem:

\begin{theorem}
\label{thm_repmoncat}
Let $\sym{C}$ be a finitely semisimple $k$-linear additive spherical category,
$k=\End(\one)$ be a field, $\omega\colon\sym{C}\to\Vect_k$ be the long forgetful
functor, $H=\coend(\sym{C},\omega)$ and $U\colon\sym{M}^H\to\Vect_k$ be the
forgetful functor. Then
\begin{equation}
\sym{M}_H\simeq\sym{M}_{(\sym{M}^H,U)}
\end{equation}
are equivalent as $k$-linear additive pivotal categories.
\end{theorem}

\subsection{Self-duality of spherical categories}

Before we can combine all our results in order to prove
Corollary~\ref{cor_intro3}, we need to relate $\sym{M}^H$ with $\sym{M}_{\hat H}$.

\begin{proposition}
\label{prop_duality}
Under the assumptions of Theorem~\ref{thm_repmoncat}, the categories
\begin{equation}
\sym{M}^H\simeq\sym{M}_{\hat H}
\end{equation}
are equivalent as $k$-linear additive pivotal categories.
\end{proposition}

\begin{proof}
We use the fact that the pairing $\left<-;\hat-\right>$ and the canonical
element $G$ of Proposition~\ref{prop_pairing} satisfy the triangle
identities. Then, $\sym{M}^H\cong\sym{M}_{\hat H}$ are isomorphic as
$k$-linear additive categories using the functor that turns every
right $H$-comodule $\beta_V\colon V\to V\otimes H$ into a right $\hat
H$-module
\begin{equation}
\gamma_V\colon V\otimes\hat H\to V,\quad
v\otimes\hat h\mapsto v_V\,\bigl<v_H;\hat h\bigr>,
\end{equation}
and the functor that turns every right $\hat H$-module $\gamma_V\colon
V\otimes\hat H\to V$ into a right $H$-comodule
\begin{equation}
\beta_V = (\gamma_V\otimes\id_H)\circ\alpha^{-1}_{V,\hat H,H}
\circ(\id_V\otimes G)\circ\rho_V^{-1}\colon
V\to V\otimes H.
\end{equation}
The categories $\sym{M}^H\simeq\sym{M}_{\hat H}$ are also equivalent as
pivotal categories, because the pairing also satisfies
\begin{equation}
\left<h;\hat m\right>=w(h)
\end{equation}
for all $h\in H$. Here $w$ denotes the pivotal form~\eqref{eq_pivotalform} of
$H$ and $\hat m$ the analogue of the pivotal element~\eqref{eq_pivotalelement} in
$\hat H$.
\end{proof}

\begin{proof}[Proof of Corollary~\ref{cor_intro3}]
Under the assumptions of Theorem~\ref{thm_repmoncat}, we have the following
equivalences of $k$-linear additive pivotal categories
\begin{equation}
\label{eq_chain}
\sym{C}\simeq\sym{M}^H\simeq\sym{M}_{\hat H}\simeq\sym{M}_H
\simeq\sym{M}_{(\sym{M}^H,U)}\simeq\sym{M}_{(\sym{C},\omega)}.
\end{equation}
The first one is from Theorem~\ref{thm_intro1}, the second is from
Proposition~\ref{prop_duality}, the third follows from the
isomorphism~\eqref{eq_isohhat}, the fourth is from
Theorem~\ref{thm_repmoncat} and the fifth from
Theorem~\ref{thm_monequiv}. Since $\sym{C}$ is spherical, all other categories
are spherical, too, and the equivalence is an equivalence of spherical
categories. In addition, $H$ and $\hat H$ are both spherical and cospherical.
\end{proof}

\begin{remark}
In order to see how the objects $X\in|\sym{C}|$ give rise to right
$(\sym{C},\omega)$-modules in~\eqref{eq_chain}, we combine the relevant
functors as follows. For each simple $X\in|\sym{C}|$, there is a right $H$-comodule
$\omega(X)=\Hom(\tilde V,\tilde V\otimes X)$, given by
\begin{eqnarray}
\label{eq_coaction}
\beta_X\colon\omega(X)&\to    &\omega(X)\otimes H,\\
e^{(X)}_j&\mapsto&\sum_\ell e^{(X)}_\ell\otimes{[e^\ell_{(X)}|e^{(X)}_j]}_X.\nn
\end{eqnarray}
It can be turned into a right $H$-module
\begin{eqnarray}
\label{eq_action}
\gamma_{\omega(X)}\colon\omega(X)\otimes H&\to&\omega(X),\\
e^{(X)}_j\otimes{[e^p_{(Y)}|e^{(Y)}_\ell]}_Y&\mapsto&
\sum_q e^{(X)}_q\,\left<e^{(X)}_{jq};\hat e^{(Y)}_{\ell p}\right>.\nn
\end{eqnarray}
Here the $e^{(X)}_{jq}$ are as in~\eqref{eq_defeml}. Note that the above
expression involves the non-canonical isomorphism~\eqref{eq_isohhat}, putting
hats on all expressions and replacing $\tilde V$ by $\tilde V^\ast$
everywhere. Finally, $\omega(X)$ forms a right $(\sym{C},\omega)$-module with
\begin{eqnarray}
\label{eq_braiding}
{(c_{\omega(X)})}_Y\colon \omega(X)\otimes\omega(Y)&\to&\omega(Y)\otimes\omega(X),\\
e^{(X)}_j\otimes e^{(Y)}_\ell&\mapsto&\sum_{p,q} e^{(Y)}_p\otimes e^{(X)}_q\,
\left<e^{(X)}_{jq};\hat e^{(Y)}_{\ell p}\right>\nn
\end{eqnarray}
for all simple $Y\in|\sym{C}|$. Note that in our construction of
${(c_{\omega(X)})}_Y$, $\omega(X)$ forms a right $H$-module in $\Vect_k$, but
not in general a right $H$-module in $\sym{M}^H$ and that the tensor product
$\omega(X)\otimes\omega(Y)$ above is in $\Vect_k$ as opposed to $\sym{M}^H$.
\end{remark}

\subsection{The modular case}

Every modular category is finitely semisimple $k$-linear additive spherical,
and so if $\sym{C}$ is modular, we know that~\eqref{eq_chain} is an
equivalence of $k$-linear additive spherical categories. In order to show that
it is actually an equivalence of additive ribbon categories, we can proceed as
follows. First, we explain the additional structure.
\begin{myenumerate}
\item
$\sym{C}\simeq\sym{M}^H$ are equivalent as additive ribbon categories
by~\cite{Pf07}. $H=\coend(\sym{C},\omega)$ is now a finite-dimensional
split-cosemisimple coribbon WHA that is weakly cofactorizable and for which
$H_s\cap H_t\cong k$. The new structure is the coquasitriangular structure
and the universal ribbon form of $H$.
\item
In order to obtain an equivalence $\sym{M}^H\simeq\sym{M}_{\hat H}$ of
additive ribbon categories, we define the notion of a quasitriangular and
ribbon WHA in such a way that the pairing~\eqref{eq_dualpair} relates the
coquasitriangular and coribbon structure of $H$ with the quasitriangular and
ribbon structure of $\hat H$.
\item
The equivalence $\sym{M}_{\hat H}\simeq\sym{M}_H$ is automatically an
equivalence of additive ribbon categories.
\item
For the equivalences
$\sym{M}_H\simeq\sym{M}_{(\sym{M}^H,U)}\simeq\sym{M}_{(\sym{C},
\omega)}$, one defines the braiding and ribbon twist of
$\sym{M}_{(\sym{M}^H,U)}$ and $\sym{M}_{(\sym{C},\omega)}$ accordingly.
\end{myenumerate}
The additional property, namely weak cofactorizability of $H$, weak
factorizability (to be defined accordingly) of $\hat H$, and the
non-degeneracy condition of $\sym{M}_{\hat H}$, $\sym{M}_H$,
$\sym{M}_{(\sym{M}^H,U)}$ and $\sym{M}_{(\sym{C},\omega)}$ then follows from
the non-degeneracy of $\sym{C}$.

\section{Example}
\label{sect_example}

In this section, we specialize the key expressions used in the present article
to the case of the modular category $\sym{C}$ associated with the quantum
group $U_q(\ssl_2)$, $q$ a root of unity. We use the diagrammatic description
of~\cite{KaLi94} and precisely follow their notation.

Let $r\in\{2,3,4,\ldots\}$ and $A$ be a primitive $4r$-th root of unity,
$q=A^2$. For simplicity, we work over the complex numbers $k=\C$. The
morphisms of $\sym{C}$ are represented by plane projections of oriented framed
tangles, drawn in blackboard framing. The coherence theorem for ribbon
categories~\cite{ReTu90} ensures that each diagram defines a morphism of
$\sym{C}$. Since $\sym{C}$ is $k$-linear, we can take formal linear
combinations of diagrams with coefficients in $k$. All our diagrams are read
from top to bottom.

The braiding of $\sym{C}$ is such that a crossing in our plane projections can
be resolved using the recursion relation for the Kauffman bracket
\begin{equation}
\begin{xy}
(-4,4);(4,-4) **\dir{-};
(-4,-4);(-1,-1) **\dir{-};
(1,1);(4,4) **\dir{-}
\end{xy}
= A\,
\begin{xy}
(-4,4);(4,4) **\crv{(-2,2)&(2,2)};
(-4,-4);(4,-4) **\crv{(-2,-2)&(2,-2)};
\end{xy}
+ A^{-1}\,
\begin{xy}
(4,-4);(4,4) **\crv{(2,-2)&(2,2)};
(-4,-4);(-4,4) **\crv{(-2,-2)&(-2,2)};
\end{xy}\,,\qquad\qquad
\begin{xy}
(0,0)*\xycircle<3mm,3mm>{-};
\end{xy}
= -(q+q^{-1}),
\end{equation}
ignoring the orientations for now. The Jones--Wenzl idempotents $P_n$,
$n\in\N_0$, are formal linear combinations of planar $(n,n)$-tangles that can
be defined recursively by
\begin{equation}
\begin{xy}
(0,0)*+{\hbox to 10mm{\strut\hfill$P_1$\hfill}}*\frm{-}!U="t" !D="b";
"t";(0, 12) **\dir{-};
"b";(0,-12) **\dir{-};
\end{xy}
:=
\begin{xy}
(0,-12);(0,12) **\dir{-};
\end{xy}\,,\qquad\qquad
\begin{xy}
(0,0)*+{\hbox to 10mm{\strut\hfill$P_{n+1}$\hfill}}*\frm{-}!U="t" !D="b";
"t"+<-4mm,0mm>;(-4, 12) **\dir{-};
"b"+<-4mm,0mm>;(-4,-12) **\dir{-};
"t"+<-3mm,0mm>;(-3, 12) **\dir{-};
"b"+<-3mm,0mm>;(-3,-12) **\dir{-};
"t"+<-2mm,0mm>;(-2, 12) **\dir{-}; ?(.5)+<3mm,0mm>*{\cdots};
"b"+<-2mm,0mm>;(-2,-12) **\dir{-}; ?(.5)+<3mm,0mm>*{\cdots};
"t"+< 4mm,0mm>;( 4, 12) **\dir{-};
"b"+< 4mm,0mm>;( 4,-12) **\dir{-};
\end{xy}
:=
\begin{xy}
(0,0)*+{\hbox to 10mm{\strut\hfill$P_n$\hfill}}*\frm{-}!U="t" !D="b";
"t"+<-4mm,0mm>;(-4, 12) **\dir{-};
"b"+<-4mm,0mm>;(-4,-12) **\dir{-};
"t"+<-3mm,0mm>;(-3, 12) **\dir{-};
"b"+<-3mm,0mm>;(-3,-12) **\dir{-};
"t"+<-2mm,0mm>;(-2, 12) **\dir{-}; ?(.5)+<3mm,0mm>*{\cdots};
"b"+<-2mm,0mm>;(-2,-12) **\dir{-}; ?(.5)+<3mm,0mm>*{\cdots};
"t"+< 4mm,0mm>;( 4, 12) **\dir{-};
"b"+< 4mm,0mm>;( 4,-12) **\dir{-};
(8,-12);(8,12) **\dir{-};
\end{xy}
+\frac{\q{n+1}}{\q{n+2}}\,
\begin{xy}
(0, 6)*+{\hbox to 10mm{\strut\hfill$P_n$\hfill}}*\frm{-}!U="t1" !D="b1";
(0,-6)*+{\hbox to 10mm{\strut\hfill$P_n$\hfill}}*\frm{-}!U="t2" !D="b2";
"t1"+<-4mm,0mm>;(-4, 12)        **\dir{-};
"b1"+<-4mm,0mm>;"t2"+<-4mm,0mm> **\dir{-};
"b2"+<-4mm,0mm>;(-4,-12)        **\dir{-};
"t1"+<-3mm,0mm>;(-3, 12)        **\dir{-}; ?(.5)+<2.5mm,0mm>*{\cdots};
"b1"+<-3mm,0mm>;"t2"+<-3mm,0mm> **\dir{-}; ?(.5)+<2.5mm,0mm>*{\cdots};
"b2"+<-3mm,0mm>;(-3,-12)        **\dir{-}; ?(.5)+<2.5mm,0mm>*{\cdots};
"t1"+< 3mm,0mm>;( 3, 12)        **\dir{-};
"b1"+< 3mm,0mm>;"t2"+< 3mm,0mm> **\dir{-};
"b2"+< 3mm,0mm>;( 3,-12)        **\dir{-};
"t1"+< 4mm,0mm>;( 4, 12)        **\dir{-};
"b2"+< 4mm,0mm>;( 4,-12)        **\dir{-};
"b1"+< 4mm,0mm>;( 8, 12) **\crv{"b1"+<4mm,-2mm>&"b1"+<8mm,-2mm>&"b1"+<8mm,0mm>};
"t2"+< 4mm,0mm>;( 8,-12) **\crv{"t2"+<4mm, 2mm>&"t2"+<8mm, 2mm>&"t2"+<8mm,0mm>};
\end{xy}\,.
\end{equation}
where $\q{n}=(q^n-q^{-n})/(q-q^{-1})$, $n\in\Z$, are the quantum
integers. The isomorphism classes of simple objects of $\sym{C}$ are indexed
by the set $I=\{0,1,\ldots,r-2\}$. The identity morphism of the object $n\in
I$ is the identity $(n,n)$-tangle with the idempotent $P_n$ inserted
somewhere (anywhere). As a shortcut, we write a single line labeled by $n$,
\begin{equation}
\begin{xy}
(0,-8);(0,8) **\dir{-} ?(.8)+<2mm,0mm>*{n};
\end{xy}
:=
\begin{xy}
(0,0)*+{\hbox to 10mm{\strut\hfill$P_n$\hfill}}*\frm{-}!U="t" !D="b";
"t"+<-4mm,0mm>;(-4, 8) **\dir{-};
"b"+<-4mm,0mm>;(-4,-8) **\dir{-};
"t"+<-3mm,0mm>;(-3, 8) **\dir{-};
"b"+<-3mm,0mm>;(-3,-8) **\dir{-};
"t"+<-2mm,0mm>;(-2, 8) **\dir{-}; ?(.5)+<3mm,0mm>*{\cdots};
"b"+<-2mm,0mm>;(-2,-8) **\dir{-}; ?(.5)+<3mm,0mm>*{\cdots};
"t"+< 4mm,0mm>;( 4, 8) **\dir{-};
"b"+< 4mm,0mm>;( 4,-8) **\dir{-};
\end{xy}\,.
\end{equation}
The object indexed by $0\in I$ is the monoidal unit and can be made invisible
in our diagrams thanks to the coherence theorem. The categorical dimension of
the simple objects is given by
\begin{equation}
\Delta_n :=
\begin{xy}
(0,0)*\xycircle<3mm,3mm>{-};
(3.5,3.5)*{n}
\end{xy}
= (-1)^n\q{n+1},
\end{equation}
which is non-zero for all $n\in I$.

Two special features of $U_q(\ssl_2)$ are exploited. First, the simple objects
are isomorphic to their duals, and the choice of representatives $V_j$, $j\in
I$, of the simple objects is such that ${(V_j)}^\ast=V_j$ are equal rather
than merely isomorphic. This allows us to omit any arrows from the diagrams
that would indicate the orientation of the ribbon tangle.

Second, there are no higher multiplicities, \ie\ for all $a,b,c\in I$,
we have $\dim_k\Hom(V_a\otimes V_b,V_c)\in\{0,1\}$. More precisely,
$\Hom(V_a\otimes V_b,V_c)\cong k$ if and only if the triple $(a,b,c)$
is \emph{admissible}. Otherwise, $\Hom(V_a\otimes V_b,V_c)=\{0\}$.

\begin{definition}
A triple $(a,b,c)\in I^3$ is called \emph{admissible} if the following
conditions hold.
\begin{myenumerate}
\item
$a+b+c\equiv 0$ mod $2$ (\emph{parity}),
\item
$a+b-c\geq 0$ and $b+c-a\geq 0$ and $c+a-b\geq 0$ (\emph{quantum triangle inequality}),
\item
$a+b+c\leq 2r-4$ (\emph{non-negligibility}).
\end{myenumerate}
\end{definition}

A special choice of basis vector of $\Hom(V_a,V_b\otimes V_c)$ is
denoted by a trivalent vertex:
\begin{equation}
\begin{xy}
(0,0)*\dir{*};
(0,0);(0,8) **\dir{-}; ?(.7)+<3mm,0mm>*{a};
(0,0);(-4,-8) **\crv{(-4,-4)&(-4,-6)}; ?(.7)+<-3mm,0mm>*{b};
(0,0);( 4,-8) **\crv{( 4,-4)&( 4,-6)}; ?(.7)+< 3mm,0mm>*{c};
\end{xy}
:=
\begin{xy}
(-2,3);( 2,3) **\dir{-};
( 2,3);( 2,5) **\dir{-};
( 2,5);(-2,5) **\dir{-};
(-2,5);(-2,3) **\dir{-};
( 0,5);( 0,8) **\dir{-}; ?(.5)+<3mm,0mm>*{a};
(-6,-5);(-2,-5) **\dir{-};
(-2,-5);(-2,-3) **\dir{-};
(-2,-3);(-6,-3) **\dir{-};
(-6,-3);(-6,-5) **\dir{-};
(-4,-5);(-4,-8) **\dir{-}; ?(.5)+<-3mm,0mm>*{b};
( 2,-5);( 6,-5) **\dir{-};
( 6,-5);( 6,-3) **\dir{-};
( 6,-3);( 2,-3) **\dir{-};
( 2,-3);( 2,-5) **\dir{-};
( 4,-5);( 4,-8) **\dir{-}; ?(.5)+<3mm,0mm>*{c};
(-1, 3);(-5,-3) **\crv{(-1,0)&(-5,0)}; ?(.5)+<-3mm,0mm>*{i};
( 1, 3);( 5,-3) **\crv{( 1,0)&( 5,0)}; ?(.5)+< 3mm,0mm>*{j};
(-3,-3);( 3,-3) **\crv{(-3,0)&( 3,0)}; ?(.5)+<0mm,-2mm>*{k};
\end{xy}\,,
\end{equation}
where $i=(a+b-c)/2$, $j=(a+c-b)/2$ and $k=(b+c-a)/2$. If we draw such a
diagram for a triple $(a,b,c)\in I^3$ that is not admissible, then by convention,
we multiply the entire diagram by zero. We also need the theta graph
\begin{equation}
\theta(a,b,c) :=
\begin{xy}
(-5,0)*\dir{*};
(5,0)*\dir{*};
(-5,0);(5,0) **\crv{(-5,5)&(5,5)}; ?(.5)+<0mm,2mm>*{a};
(-5,0);(5,0) **\dir{-}; ?(.5)+<0mm,2mm>*{b};
(-5,0);(5,0) **\crv{(-5,-5)&(5,-5)}; ?(.5)+<0mm,-2mm>*{c};
\end{xy}\,,
\end{equation}
which is non-zero for all admissible triples $(a,b,c)$. When we compose the
morphisms associated with such diagrams, the composition is zero unless the
labels at the open ends of the tangles match, \ie\ putting
\begin{equation}
\begin{xy}
(0,0)*\dir{*};
(0,0);(-4, 4) **\dir{-}; ?(.8)+<-3mm,0mm>*{r};
(0,0);( 4, 4) **\dir{-}; ?(.8)+< 3mm,0mm>*{j};
(0,0);( 0,-6) **\dir{-}; ?(.8)+< 3mm,0mm>*{s};
\end{xy}\qquad\mbox{below}\qquad
\begin{xy}
(0,0)*\dir{*};
(0,0);(-4,-4) **\dir{-}; ?(.8)+<-3mm,0mm>*{q};
(0,0);( 4,-4) **\dir{-}; ?(.8)+< 3mm,0mm>*{k};
(0,0);( 0, 6) **\dir{-}; ?(.8)+< 3mm,0mm>*{p};
\end{xy}\qquad\mbox{gives}\qquad
\delta_{qr}\delta_{kj}\,\,
\begin{xy}
(0, 4)*\dir{*};
(0,-4)*\dir{*};
(0,-4);(0, 4) **\crv{(-4,-4)&(-4, 4)}; ?(.5)+<-2mm,0mm>*{q};
(0,-4);(0, 4) **\crv{( 4,-4)&( 4, 4)}; ?(.5)+< 2mm,0mm>*{k};
(0, 4);(0, 8) **\dir{-}; ?(.8)+< 3mm,0mm>*{p};
(0,-4);(0,-8) **\dir{-}; ?(.8)+< 3mm,0mm>*{s};
\end{xy}
\end{equation}
We use the following pair of dual bases of $\omega(V_j)=\Hom(\tilde V,\tilde
V\otimes V_j)$ and ${\omega(V_j)}^\ast=\Hom(\tilde V\otimes V_j,\tilde V)$,
$j\in I$,
\begin{equation}
e^{(V_j)}_{pq} =
\begin{xy}
(0,0)*\dir{*};
(0,0);(0,6) **\dir{-}; ?(.7)+<3mm,0mm>*{p};
(0,0);(-4,-4) **\dir{-}; ?(.7)+<-3mm,0mm>*{q};
(0,0);( 4,-4) **\dir{-}; ?(.7)+< 3mm,0mm>*{j};
\end{xy}\qquad\mbox{and}\qquad
e_{(V_j)}^{qp} = \frac{\Delta_p}{\theta(q,p,j)}\,
\begin{xy}
(0,0)*\dir{*};
(0,0);( 0,-6) **\dir{-}; ?(.7)+<3mm,0mm>*{p};
(0,0);(-4, 4) **\dir{-}; ?(.7)+<-3mm,0mm>*{q};
(0,0);( 4, 4) **\dir{-}; ?(.7)+< 3mm,0mm>*{j};
\end{xy}\,,
\end{equation}
where $p,q\in I$ are such that $(p,q,j)$ is admissible. Then the
basis~\eqref{eq_defeml} of $H=\coend(\sym{C},\omega)$ is given by
\begin{equation}
e^{(V_j)}_{pq,rs} = \frac{1}{\theta(j,r,s)}
\begin{xy}
(0,4)*\dir{*};
(0,-4)*\dir{*};
(0,-4);( 0, 4) **\dir{-}; ?(.5)+< 3mm,0mm>*{j};
(0, 4);(-4, 8) **\dir{-}; ?(.8)+<-3mm,0mm>*{q};
(0, 4);( 4, 8) **\dir{-}; ?(.8)+< 3mm,0mm>*{p};
(0,-4);(-4,-8) **\dir{-}; ?(.8)+<-3mm,0mm>*{r};
(0,-4);( 4,-8) **\dir{-}; ?(.8)+< 3mm,0mm>*{s};
\end{xy}\,.
\end{equation}
The quantum $6j$-symbol is defined as
\begin{equation}
\left\{\begin{matrix}a&b&i\\ c&d&j\end{matrix}\right\}_q\,
:=\frac{\Delta_i}{\theta(a,d,i)\theta(b,c,i)}\,
\begin{xy}
( 0, 6)*\dir{*};
(-6, 0)*\dir{*};
( 6, 0)*\dir{*};
( 0,-6)*\dir{*};
(-6, 0);( 6, 0) **\dir{-}; ?(.5)+< 0mm,-2mm>*{j};
(-6, 0);( 0,-6) **\dir{-}; ?(.5)+<-2mm,-2mm>*{a};
(-6, 0);( 0, 6) **\dir{-}; ?(.5)+<-2mm, 2mm>*{b};
( 0, 6);( 6, 0) **\dir{-}; ?(.5)+< 2mm, 2mm>*{c};
( 0,-6);( 6, 0) **\dir{-}; ?(.5)+< 2mm,-2mm>*{d};
(12, 6);(12,-6) **\dir{-}; ?(.5)+< 3mm, 0mm>*{i};
( 0, 6);(12, 6) **\crv{(0,12)&(12,12)};
( 0,-6);(12,-6) **\crv{(0,-12)&(12,-12)};
\end{xy}\,.
\end{equation}
It is used in the recoupling identity,
\begin{equation}
\begin{xy}
(-4, 0)*\dir{*};
( 4, 0)*\dir{*};
(-4, 0);( 4, 0) **\dir{-}; ?(.5)+<0mm, 3mm>*{j};
( 4, 0);( 8, 4) **\dir{-}; ?(.8)+<0mm, 3mm>*{c};
( 4, 0);( 8,-4) **\dir{-}; ?(.8)+<0mm,-3mm>*{d};
(-4, 0);(-8,-4) **\dir{-}; ?(.8)+<0mm,-3mm>*{a};
(-4, 0);(-8, 4) **\dir{-}; ?(.8)+<0mm, 3mm>*{b};
\end{xy} = \sum_i\,
\left\{\begin{matrix}a&b&i\\ c&d&j\end{matrix}\right\}_q\,
\begin{xy}
(0,4)*\dir{*};
(0,-4)*\dir{*};
(0,-4);( 0, 4) **\dir{-}; ?(.5)+< 3mm,0mm>*{i};
(0, 4);(-4, 8) **\dir{-}; ?(.8)+<-3mm,0mm>*{b};
(0, 4);( 4, 8) **\dir{-}; ?(.8)+< 3mm,0mm>*{c};
(0,-4);(-4,-8) **\dir{-}; ?(.8)+<-3mm,0mm>*{a};
(0,-4);( 4,-8) **\dir{-}; ?(.8)+< 3mm,0mm>*{d};
\end{xy}\,.
\end{equation}
Diagrammatically, the WHA structure of $H$ is:
\begin{equation}
\eta(1) = \sum_{p,q}\frac{1}{\Delta_q}\,
\begin{xy}
(-4, 4);(4, 4) **\crv{(-4,0)&(4,0)}; ?(.5)+<0mm, 2mm>*{p};
(-4,-4);(4,-4) **\crv{(-4,0)&(4,0)}; ?(.5)+<0mm,-2mm>*{q};
\end{xy}\,,
\end{equation}
\begin{gather}
\mu\Biggl(\,
\begin{xy}
(0,4)*\dir{*};
(0,-4)*\dir{*};
(0,-4);( 0, 4) **\dir{-}; ?(.5)+< 3mm,0mm>*{j};
(0, 4);(-4, 8) **\dir{-}; ?(.8)+<-3mm,0mm>*{q};
(0, 4);( 4, 8) **\dir{-}; ?(.8)+< 3mm,0mm>*{p};
(0,-4);(-4,-8) **\dir{-}; ?(.8)+<-3mm,0mm>*{r};
(0,-4);( 4,-8) **\dir{-}; ?(.8)+< 3mm,0mm>*{s};
\end{xy}\,\otimes\,
\begin{xy}
(0,4)*\dir{*};
(0,-4)*\dir{*};
(0,-4);( 0, 4) **\dir{-}; ?(.5)+< 3mm,0mm>*{\ell};
(0, 4);(-4, 8) **\dir{-}; ?(.8)+<-3mm,0mm>*{b};
(0, 4);( 4, 8) **\dir{-}; ?(.8)+< 3mm,0mm>*{a};
(0,-4);(-4,-8) **\dir{-}; ?(.8)+<-3mm,0mm>*{c};
(0,-4);( 4,-8) **\dir{-}; ?(.8)+< 3mm,0mm>*{d};
\end{xy}\,\Biggr) = \delta_{pb}\delta_{sc}\,\theta(j,b,q)\,\frac{\Delta_c}{\Delta_q}\,\sum_{u\in I}\,
\left\{\begin{matrix}r&j&u\\ \ell&d&c\end{matrix}\right\}_q\,
\left\{\begin{matrix}u&j&q\\ b&a&\ell\end{matrix}\right\}_q\,
\begin{xy}
(0,4)*\dir{*};
(0,-4)*\dir{*};
(0,-4);( 0, 4) **\dir{-}; ?(.5)+< 3mm,0mm>*{u};
(0, 4);(-4, 8) **\dir{-}; ?(.8)+<-3mm,0mm>*{q};
(0, 4);( 4, 8) **\dir{-}; ?(.8)+< 3mm,0mm>*{a};
(0,-4);(-4,-8) **\dir{-}; ?(.8)+<-3mm,0mm>*{r};
(0,-4);( 4,-8) **\dir{-}; ?(.8)+< 3mm,0mm>*{d};
\end{xy}\,,\\
\Delta\Biggl(\,
\begin{xy}
(0,4)*\dir{*};
(0,-4)*\dir{*};
(0,-4);( 0, 4) **\dir{-}; ?(.5)+< 3mm,0mm>*{j};
(0, 4);(-4, 8) **\dir{-}; ?(.8)+<-3mm,0mm>*{q};
(0, 4);( 4, 8) **\dir{-}; ?(.8)+< 3mm,0mm>*{p};
(0,-4);(-4,-8) **\dir{-}; ?(.8)+<-3mm,0mm>*{r};
(0,-4);( 4,-8) **\dir{-}; ?(.8)+< 3mm,0mm>*{s};
\end{xy}\,\Biggr) = \sum_{t,u}\frac{1}{\theta(j,u,t)}
\begin{xy}
(0,4)*\dir{*};
(0,-4)*\dir{*};
(0,-4);( 0, 4) **\dir{-}; ?(.5)+< 3mm,0mm>*{j};
(0, 4);(-4, 8) **\dir{-}; ?(.8)+<-3mm,0mm>*{u};
(0, 4);( 4, 8) **\dir{-}; ?(.8)+< 3mm,0mm>*{t};
(0,-4);(-4,-8) **\dir{-}; ?(.8)+<-3mm,0mm>*{r};
(0,-4);( 4,-8) **\dir{-}; ?(.8)+< 3mm,0mm>*{s};
\end{xy}\,\otimes\,
\begin{xy}
(0,4)*\dir{*};
(0,-4)*\dir{*};
(0,-4);( 0, 4) **\dir{-}; ?(.5)+< 3mm,0mm>*{j};
(0, 4);(-4, 8) **\dir{-}; ?(.8)+<-3mm,0mm>*{q};
(0, 4);( 4, 8) **\dir{-}; ?(.8)+< 3mm,0mm>*{p};
(0,-4);(-4,-8) **\dir{-}; ?(.8)+<-3mm,0mm>*{u};
(0,-4);( 4,-8) **\dir{-}; ?(.8)+< 3mm,0mm>*{t};
\end{xy}\,,\\
\epsilon\Biggl(\,
\begin{xy}
(0,4)*\dir{*};
(0,-4)*\dir{*};
(0,-4);( 0, 4) **\dir{-}; ?(.5)+< 3mm,0mm>*{j};
(0, 4);(-4, 8) **\dir{-}; ?(.8)+<-3mm,0mm>*{q};
(0, 4);( 4, 8) **\dir{-}; ?(.8)+< 3mm,0mm>*{p};
(0,-4);(-4,-8) **\dir{-}; ?(.8)+<-3mm,0mm>*{r};
(0,-4);( 4,-8) **\dir{-}; ?(.8)+< 3mm,0mm>*{s};
\end{xy}\,\Biggr) = \delta_{qr}\delta_{ps}\theta(j,r,s),\\
S\Biggl(\,
\begin{xy}
(0,4)*\dir{*};
(0,-4)*\dir{*};
(0,-4);( 0, 4) **\dir{-}; ?(.5)+< 3mm,0mm>*{j};
(0, 4);(-4, 8) **\dir{-}; ?(.8)+<-3mm,0mm>*{q};
(0, 4);( 4, 8) **\dir{-}; ?(.8)+< 3mm,0mm>*{p};
(0,-4);(-4,-8) **\dir{-}; ?(.8)+<-3mm,0mm>*{r};
(0,-4);( 4,-8) **\dir{-}; ?(.8)+< 3mm,0mm>*{s};
\end{xy}\,\Biggr) = \frac{\Delta_q\,\theta(j,r,s)}{\Delta_p\,\theta(j,p,q)}\,
\begin{xy}
(0,4)*\dir{*};
(0,-4)*\dir{*};
(0,-4);( 0, 4) **\dir{-}; ?(.5)+< 3mm,0mm>*{j};
(0, 4);(-4, 8) **\dir{-}; ?(.8)+<-3mm,0mm>*{s};
(0, 4);( 4, 8) **\dir{-}; ?(.8)+< 3mm,0mm>*{r};
(0,-4);(-4,-8) **\dir{-}; ?(.8)+<-3mm,0mm>*{p};
(0,-4);( 4,-8) **\dir{-}; ?(.8)+< 3mm,0mm>*{q};
\end{xy}\,.
\end{gather}
The pivotal form of $H$ is given by
\begin{equation}
w\Biggl(\,
\begin{xy}
(0,4)*\dir{*};
(0,-4)*\dir{*};
(0,-4);( 0, 4) **\dir{-}; ?(.5)+< 3mm,0mm>*{j};
(0, 4);(-4, 8) **\dir{-}; ?(.8)+<-3mm,0mm>*{q};
(0, 4);( 4, 8) **\dir{-}; ?(.8)+< 3mm,0mm>*{p};
(0,-4);(-4,-8) **\dir{-}; ?(.8)+<-3mm,0mm>*{r};
(0,-4);( 4,-8) **\dir{-}; ?(.8)+< 3mm,0mm>*{s};
\end{xy}\,\Biggr) = \delta_{ps}\delta_{qr}\,\theta(j,r,s)\,\frac{\Delta_p}{\Delta_q}.
\end{equation}
Since $\tilde V=\tilde V^\ast$, the dual WHA $\hat H$ has precisely the same
description as $H$. The pairing $\left<-|\hat-\right>$ of~\eqref{eq_dualpair}
reads:
\begin{equation}
\Biggl<\,
\begin{xy}
(0,4)*\dir{*};
(0,-4)*\dir{*};
(0,-4);( 0, 4) **\dir{-}; ?(.5)+< 3mm,0mm>*{j};
(0, 4);(-4, 8) **\dir{-}; ?(.8)+<-3mm,0mm>*{b};
(0, 4);( 4, 8) **\dir{-}; ?(.8)+< 3mm,0mm>*{a};
(0,-4);(-4,-8) **\dir{-}; ?(.8)+<-3mm,0mm>*{c};
(0,-4);( 4,-8) **\dir{-}; ?(.8)+< 3mm,0mm>*{d};
\end{xy}\,;\,
\begin{xy}
(0,4)*\dir{*};
(0,-4)*\dir{*};
(0,-4);( 0, 4) **\dir{-}; ?(.5)+< 3mm,0mm>*{\ell};
(0, 4);(-4, 8) **\dir{-}; ?(.8)+<-3mm,0mm>*{q};
(0, 4);( 4, 8) **\dir{-}; ?(.8)+< 3mm,0mm>*{p};
(0,-4);(-4,-8) **\dir{-}; ?(.8)+<-3mm,0mm>*{r};
(0,-4);( 4,-8) **\dir{-}; ?(.8)+< 3mm,0mm>*{s};
\end{xy}\,\Biggr> = \delta_{pa}\delta_{qd}\delta_{rc}\delta_{sb}\,\theta(j,c,d)\theta(\ell,r,s)\,
\left\{\begin{matrix}\ell&d&c\\ j&b&a\end{matrix}\right\}_q.
\end{equation}
The other pairing $\left<\hat-|-\right>$ is different. The canonical element
$G(1)\in\hat H\otimes H$ of~\eqref{eq_canelement} is given by
\begin{equation}
G(1)=\sum_{j,a,b,c,d}\frac{\Delta_j}{\Delta_c\,\theta(j,a,b)\theta(j,c,d)}\,
\begin{xy}
(-4, 0)*\dir{*};
( 4, 0)*\dir{*};
(-4, 0);( 4, 0) **\dir{-}; ?(.5)+<0mm, 3mm>*{j};
( 4, 0);( 8, 4) **\dir{-}; ?(.8)+<0mm, 3mm>*{a};
( 4, 0);( 8,-4) **\dir{-}; ?(.8)+<0mm,-3mm>*{b};
(-4, 0);(-8,-4) **\dir{-}; ?(.8)+<0mm,-3mm>*{c};
(-4, 0);(-8, 4) **\dir{-}; ?(.8)+<0mm, 3mm>*{d};
\end{xy}\,\otimes\,
\begin{xy}
(0,4)*\dir{*};
(0,-4)*\dir{*};
(0,-4);( 0, 4) **\dir{-}; ?(.5)+< 3mm,0mm>*{j};
(0, 4);(-4, 8) **\dir{-}; ?(.8)+<-3mm,0mm>*{b};
(0, 4);( 4, 8) **\dir{-}; ?(.8)+< 3mm,0mm>*{a};
(0,-4);(-4,-8) **\dir{-}; ?(.8)+<-3mm,0mm>*{c};
(0,-4);( 4,-8) **\dir{-}; ?(.8)+< 3mm,0mm>*{d};
\end{xy}\,,
\end{equation}
from which we can read off a pair of dual bases of $\hat H$ and $H$ with
respect to $\left<-|\hat-\right>$. The coaction~\eqref{eq_coaction} of $H$ on
$\omega(V_j)$ and the action~\eqref{eq_action} of $H$ on $\omega(V_j)$ are as
follows:
\begin{equation}
\beta_{V_j}\Biggl(\,
\begin{xy}
(0,0)*\dir{*};
(0,0);(0,6) **\dir{-}; ?(.7)+<3mm,0mm>*{p};
(0,0);(-4,-4) **\dir{-}; ?(.7)+<-3mm,0mm>*{q};
(0,0);( 4,-4) **\dir{-}; ?(.7)+< 3mm,0mm>*{j};
\end{xy}\,\Biggr) = \sum_{r,s}\frac{1}{\theta(j,r,s)}\,
\begin{xy}
(0,0)*\dir{*};
(0,0);(0,6) **\dir{-}; ?(.7)+<3mm,0mm>*{s};
(0,0);(-4,-4) **\dir{-}; ?(.7)+<-3mm,0mm>*{r};
(0,0);( 4,-4) **\dir{-}; ?(.7)+< 3mm,0mm>*{j};
\end{xy}\,\otimes\,
\begin{xy}
(0,4)*\dir{*};
(0,-4)*\dir{*};
(0,-4);( 0, 4) **\dir{-}; ?(.5)+< 3mm,0mm>*{j};
(0, 4);(-4, 8) **\dir{-}; ?(.8)+<-3mm,0mm>*{q};
(0, 4);( 4, 8) **\dir{-}; ?(.8)+< 3mm,0mm>*{p};
(0,-4);(-4,-8) **\dir{-}; ?(.8)+<-3mm,0mm>*{r};
(0,-4);( 4,-8) **\dir{-}; ?(.8)+< 3mm,0mm>*{s};
\end{xy}\,,
\end{equation}
\begin{equation}
\gamma_{\omega(V_j)}\Biggl(\,
\begin{xy}
(0,0)*\dir{*};
(0,0);(0,6) **\dir{-}; ?(.7)+<3mm,0mm>*{p};
(0,0);(-4,-4) **\dir{-}; ?(.7)+<-3mm,0mm>*{q};
(0,0);( 4,-4) **\dir{-}; ?(.7)+< 3mm,0mm>*{j};
\end{xy}\,\otimes\,
\begin{xy}
(0,4)*\dir{*};
(0,-4)*\dir{*};
(0,-4);( 0, 4) **\dir{-}; ?(.5)+< 3mm,0mm>*{\ell};
(0, 4);(-4, 8) **\dir{-}; ?(.8)+<-3mm,0mm>*{b};
(0, 4);( 4, 8) **\dir{-}; ?(.8)+< 3mm,0mm>*{a};
(0,-4);(-4,-8) **\dir{-}; ?(.8)+<-3mm,0mm>*{c};
(0,-4);( 4,-8) **\dir{-}; ?(.8)+< 3mm,0mm>*{d};
\end{xy}\,\Biggr) = \delta_{ap}\delta_{dq}\,\theta(\ell,c,d)\,
\left\{\begin{matrix}\ell&b&c\\ j&q&p\end{matrix}\right\}_q\,
\begin{xy}
(0,0)*\dir{*};
(0,0);(0,6) **\dir{-}; ?(.7)+<3mm,0mm>*{b};
(0,0);(-4,-4) **\dir{-}; ?(.7)+<-3mm,0mm>*{c};
(0,0);( 4,-4) **\dir{-}; ?(.7)+< 3mm,0mm>*{j};
\end{xy}\,.
\end{equation}
Finally, we have the structure of $\omega(V_j)$ as a right
$(\sym{C},\omega)$-module,
\begin{equation}
{(c_{\omega(V_j)})}_{V_\ell}\Biggl(\,
\begin{xy}
(0,0)*\dir{*};
(0,0);(0,6) **\dir{-}; ?(.7)+<3mm,0mm>*{p};
(0,0);(-4,-4) **\dir{-}; ?(.7)+<-3mm,0mm>*{q};
(0,0);( 4,-4) **\dir{-}; ?(.7)+< 3mm,0mm>*{j};
\end{xy}\,\otimes\,
\begin{xy}
(0,0)*\dir{*};
(0,0);(0,6) **\dir{-}; ?(.7)+<3mm,0mm>*{r};
(0,0);(-4,-4) **\dir{-}; ?(.7)+<-3mm,0mm>*{s};
(0,0);( 4,-4) **\dir{-}; ?(.7)+< 3mm,0mm>*{\ell};
\end{xy}\,\Biggr) = \delta_{rp}\,\sum_a\,
\left\{\begin{matrix}\ell&s&a\\ j&q&p\end{matrix}\right\}_q\,
\begin{xy}
(0,0)*\dir{*};
(0,0);(0,6) **\dir{-}; ?(.7)+<3mm,0mm>*{q};
(0,0);(-4,-4) **\dir{-}; ?(.7)+<-3mm,0mm>*{a};
(0,0);( 4,-4) **\dir{-}; ?(.7)+< 3mm,0mm>*{\ell};
\end{xy}\,\otimes\,
\begin{xy}
(0,0)*\dir{*};
(0,0);(0,6) **\dir{-}; ?(.7)+<3mm,0mm>*{s};
(0,0);(-4,-4) **\dir{-}; ?(.7)+<-3mm,0mm>*{a};
(0,0);( 4,-4) **\dir{-}; ?(.7)+< 3mm,0mm>*{j};
\end{xy}\,.
\end{equation}

\appendix
\section{Background on tensor categories}

In this appendix, we collect the relevant definitions and properties of
autonomous, pivotal and spherical categories, following
Schauenburg~\cite{Sc92}, Freyd--Yetter~\cite{FrYe92},
Barrett--Westbury~\cite{BaWe99} and Turaev~\cite{Tu94}, and of abelian
categories following MacLane~\cite{Ma73}.

\subsection{Monoidal categories}

\begin{definition}
A \emph{monoidal category}
$(\sym{C},\otimes,\one,\alpha,\lambda,\rho)$ is a category $\sym{C}$
with a bifunctor $\otimes\colon\sym{C}\times\sym{C}\to\sym{C}$
(\emph{tensor product}), an object $\one\in|\sym{C}|$ (\emph{monoidal
unit}) and natural isomorphisms $\alpha_{X,Y,Z}\colon(X\otimes
Y)\otimes Z\to X\otimes(Y\otimes Z)$ (\emph{associator}),
$\lambda_X\colon \one\otimes X\to X$ (\emph{left-unit constraint}) and
$\rho_X\colon X\otimes\one\to X$ (\emph{right-unit constraint}) for
all $X,Y,Z\in|\sym{C}|$, subject to the pentagon axiom
\begin{equation}
\alpha_{X,Y,Z\otimes W}\circ\alpha_{X\otimes Y,Z,W}
= (\id_X\otimes\alpha_{Y,Z,W})\circ\alpha_{X,Y\otimes Z,W}\circ(\alpha_{X,Y,Z}\otimes\id_W)
\end{equation}
and the triangle axiom
\begin{equation}
\rho_X\otimes\id_Y=(\id_X\otimes\lambda_Y)\circ\alpha_{X,\one,Y}
\end{equation}
for all $X,Y,Z,W\in|\sym{C}|$.
\end{definition}

\begin{definition}
\label{def_lax}
Let $(\sym{C},\otimes,\one,\alpha,\lambda,\rho)$ and
$(\sym{C}^\prime,\otimes^\prime,\one^\prime,\alpha^\prime,\lambda^\prime,\rho^\prime)$
be monoidal categories.
\begin{enumerate}
\item
A \emph{lax monoidal functor}
$(F,F_{X,Y},F_0)\colon\sym{C}\to\sym{C}^\prime$ consists of a
functor $F\colon\sym{C}\to\sym{C}^\prime$, morphisms $F_{X,Y}\colon
FX\otimes^\prime FY\to F(X\otimes Y)$ that are natural in
$X,Y\in|\sym{C}|$, and of a morphism $F_0\colon\one^\prime\to
F\one$, subject to the hexagon axiom
\begin{equation}
F_{X,Y\otimes Z}\circ(\id_{FX}\otimes^\prime F_{Y,Z})\circ\alpha^\prime_{FX,FY,FZ}
= F\alpha_{X,Y,Z}\circ F_{X\otimes Y,Z}\circ(F_{X,Y}\otimes^\prime\id_{FZ})
\end{equation}
and the two squares
\begin{eqnarray}
\lambda^\prime_{FX} &=& F\lambda_X\circ F_{\one,X}\circ(F_0\otimes^\prime\id_{FX}),\\
\rho^\prime_{FX}    &=& F\rho_X\circ F_{X,\one}\circ(\id_{FX}\otimes^\prime F_0)
\end{eqnarray}
for all $X,Y,Z\in|\sym{C}|$.
\item
An \emph{oplax monoidal functor}
$(F,F^{X,Y},F^0)\colon\sym{C}\to\sym{C}^\prime$ consists of a
functor $F\colon\sym{C}\to\sym{C}^\prime$, morphisms $F^{X,Y}\colon
F(X\otimes Y)\to FX\otimes^\prime FY$ that are natural in
$X,Y\in|\sym{C}|$, and of a morphism $F^0\colon
F\one\to\one^\prime$, subject to the hexagon axiom
\begin{equation}
(\id_{FX}\otimes^\prime F^{Y,Z})\circ F^{X,Y\otimes Z}\circ F\alpha_{X,Y,Z}
= \alpha^\prime_{FX,FY,FZ}\circ(F^{X,Y}\otimes^\prime\id_{FZ})\circ F^{X\otimes Y,Z}
\end{equation}
and the two squares
\begin{eqnarray}
F\lambda_X &=& \lambda^\prime_{FX}\circ(F^0\otimes^\prime\id_{FX})\circ F^{\one,X},\\
F\rho_X    &=& \rho^\prime_{FX}\circ(\id_{FX}\otimes^\prime F^0)\circ F^{X,\one}
\end{eqnarray}
for all $X,Y,Z\in|\sym{C}|$.
\item
A \emph{strong monoidal functor}
$(F,F_{X,Y},F_0)\colon\sym{C}\to\sym{C}^\prime$ is a lax monoidal
functor such that $F_0$ and all $F_{X,Y}$, $X,Y\in|\sym{C}|$, are
isomorphisms.
\item
A \emph{strict monoidal functor} $(F,F_{X,Y},F_0)$ is a strong monoidal
functor for which $F_0$ and all $F_{X,Y}$, $X,Y\in|\sym{C}|$, are identity
morphisms.
\end{enumerate}
\end{definition}

\begin{definition}
Let $(F,F_{X,Y},F_0)\colon\sym{C}\to\sym{C}^\prime$ and
$(G,G_{X,Y},G_0)\colon\sym{C}\to\sym{C}^\prime$ be lax monoidal
functors between monoidal categories $\sym{C}$ and $\sym{C}^\prime$. A
\emph{monoidal natural transformation} $\eta\colon F\Rightarrow G$ is
a natural transformation such that
\begin{equation}
\eta_{X\otimes Y}\circ F_{X,Y} = G_{X,Y}\circ(\eta_X\otimes^\prime\eta_Y)
\end{equation}
for all $X,Y\in\sym{C}$.
\end{definition}

There is a similar notion of monoidal natural transformation if the functors
are oplax rather than lax monoidal. Compositions of [lax, oplax, strong]
monoidal functors are again [lax, oplax, strong] monoidal. The following
result is well known, but quite laborious to verify.

\begin{proposition}
\label{prop_monadjunction}
Let $\sym{C}$ and $\sym{C}^\prime$ be monoidal categories and $F\dashv
G\colon\sym{C}^\prime\to\sym{C}$ be an adjunction with unit
$\eta\colon 1_{\sym{C}}\Rightarrow G\circ F$ and counit
$\epsilon\colon F\circ G\Rightarrow 1_{\sym{C}^\prime}$.
\begin{enumerate}
\item
If $F$ has an oplax monoidal structure $(F,F^{C_1,C_2},F^0)$, then
$G$ has a lax monoidal structure $(G,G_{D_1,D_2},G_0)$ as follows,
\begin{eqnarray}
G_{D_1,D_2} &=& G(\epsilon_{D_1}\otimes\epsilon_{D_2})
\circ G(F^{G(D_1),G(D_2)})\circ\eta_{G(D_1)\otimes G(D_2)},\\
G_0 &=& G(F^0)\circ\eta_\one.
\end{eqnarray}
\item
If $F$ is strong monoidal, then both $\eta$ and $\epsilon$ are monoidal
natural transformations.
\item
If $F$ is strong monoidal and the adjunction is an equivalence, then
$G$ is strong monoidal.
\end{enumerate}
\end{proposition}

By an equivalence of monoidal categories, we mean an equivalence of
categories such that one of the functors is strong monoidal and write
$\sym{C}\simeq_\otimes\sym{D}$ in this case.

\subsection{Duality}
\label{app_duals}

\begin{definition}
Let $(\sym{C},\otimes,\one,\alpha,\lambda,\rho)$ be a monoidal
category.
\begin{enumerate}
\item
A \emph{left-dual} $(X^\ast,\ev_X,\coev_X)$ of an object
$X\in|\sym{C}|$ consists of an object $X^\ast\in|\sym{C}|$ and
morphisms $\ev_X\colon X^\ast\otimes X\to\one$
(\emph{left-evaluation}) and $\coev_X\colon\one\to X\otimes X^\ast$
(\emph{left-coevaluation}) that satisfy the triangle identities
\begin{eqnarray}
\label{eq_triangle1}
\rho_X\circ(\id_X\otimes\ev_X)\circ\alpha_{X,X^\ast,X}
\circ(\coev_X\otimes\id_X)\circ\lambda_X^{-1} &=& \id_X,\\
\label{eq_triangle2}
\lambda_{X^\ast}\circ(\ev_X\otimes\id_{X^\ast})\circ\alpha^{-1}_{X^\ast,X,X^\ast}
\circ(\id_{X^\ast}\otimes\coev_X)\circ\rho_{X^\ast}^{-1} &=& \id_{X^\ast}.
\end{eqnarray}
\item
A \emph{right-dual} $(\bar X,\bar\ev_X,\bar\coev_X)$ of
$X\in|\sym{C}|$ consists of an object $\bar X\in|\sym{C}|$ and
morphisms $\bar\ev_X\colon X\otimes\bar X\to\one$
(\emph{right-evaluation}) and $\bar\coev_X\colon\one\to\bar X\otimes
X$ (\emph{right-coevaluation}) that satisfy the triangle identities
\begin{eqnarray}
\lambda_X\circ(\bar\ev_X\otimes\id_X)\circ\alpha^{-1}_{X,\bar X,X}
\circ(\id_X\otimes\bar\coev_X)\circ\rho_X^{-1} &=& \id_X,\\
\rho_{\bar X}\circ(\id_{\bar X}\otimes\bar\ev_X)\circ\alpha_{\bar X,X,\bar X}
\circ(\bar\coev_X\otimes\id_{\bar X})\circ\lambda^{-1}_{\bar X} &=& \id_{\bar X}.
\end{eqnarray}
\end{enumerate}
\end{definition}

\begin{definition}
Let $(\sym{C},\otimes,\one,\alpha,\lambda,\rho)$ be a monoidal category
and $f\colon X\to Y$ be a morphism of $\sym{C}$.
\begin{enumerate}
\item
If both $X$ and $Y$ have left-duals, the \emph{left-dual} of $f$ is defined as
\begin{equation}
\label{eq_leftdual}
f^\ast := \lambda_{X^\ast}\circ(\ev_Y\otimes\id_{X^\ast})\circ\alpha^{-1}_{Y^\ast,Y,X^\ast}
\circ(\id_{Y^\ast}\otimes(f\otimes\id_{X^\ast}))\circ(\id_{Y^\ast}\otimes\coev_X)\circ\rho^{-1}_{Y^\ast}.
\end{equation}
\item
If both $X$ and $Y$ have right-duals, the \emph{right-dual} of $f$ is defined as
\begin{equation}
\bar f:=\rho_{\bar X}\circ(\id_{\bar X}\otimes\bar\ev_Y)\circ\alpha_{\bar X,Y,\bar Y}
\circ((\id_{\bar X}\otimes f)\otimes\id_{\bar Y})\circ(\bar\coev_X\otimes\id_{\bar Y})\circ\lambda^{-1}_{\bar Y}.
\end{equation}
\end{enumerate}
\end{definition}

\begin{proposition}
\label{prop_dualisic}
Let $(\sym{C},\otimes,\one,\alpha,\lambda,\rho)$ be a monoidal category and
$(X^\ast,\ev_X,\coev_X)$ and $(X^\prime,\ev_X^\prime,\coev_X^\prime)$ both be
left-duals of $X\in|\sym{C}|$. Then there is a natural isomorphism
\begin{equation}
u_X=\lambda_{X^\prime}\circ(\ev_X\otimes\id_{X^\prime})
\circ\alpha^{-1}_{X^\ast,X,X^\prime}\circ(\id_{X^\ast}\otimes\coev^\prime_X)
\circ\rho^{-1}_{X^\ast}\colon X^\ast\to X^\prime
\end{equation}
with inverse
\begin{equation}
u_X^{-1}=\lambda_{X^\ast}\circ(\ev_X^\prime\otimes\id_{X^\ast})
\circ\alpha^{-1}_{X^\prime,X,X^\ast}\circ(\id_{X^\prime}\otimes\coev_X)
\circ\rho^{-1}_{X^\prime}\colon X^\prime\to X^\ast.
\end{equation}
They satisfy $\ev_X = \ev^\prime_X\circ(u_X\otimes\id_X)$ and $\coev_X =
(\id_X\otimes u_X^{-1})\circ\coev^\prime_X$.
\end{proposition}

\begin{definition}
A \emph{[left-, right-]autonomous category} is a monoidal category in
which each object is equipped with a specified [left-, right-]dual. An
\emph{autonomous category} is a monoidal category that is both left-
and right-autonomous.
\end{definition}

\begin{definition}
A \emph{pivotal category}
$(\sym{C},\otimes,\one,\alpha,\lambda,\rho,{(-)}^\ast,\ev,\coev,\tau)$
is a left-autonomous category
$(\sym{C},\otimes,\one,\alpha,\lambda,\rho,{(-)}^\ast,\ev,\coev)$ with
natural isomorphisms $\tau_X\colon X\to {X^\ast}^\ast$ such that
\begin{equation}
\label{eq_pivotal}
(\tau_X)\ast = \tau_{X^\ast}^{-1}.
\end{equation}
for all $X\in|\sym{C}|$.
\end{definition}

\noindent
Note that every pivotal category is also right-autonomous with $\bar
X=X^\ast$ and
\begin{eqnarray}
\bar\ev_X   &=& \ev_{X^\ast}\circ(\tau_X\otimes\id_{X^\ast}),\\
\label{eq_barcoev}
\bar\coev_X &=& (\id_{X^\ast}\otimes\tau_X^{-1})\circ\coev_{X^\ast}
\end{eqnarray}
for all $X\in|\sym{C}|$. In a pivotal category, $\bar f=f^\ast$ for
all morphisms $f\colon X\to Y$.

\begin{definition}
\label{def_traces}
Let
$(\sym{C},\otimes,\one,\alpha,\lambda,\rho,{(-)}^\ast,\ev,\coev,\tau)$
be a pivotal category and $f\colon X\to X$ be a morphism of
$\sym{C}$. We define the \emph{left-trace} of $f$
\begin{equation}
\tr^{(L)}_X(f) = \ev_X\circ(\id_{X^\ast}\otimes f)\circ\bar\coev_X\colon\one\to\one
\end{equation}
and the \emph{right-trace} of $f$
\begin{equation}
\tr^{(R)}_X(f) = \bar\ev_X\circ(f\otimes\id_{X^\ast})\circ\coev_X\colon\one\to\one.
\end{equation}
\end{definition}

Note that in a pivotal category, both left- and right-traces are
cyclic, \ie\ $\tr^{(L)}_X(g\circ f)=\tr^{(L)}_Y(f\circ g)$ for all
$f\colon X\to Y$ and $g\colon Y\to X$ and similarly for the
right-trace.

\begin{definition}
A \emph{spherical category}
$(\sym{C},\otimes,\one,\alpha,\lambda,\rho,{(-)}^\ast,\ev,\coev,\tau)$
is a pivotal category in which
\begin{equation}
\tr^{(L)}_X(f) = \tr^{(R)}_X(f)
\end{equation}
for all morphisms $f\colon X\to X$ in $\sym{C}$. In this case, the above
expression is just called the \emph{trace} of $f$ and denoted by $\tr_X(f)$,
and
\begin{equation}
\dim(X) = \tr_X(\id_X)
\end{equation}
is called the \emph{dimension} of $X$.
\end{definition}

Note that in a spherical category, $\tr_X(f)=\tr_{X^\ast}(f^\ast)$ for every
morphism $f\colon X\to X$ and thus $\dim(X)=\dim(X^\ast)$.

\begin{example}
\label{ex_ribboncat}
Every ribbon category (see, for example~\cite[Appendix A.3]{Pf07}) is spherical.
\end{example}

\begin{proof}
Every ribbon category is a pivotal category because
of~\cite[eq.~(A.24)]{Pf07}. The spherical property follows
from~\cite[eq.~(A.25) and eq.~(A.26)]{Pf07} and from the naturality of the braiding.
\end{proof}

\subsection{Abelian and semisimple categories}

A category $\sym{C}$ is called \emph{$\mathbf{Ab}$-enriched} if it is enriched
in the category $\mathbf{Ab}$ of abelian groups, \ie\ if $\Hom(X,Y)$ is an
abelian group for all objects $X,Y\in|\sym{C}|$ and if the composition of
morphisms is $\Z$-bilinear. If $k$ is a commutative ring, a category $\sym{C}$
is called \emph{$k$-linear} if it is enriched in ${}_k\sym{M}$, the category
of $k$-modules, \ie\ if $\Hom(X,Y)$ is a $k$-module for all $X,Y\in|\sym{C}|$
and if the composition of morphisms is $k$-bilinear.

A functor $F\colon\sym{C}\to\sym{C}^\prime$ between [$\mathbf{Ab}$-enriched,
$k$-linear] categories is called [\emph{additive}, $k$-\emph{linear}] if it
induces homomorphisms of [additive groups, $k$-modules]
\begin{equation}
\Hom(X,Y)\to\Hom(FX,FY)
\end{equation}
for all $X,Y\in|\sym{C}|$.

An \emph{additive category} is an $\mathbf{Ab}$-enriched category that has a
terminal object and all binary products. A \emph{preabelian category} is an
$\mathbf{Ab}$-enriched category that has all finite limits. An \emph{abelian
category} is a preabelian category in which every monomorphism is a kernel
and in which every epimorphism is a cokernel. A functor
$F\colon\sym{C}\to\sym{C}^\prime$ between preabelian categories is called
\emph{exact} if it preserves all finite limits. An equivalence of
[$\mathbf{Ab}$-enriched, $k$-linear] categories is an equivalence of
categories, one functor of which is [additive, $k$-linear].

\begin{definition}
\label{def_semisimple}
Let $\sym{C}$ be a $k$-linear category and $k$ be a commutative ring.
\begin{myenumerate}
\item
An object $X\in|\sym{C}|$ is called \emph{simple} if $\End(X)\cong
k$ are isomorphic as $k$-modules.
\item
An object $X\in|\sym{C}|$ is called \emph{null} if
$\End(X)\cong\{0\}$.
\item
The category $\sym{C}$ is called \emph{semisimple} if there exists a
family ${\{V_j\}}_{j\in I}$ of objects $V_j\in|\sym{C}|$, $I$ some
index set, such that
\begin{myenumerate}
\item
$V_j$ is simple for all $j\in I$.
\item
$\Hom(V_j,V_\ell)=\{0\}$ for all $j,\ell\in I$ for which $j\neq\ell$.
\item
For each object $X\in|\sym{C}|$, there is a finite sequence
$j_1^{(X)},\ldots,j_{n^X}^{(X)}\in I$, $n^X\in\N_0$, and morphisms
$\imath_\ell^{(X)}\colon V_{j_\ell}\to X$ and $\pi_\ell^{(X)}\colon
X\to V_{j_\ell}$ such that
\begin{equation}
\label{eq_dominance}
\id_X = \sum_{\ell=1}^{n^X}\imath^X_\ell\circ\pi^X_\ell.
\end{equation}
and
\begin{equation}
\label{eq_semisimpledualbases}
\pi^X_\ell\circ\imath^X_m = \left\{
\begin{matrix}
\id_{V_{j^X_\ell}},&\mbox{if}\quad \ell=m,\\
0,                 &\mbox{else}
\end{matrix}
\right.
\end{equation}
\end{myenumerate}
\item
The category is called \emph{finitely semisimple} (also
\emph{Artinian semisimple}) if it is semisimple with a finite index
set $I$ in condition~(3).
\end{myenumerate}
\end{definition}

\begin{proposition}[see {\cite[Lemma II.4.2.2]{Tu94}}]
\label{prop_semisimplehom}
Let $\sym{C}$ be a $k$-linear category and $k$ a commutative ring. If
$\sym{C}$ is [finitely] semisimple, then there is a [finite] set
$\sym{J}\subseteq|\sym{C}|$ of non-null objects such
that
\begin{eqnarray}
\label{eq_semisimplehom}
\Phi\colon\bigoplus_{J\in\sym{J}}\Hom(X,J)\otimes\Hom(J,Y) &\to& \Hom(X,Y),\nn\\
f\otimes g &\mapsto& g\circ f,
\end{eqnarray}
is an isomorphism for all $X,Y\in|\sym{C}|$.
\end{proposition}

If $\sym{C}$ is a semisimple $k$-linear category, then by~\cite[Proposition
II.4.2.1]{Tu94}, $\Hom(X,Y)$ is a finitely generated projective $k$-module for
all $X,Y\in|\sym{C}|$. If $\sym{C}$ is a semisimple $k$-linear category with
family ${\{V_j\}}_{j\in I}$ of simple objects as in
Definition~\ref{def_semisimple}(3) and $k$ is a field, then for each simple
$X\in|\sym{C}|$, there is some $j\in I$ such that $X\cong V_j$.

\subsection{Additive and non-degenerate spherical categories}

A monoidal category  is called
[\emph{$\mathbf{Ab}$-enriched}, $k$-\emph{linear}] if $\sym{C}$ is
[$\mathbf{Ab}$-enriched, $k$-linear] and if the tensor product of morphisms is
[$\Z$-bilinear, $k$-bilinear].

In a monoidal category $(\sym{C},\otimes,\one,\alpha,\lambda,\rho)$, the set
$k:=\End(\one)$ forms a commutative monoid with respect to composition. If
$\sym{C}$ is $\mathbf{Ab}$-enriched as a monoidal category, then $k$ is a
unital commutative ring and $\sym{C}$ is $k$-linear as an ordinary category,
but not necessarily as a monoidal category.

In a $k$-linear pivotal category, the left- and right-traces
\begin{equation}
\tr^{(L)}_X\colon\End(X)\to k\qquad\mbox{and}\qquad
\tr^{(R)}_X\colon\End(X)\to k
\end{equation}
are $k$-linear for all $X\in|\sym{C}|$.

If we work with semisimple pivotal categories, we also require the set of
representatives of the simple objects to contain the monoidal unit and to be
closed under duality as follows.

\begin{definition}
\label{def_pivotalss}
A $k$-linear pivotal category
$(\sym{C},\otimes,\one,\alpha,\lambda,\rho,{(-)}^\ast,\ev,\coev,\tau)$,
$k=\End(\one)$, is called [\emph{finitely}] \emph{semisimple} if the
underlying $k$-linear category is [finitely] semisimple and the family
${\{V_j\}}_j$ of Definition~\ref{def_semisimple}(3) satisfies the following
conditions:
\begin{myenumerate}
\item
There is an element $0\in I$ such that $V_0\cong\one$.
\item
For each $j\in I$, there is some $j^\ast\in I$ such that $V_{j^\ast}\cong V_j^\ast$.
\end{myenumerate}
\end{definition}

\begin{definition}
\label{def_nondegenerate}
A $k$-linear spherical category
$(\sym{C},\otimes,\one,\alpha,\lambda,\rho,{(-)}^\ast,\ev,\coev,\tau)$,
$k=\End(\one)$, is called \emph{non-degenerate} if the $k$-bilinear maps
\begin{equation}
\Hom(Y,X)\otimes\Hom(X,Y)\to k,\quad
f\otimes g\mapsto \tr_X(f\circ g)
\end{equation}
are non-degenerate for all objects $X,Y\in|\sym{C}|$, \ie\ if $\tr_X(f\circ
g)=0$ for all $g\colon X\to Y$ implies $f=0$.
\end{definition}

In a $k$-linear spherical category, the trace is multiplicative with respect
to the tensor product, \ie\ $\tr_{X_1\otimes X_2}(h_1\otimes
h_2)=\tr_{X_1}(h_1)\cdot \tr_{X_2}(h_2)$ for all $h_1\colon X_1\to X_1$ and
$h_2\colon X_2\to X_2$.

\begin{proposition}[analogous to {\cite[Lemma II.4.2.3]{Tu94}}]
\label{prop_nondegenerate}
Every semisimple $k$-linear spherical category is non-degenerate.
\end{proposition}

\begin{proposition}[analogous to {\cite[Lemma II.4.2.4]{Tu94}}]
Let $\sym{C}$ be a semisimple $k$-linear spherical category with family
${\{V_j\}}_{j\in I}$ as in Definition~\ref{def_semisimple}(3). Then for all
$j\in I$, $\dim V_j$ is invertible in $k$.
\end{proposition}

\begin{proposition}[see {\cite[Proposition A.29]{Pf07}}]
Let $\sym{C}$ be a $k$-linear spherical category and $k=\End(\one)$ be a
field. If $\sym{C}$ satisfies all conditions of a finitely semisimple category
of Definition~\ref{def_semisimple}(3) except maybe
for~\eqref{eq_semisimpledualbases}, then the $\imath_\ell^{(X)}$ and
$\pi_\ell^{(X)}$ can be chosen in such a way
that~\eqref{eq_semisimpledualbases} holds as well.
\end{proposition}

\subsubsection*{Acknowledgements}

The author would like to thank Gabriella B{\"o}hm, Shahn Majid and
Korn{\'e}l Szlach{\'a}nyi for valuable discussions and everyone at
RMKI Budapest for their hospitality.

\newenvironment{hpabstract}{%
\renewcommand{\baselinestretch}{0.2}
\begin{footnotesize}%
}{\end{footnotesize}}%
\newcommand{\hpeprint}[2]{%
\href{http://www.arxiv.org/abs/#1}{\texttt{arxiv:#1#2}}}%
\newcommand{\hpspires}[1]{%
\href{http://www.slac.stanford.edu/spires/find/hep/www?#1}{SPIRES Link}}%
\newcommand{\hpmathsci}[1]{%
\href{http://www.ams.org/mathscinet-getitem?mr=#1}{\texttt{MR #1}}}%
\newcommand{\hpdoi}[1]{%
\href{http://dx.doi.org/#1}{\texttt{DOI #1}}}%
\newcommand{\hpjournal}[2]{%
\href{http://dx.doi.org/#2}{\textsl{#1\/}}}%

\end{document}